\documentclass[a4paper,11pt,twoside,english]{article}
\usepackage{babel}
\usepackage{latexsym,amsfonts}
\usepackage{amsthm}
\usepackage{amsmath}
\usepackage{mathtools}
\usepackage{paralist}
\usepackage{ifthen}
\usepackage{moredefs}
\usepackage{todonotes}
\usepackage{tikz}
\usepackage{booktabs}
\usepackage{float}

\usepackage{caption}
\providecommand{\CAP}[1][-30pt]{
  \captionsetup{font=small,name=Fig.,justification=raggedleft,singlelinecheck=false,skip=#1}
  \caption{}
  }

\usepackage{lastpage}
\usepackage{fancyhdr}
\fancypagestyle{plain}{

  \fancyhead[L,C]{}
  \fancyhead[R]{\today}
  \fancyfoot[R,L]{}
  \fancyfoot[C]{\thepage \ of \pageref{LastPage}}
  }

\title{Absolute sums of Banach spaces and some geometric properties related to rotundity and smoothness}
\author{Jan-David Hardtke}
\date{}

\DeclareMathOperator{\lin}{span}

\DeclareMathOperator{\ran}{ran}

\providecommand{\shortversion}[1]{#1}
\providecommand{\longversion}[2]{#2 {\rm(#1 in short)}}

\providecommand{\acronym}[2]{
  \newboolean{acronym#1}
  \setboolean{acronym#1}{true}
  \expandafter\providecommand\expandafter{\csname acronym#1\endcsname}{#2}
  }

\providecommand{\ac}[1]{\ifthenelse{\boolean{acronym#1}}
{\longversion{#1}{\csname acronym#1\endcsname}\Global\ToggleBoolean{acronym#1}}
{\shortversion{#1}}}

\providecommand{\ifif}{iff }
\providecommand{\sm}{\setminus}
\providecommand{\ssq}{\subseteq}
\providecommand{\id}{\ensuremath{\mathrm{id}}}
\providecommand{\N}{\ensuremath{\mathbb{N}}}

\providecommand{\R}{\ensuremath{\mathbb{R}}}

\providecommand{\T}{\ensuremath{\mathbb{T}}}

\providecommand{\eps}{\ensuremath{\varepsilon}}

\providecommand{\Norm}[1]{\ensuremath{\left|\hspace{-1.2pt}\left|\hspace{-1.2pt}\left|#1\right|\hspace{-1.2pt}\right|\hspace{-1.2pt}\right|}}
\providecommand{\Todo}[2][]{\todo[noline,backgroundcolor=white,size=\small,#1]{#2}}

\providecommand{\keywords}[1]{
  {\let\thefootnote=\relax
  \footnote{{\em Keywords}: #1}}
  \addtocounter{footnote}{-1}
  }

\providecommand{\AMS}[1]{
  {\let\thefootnote=\relax
  \footnote{{\em AMS Subject Classification} (2010): #1}}
  \addtocounter{footnote}{-1}
  }

\providecommand{\address}{
  {\sc \noindent Department of Mathematics \\
  Freie Universit\"at Berlin \\
  Arnimallee 6, 14195 berlin \\
  Germany \\}
  }

\DeclarePairedDelimiter{\set}{\lbrace}{\rbrace}
\DeclarePairedDelimiter{\paren}{\lparen}{\rparen}

\DeclarePairedDelimiter{\abs}{\lvert}{\rvert}
\DeclarePairedDelimiter{\norm}{\lVert}{\rVert}

\theoremstyle{definition}
\newtheorem{definition}{Definition}[section]
\newtheorem*{definition*}{Definition}
\newtheorem{example}[definition]{Example}
\newtheorem*{example*}{Example}

\newtheorem*{remark*}{Remark}

\newtheorem*{problem*}{Problem}

\newtheorem*{question*}{Question}

\newtheorem*{conjecture*}{Conjecture}
\theoremstyle{remark}

\newtheorem*{claim*}{Claim}

\newtheorem*{fact*}{Fact}
\theoremstyle{plain}
\newtheorem{lemma}[definition]{Lemma}
\newtheorem*{lemma*}{Lemma}
\newtheorem{proposition}[definition]{Proposition}
\newtheorem*{proposition*}{Proposition}
\newtheorem{theorem}[definition]{Theorem}
\newtheorem*{theorem*}{Theorem}
\newtheorem{corollary}[definition]{Corollary}
\newtheorem*{corolary*}{Corollary}

\newenvironment{Proof}[1][\proofname]{\begin{proof}[#1] \setlength{\parindent}{0pt}}{\end{proof}}
\newenvironment{Abstract}{\centering\begin{minipage}{0.8\textwidth} \noindent \small {\sc Abstract.}}{\end{minipage}\par}

\usepackage{color}
\definecolor{darkgreen}{rgb}{0,0.5,0}

\numberwithin{equation}{section}
\newtagform{colored}[\color{blue}]{\color{blue}(}{\color{blue})}
\usetagform{colored}

\usepackage[colorlinks,linkcolor=blue,citecolor=red,urlcolor=darkgreen]{hyperref}
\providecommand{\email}{{\it E-mail address:} \href{mailto:hardtke@math.fu-berlin.de}{\tt hardtke@math.fu-berlin.de}}
\providecommand{\mr}[1]{\href{http://www.ams.org/mathscinet-getitem?mr=#1}{MR#1}}

\usepackage{amsrefs}

\begin{document}

\maketitle

\begin{Abstract}
\noindent We study the notions of acs, luacs and uacs Banach spaces which were introduced in \cite{kadets} and 
form common generalisations of the usual rotundity and smoothness properties of Banach spaces. In particular, 
we are interested in (mainly infinite) absolute sums of such spaces. We also introduce some new classes of spaces 
that lie inbetween those of acs and uacs spaces and study their behaviour under the formation of absolute sums as well.
\end{Abstract}
\AMS{46B20}

\section{Introduction}\label{sec:intro}
\acronym{UR}{uniformly rotund}
\acronym{LUR}{locally uniformly rotund}
\acronym{WLUR}{weakly locally uniformly rotund}
\acronym{WUR}{weakly uniformly rotund}
\acronym{URED}{uniformly rotund in every direction}
\acronym{R}{rotund}
\acronym{S}{smooth}
\acronym{FS}{Fr\'echet-smooth}
\acronym{UG}{uniformly G\^ateaux-differentiable}
\acronym{US}{uniformly smooth}
\acronym{acs}{alternatively convex or smooth}
\acronym{luacs}{locally uniformly alternatively convex or smooth}
\acronym{uacs}{uniformly alternatively convex or smooth}
\acronym{sluacs}{strongly locally uniformly alternatively convex or smooth}
\acronym{wuacs}{weakly uniformly alternatively convex or smooth}
\acronym{uacsed}{uniformly alternatively convex or smooth in every direction}
\acronym{MLUR}{midpoint locally uniformly rotund}
\acronym{WMLUR}{weakly midpoint locally uniformly rotund}
\acronym{mluacs}{midpoint locally uniformly alternatively convex or smooth}
\acronym{msluacs}{midpoint strongly locally uniformly alternatively convex or smooth}
\acronym{aDP}{alternative Daugavet property}
\acronym{aDE}{alternative Daugavet equation}

First let us fix some notation. Where not otherwise stated, $X$ denotes a real Banach space, $X^*$ its 
dual, $B_X$ its unit ball and $S_X$ its unit sphere.\par
Since we will deal with various generalisations of rotundity and smoothness properties for Banach spaces, 
we start by recalling the most important of these notions.
\begin{definition}\label{def:rotundity}
A Banach space $X$ is called
  \begin{enumerate}[(i)]
    \item {\em\ac{R}} if for any two elements $x,y\in S_X$ the equality $\norm*{x+y}=2$ implies $x=y$,
    \item {\em\ac{LUR}} if for every $x\in S_X$ the implication
    \begin{equation*}
    \norm*{x_n+x}\to 2 \ \Rightarrow \ \norm*{x_n-x}\to 0
    \end{equation*}
    holds for every sequence $(x_n)_{n\in \N}$ in $S_X$,
    \item {\em\ac{WLUR}} if for every $x\in S_X$ and every sequence $(x_n)_{n\in \N}$ in $S_X$ we have
    \begin{equation*}
    \norm*{x_n+x}\to 2 \ \Rightarrow \ x_n\to x \ \mathrm{weakly},
    \end{equation*}
    \item {\em\ac{UR}} if for any two sequences $(x_n)_{n\in \N}$ and $(y_n)_{n\in \N}$ in $S_X$ the implication 
    \begin{equation*}
    \norm*{x_n+y_n}\to 2 \ \Rightarrow \ \norm*{x_n-y_n}\to 0
    \end{equation*}
    holds,
    \item {\em\ac{WUR}} if for any two sequences $(x_n)_{n\in \N}$ and $(y_n)_{n\in \N}$ the following implication holds
    \begin{equation*}
    \norm*{x_n+y_n}\to 2 \ \Rightarrow \ x_n-y_n\to 0 \ \mathrm{weakly}.
    \end{equation*}
  \end{enumerate}
\end{definition}
\noindent The obvious implications between these notions are summarised in the chart below and no other implications 
are valid, as is shown by the examples in \cite{smith}.

\begin{figure}[H]
\begin{center}
  \begin{tikzpicture}
  \node (UR) at (-2,0) {UR};
  \node (WUR) at (0,1) {WUR};
  \node (LUR) at (0,-1) {LUR};
  \node (WLUR) at (2,0) {WLUR};
  \node (R) at (4,0) {R};
  \draw[->] (UR)--(LUR);
  \draw[->] (UR)--(WUR);
  \draw[->] (LUR)--(WLUR);
  \draw[->] (WUR)--(WLUR);
  \draw[->] (WLUR)--(R);
  \end{tikzpicture}
\end{center}
\CAP\label{fig:1}
\end{figure}

Note that, by standard normalisation arguments, $X$ is \ac{UR} \ifif for all bounded sequences $(x_n)_{n\in \N}$ 
and $(y_n)_{n\in \N}$ in $X$ which fulfil the conditions $\norm*{x_n+y_n}-\norm*{x_n}-\norm*{y_n}\to 0$ and 
$\norm*{x_n}-\norm*{y_n}\to 0$ we have that $\norm*{x_n-y_n}\to 0$ and further that the two conditions 
$\norm*{x_n+y_n}-\norm*{x_n}-\norm*{y_n}\to 0$ and $\norm*{x_n}-\norm*{y_n}\to 0$ can be replaced by the single 
equivalent condition $2\norm*{x_n}^2+2\norm*{y_n}^2-\norm*{x_n+y_n}^2\to 0$. Similar remarks apply to the definitions 
of \ac{LUR}, \ac{WUR} and \ac{WLUR} spaces. Also, for a finite-dimensional space $X$ all the above notions coincide 
(by compactness of the unit ball).\par
Recall also that the modulus of convexity of the space $X$ is defined by 
$\delta_X(\eps)=\inf\set*{1-1/2\norm*{x+y}:x,y\in B_X \ \mathrm{and} \ \norm*{x-y}\geq\eps}$ for every $\eps$ in the 
interval $]0,2]$. Then $X$ is \ac{UR} \ifif $\delta_X(\eps)>0$ for all $0<\eps\leq2$.\par
Concerning notions of smoothness, the space $X$ is called {\em\ac{S}} if its norm is G\^ateaux-differentiable at every non-zero 
point (equivalently at every point of $S_X$), which is the case \ifif for every $x\in S_X$ there is a unique functional 
$x^*\in S_{X^*}$ with $x^*(x)=1$ (cf. \cite{fabian}*{Lemma 8.4 (ii)}). $X$ is called {\em\ac{FS}} if the norm is 
Fr\'echt-differentiable at every non-zero point. Finally, $X$ is called {\em\ac{US}} if $\lim_{\tau \to 0}\rho_X(\tau)/\tau=0$,
where $\rho_X$ denotes the modulus of smoothness of $X$ defined by 
$\rho_X(\tau)=\sup\set*{1/2(\norm*{x+\tau y}+\norm*{x-\tau y}-2):x,y\in S_X}$ for every $\tau>0$.\par 
Obviously, \ac{FS} implies \ac{S} and from \cite{fabian}*{Fact 9.7} it follows that \ac{US} implies \ac{FS}. It is also 
well known that $X$ is \ac{US} \ifif $X^*$ is \ac{UR} and $X$ is  \ac{UR} \ifif $X^*$ is \ac{US} (cf. \cite{fabian}*{Theorem 9.10}).\par
There is yet another notion of smoothness, namely the norm of the space $X$ is said to be {\em\ac{UG}} if for each $y\in S_X$
the limit $\lim_{\tau\to 0}\paren*{\norm{x+\tau y}-1}/\tau$ exists uniformly in $x\in S_X$. The property \ac{UG} lies between
\ac{US} and \ac{S}. It is known (cf. \cite{deville}*{Theorem II.6.7}) that $X^*$ is \ac{UG} \ifif $X$ is \ac{WUR} and $X$ is \ac{UG} 
\ifif $X^*$ is $\mathrm{WUR}^*$ (which means that $X^*$ fulfils the definition of \ac{WUR} with weak- replaced by weak*-convergence).\par
In \cite{kadets} the following notions were introduced (in connection with the so called Anti-Daugavet property).
\begin{definition}\label{def:acs-luacs-uacs}
A Banach space $X$ is called
  \begin{enumerate}[(i)]
  \item {\em\ac{acs}} if for every $x,y\in S_X$ with $\norm*{x+y}=2$ and every $x^*\in S_{X^*}$ with $x^*(x)=1$ we have $x^*(y)=1$ as well,
  \item {\em\ac{luacs}} if for every $x\in S_X$, every sequence $(x_n)_{n\in \N}$ in $S_X$ and every functional $x^*\in S_{X^*}$ we have
  \begin{equation*}
  \norm*{x_n+x}\to 2 \ \mathrm{and} \ x^*(x_n)\to 1 \ \Rightarrow \ x^*(x)=1,
  \end{equation*}
  \item {\em\ac{uacs}} if for all sequences $(x_n)_{n\in \N}$, $(y_n)_{n\in \N}$ in $S_X$ and $(x_n^*)_{n\in \N}$ in $S_{X^*}$ we have
  \begin{equation*}
  \norm*{x_n+y_n}\to 2 \ \mathrm{and} \ x_n^*(x_n)\to 1 \ \Rightarrow \ x_n^*(y_n)\to 1.
  \end{equation*}
  \end{enumerate}
\end{definition}
\noindent Clearly, \ac{R} and \ac{S} both imply \ac{acs}, \ac{WLUR} implies \ac{luacs} and \ac{UR} and \ac{US} both
imply \ac{uacs}. Again by standard normalisation arguments one can easily check that $X$ is \ac{uacs} \ifif for all 
bounded sequences $(x_n)_{n\in \N}$, $(y_n)_{n\in \N}$ in $X$ and $(x_n^*)_{n\in \N}$ in $X^*$ with 
$x_n^*(x_n)-\norm*{x_n^*}\norm*{x_n}\to 0$, $\norm*{x_n+y_n}-\norm*{x_n}-\norm*{y_n}\to 0$ and $\norm*{x_n}-\norm*{y_n}\to 0$ 
(or equivalently $2\norm*{x_n}^2+2\norm*{y_n}^2-\norm*{x_n+y_n}^2\to 0$) we also have $x_n^*(y_n)-\norm*{x_n^*}\norm*{y_n}\to 0$ 
and a similar characterisation holds for \ac{luacs} spaces.\par
Note also the following reformulation of the definition of \ac{acs} spaces, which was observed in \cite{kadets}:
A Banach space $X$ is \ac{acs} \ifif whenever $x,y\in S_X$ such that $\norm{x+y}=2$ then the norm of $\lin\set{x,y}$ is
G\^ateaux-differentiable at $x$ and $y$.\par
Finally let us note that, again by compactness, in the case $\dim{X}<\infty$ the notions of \ac{acs}, \ac{luacs} and \ac{uacs} 
spaces coincide.\par 
Recall that a Banach space $X$ is said to be uniformly non-square if there is some $\delta>0$ such that for all $x,y\in B_X$
we have $\norm*{x+y}\leq2(1-\delta)$ or $\norm*{x-y}\leq2(1-\delta)$. It is easily seen that \ac{uacs} spaces are uniformly 
non-square and hence by a well-known theorem of James (cf. \cite{beauzamy}*{p.261}) they are superreflexive, as was observed
in \cite{kadets}*{Lemma 4.4}.\par
Actually, to prove the superreflexivity of \ac{uacs} spaces it is not necessary to employ the rather deep theorem of James,
as we will see in the next section.\par
In \cite{sirotkin} it is shown by G. Sirotkin that for every $1<p<\infty$ and every measure space $(\Omega,\Sigma,\mu)$ the 
Lebesgue-Bochner space $L^p(\Omega,\Sigma,\mu;X)$ is \ac{uacs} (resp. \ac{luacs}, resp. \ac{acs}) whenever $X$ is an \ac{uacs}
(resp. \ac{luacs}, resp. \ac{acs}) Banach space. To get this result, Sirotkin first proves  the following characterisation of 
\ac{uacs} spaces.
\begin{proposition}[Sirotkin, cf. \cite{sirotkin}]\label{prop:char-uacs}
A Banach space $X$ is \ac{uacs} \ifif for any two sequences $(x_n)_{n\in \N}$ and $(y_n)_{n\in \N}$ in $S_X$ and every 
sequence $(x_n^*)_{n\in \N}$ in $S_{X^*}$ we have
  \begin{equation*}
  \norm*{x_n+y_n}\to 2 \ \mathrm{and} \ x_n^*(x_n)=1 \ \forall n\in \N \ \Rightarrow \ x_n^*(y_n)\to 1.
  \end{equation*}
\end{proposition}
\noindent Instead of repeating the proof from \cite{sirotkin} here, we shall give a slightly different proof below 
(see Proposition \ref{prop:char-uacs-sluacs}), which---unlike Sirotkin's proof---does not use any reflexivity arguments
(but see also the proof of Lemma \ref{lemma:delta deltatilde}).\par 
Now with this characterisation we can define a kind of `\ac{uacs}-modulus' of a given Banach space.
\begin{definition}\label{def:uacs-mod}
For a Banach space $X$ we define
  \begin{align*}
  & D_X(\eps)=\set*{(x,y)\in S_X\times S_X:\exists x^*\in S_{X^*} \ x^*(x)=1 \ \mathrm{and} \ x^*(y)\leq 1-\eps} \\
  & \mathrm{and} \ \delta_{\mathrm{uacs}}^X(\eps)=\inf\set*{1-\norm*{\frac{x+y}{2}}:(x,y)\in D_X(\eps)} \ \forall \eps\in ]0,2].
  \end{align*}
\end{definition}
\noindent Then by Proposition \ref{prop:char-uacs} $X$ is \ac{uacs} \ifif $\delta_{\mathrm{uacs}}^X(\eps)>0$ for every 
$\eps\in ]0,2]$ and we clearly have $\delta_X(\eps)\leq\delta_{\mathrm{uacs}}^X(\eps)$ for each $\eps\in ]0,2]$. For the 
connection to the modulus of smoothness see Lemma \ref{lemma:delta rho}.\par
The characterisation of \ac{uacs} spaces given above coincides with the notion of $U$-spaces 
introduced by Lau in \cite{lau} and our modulus $\delta_{\mathrm{uacs}}^X$ is the same as the 
modulus of $u$-convexity from \cite{gao1}. Also, the notion of $u$-spaces which was introduced in \cite{dhompongsa}
coincides with the notion of \ac{acs} spaces. The interested reader may also have a look at \cite{dutta}, where two notions
of local $U$-convexity are introduced and studied quantitatively.\par
The $U$-spaces ($=$ \ac{uacs} spaces) are of particular interest, because they possess normal structure (cf. \cite{gao2}*{Theorem 3.2}
or \cite{sirotkin}*{Theorem 3.1}) and hence (since they are also reflexive) they enjoy the fixed point property.\footnote{The 
reader is referred to \cite{goebel}*{Section 2} for definitions and background.}\par
It seems natural to introduce two more notions related to \ac{uacs} spaces, namely the following.
\begin{definition}\label{def:sluacs-wuacs}
A Banach space $X$ is called 
  \begin{enumerate}[(i)]
  \item {\em\ac{sluacs}} if for every $x\in S_X$ and all sequences $(x_n)_{n\in \N}$ in $S_X$ and $(x_n^*)_{n\in \N}$ in $S_{X^*}$ we have
  \begin{equation*}
  \norm*{x_n+x}\to 2 \ \mathrm{and} \ x_n^*(x_n)\to 1 \ \Rightarrow \ x_n^*(x)\to 1,
  \end{equation*}
  \item {\em\ac{wuacs}} if for any two sequences $(x_n)_{n\in \N}$, $(y_n)_{n\in \N}$ in $S_X$ and every functional $x^*\in S_{X^*}$ we have
  \begin{equation*}
  \norm*{x_n+y_n}\to 2 \ \mathrm{and} \ x^*(x_n)\to 1 \ \Rightarrow \ x^*(y_n)\to 1.
  \end{equation*}
  \end{enumerate}
\end{definition}
\noindent With these definitions we get the following implication chart.

\begin{figure}[H]
\begin{center}
  \begin {tikzpicture}
  \node (uacs) at (-2,0) {uacs};
  \node (wuacs) at (0,1) {wuacs};
  \node (sluacs) at (0,-1) {sluacs};
  \node (luacs) at (2,0) {luacs};
  \node (acs) at (4,0) {acs};
  \draw[->] (uacs)--(sluacs);
  \draw[->] (uacs)--(wuacs);
  \draw[->] (sluacs)--(luacs);
  \draw[->] (wuacs)--(luacs);
  \draw[->] (luacs)--(acs);
  \end{tikzpicture}
\end{center}
\CAP\label{fig:2}
\end{figure}

\noindent Including the rotundity properties finally leaves us with the diagram below.

\begin{figure}[H]
\begin{center}
\begin{tikzpicture}
  \node (UR) at (-2,0) {UR};
  \node (WUR) at (0,1) {WUR};
  \node (LUR) at (0,-1) {LUR};
  \node (WLUR) at (2,0) {WLUR};
  \node (R) at (4,0) {R};
  \node (uacs) at (-2,-1) {uacs};
  \node (wuacs) at (0,0) {wuacs};
  \node (sluacs) at (0,-2) {sluacs};
  \node (luacs) at (2,-1) {luacs};
  \node (acs) at (4,-1) {acs};
  \draw[->] (UR)--(LUR);
  \draw[->] (UR)--(WUR);
  \draw[->] (LUR)--(WLUR);
  \draw[->] (WUR)--(WLUR);
  \draw[->] (WLUR)--(R);
  \draw[->] (uacs)--(sluacs);
  \draw[->] (uacs)--(wuacs);
  \draw[->] (sluacs)--(luacs);
  \draw[->] (wuacs)--(luacs);
  \draw[->] (luacs)--(acs);
  \draw[->] (UR)--(uacs);
  \draw[->] (LUR)--(sluacs);
  \draw[->] (WUR)--(wuacs);
  \draw[->] (WLUR)--(luacs);
  \draw[->] (R)--(acs);
\end{tikzpicture}
\end{center}
\CAP\label{fig:3}
\end{figure}

In Section \ref{sec:examples} we will see some examples which show that no other implications 
are valid in general.\par
Let us further remark that every space whose norm is \ac{UG} is also \ac{sluacs}, thus we have
the following diagram illustrating the connection to smoothness properties.

\begin{figure}[H]
\begin{center}
\begin{tikzpicture}
  \node (US) at (-2,1) {US};
  \node (UG) at (0,1) {UG};
  \node (S) at (2,1) {S};
  \node (uacs) at (-2,-0.5) {uacs};
  \node (sluacs) at (0,-0.5) {sluacs};
  \node (acs) at (2,-0.5) {acs};
  \draw[->] (US)--(UG);
  \draw[->] (UG)--(S);
  \draw[->] (uacs)--(sluacs);
  \draw[->] (sluacs)--(acs);
  \draw[->] (US)--(uacs);
  \draw[->] (UG)--(sluacs);
  \draw[->] (S)--(acs);
\end{tikzpicture}
\end{center}
\CAP\label{fig:4}
\end{figure}

In the next section we collect some general results on \ac{uacs} spaces and their relatives.

\section{Some general facts}\label{sec:general facts}
We start with the promised alternative proof of Proposition \ref{prop:char-uacs} which does
not rely on reflexivity. Instead, we shall employ the Bishop--Phelps--Bollob\'as theorem 
(cf. \cite{bollobas}*{Chap. 8, Theorem 11}), an argument that will also work for the case of \ac{sluacs} spaces.
This idea was suggested to the author by Dirk Werner.
\begin{proposition}\label{prop:char-uacs-sluacs}
A Banach space $X$ is \ac{uacs} \ifif for any two sequences $(x_n)_{n\in \N}$, $(y_n)_{n\in \N}$ in $S_X$ and 
every sequence $(x_n^*)_{n\in \N}$ in $S_{X^*}$ we have
\begin{equation}\label{eq:1}
\norm*{x_n+y_n}\to 2 \ \mathrm{and} \ x_n^*(x_n)=1 \ \forall n\in \N \ \Rightarrow \ x_n^*(y_n)\to 1.
\end{equation}
$X$ is \ac{sluacs} \ifif for every $x\in S_X$ and all sequences $(x_n)_{n\in \N}$, $(x_n^*)_{n\in \N}$ in $S_X$ 
resp. $S_{X^*}$ we have
\begin{equation}\label{eq:2}
\norm{x_n+x}\to 2 \ \mathrm{and} \ x_n^*(x_n)=1 \ \forall n\in \N \ \Rightarrow \ x_n^*(x)\to 1.
\end{equation}
\end{proposition}

\begin{Proof}
We only prove the statement for \ac{uacs} spaces, the proof for the \ac{sluacs} case is completely analogous. 
Furthermore, only the `if' part of the stated equivalence requires proof. So suppose  \eqref{eq:1} holds for any
two sequences in $S_X$ and all sequences in $S_{X^*}$.\par
Now if $(x_n)_{n\in \N}$ and $(y_n)_{n\in \N}$ are sequences in $S_X$ and $(x_n^*)_{n\in \N}$ is a sequence in 
$S_{X^*}$ such that $\norm*{x_n+y_n}\to 2$ and $x_n^*(x_n)\to 1$ we can choose a strictly increasing sequence 
$(n_k)_{k\in \N}$ in $\N$ such that $x_{n_k}^*(x_{n_k})>1-2^{-2k-2}$ holds for all $k\in \N$. By the already 
cited Bishop--Phelps--Bollob\'as theorem we can find sequences $(\tilde{x}_k)_{k\in \N}$ in $S_X$ and 
$(\tilde{x}_k^*)_{k\in \N}$ in $S_{X^*}$ such that $\tilde{x}_k^*(\tilde{x}_k)=1$, $\norm*{\tilde{x}_k-x_{n_k}}\leq 2^{-k}$ 
and $\norm*{\tilde{x}_k^*-x_{n_k}^*}\leq 2^{-k}$ for all $k\in \N$.\par
It follows that $\norm*{\tilde{x}_k-x_{n_k}}\to 0$ and $\norm*{\tilde{x}_k^*-x_{n_k}^*}\to 0$ and since 
$\norm*{x_n+y_n}\to 2$ we get that $\norm*{\tilde{x}_k+y_{n_k}}\to 2$.\par
But then we also have $\tilde{x}_k^*(y_{n_k})\to 1$, by our assumption, which in turn implies $x_{n_k}^*(y_{n_k})\to 1$.\par
In the same way we can show that every subsequence of $(x_n^*(y_n))_{n\in \N}$ has another subsequence that tends 
to one and hence $x_n^*(y_n)\to 1$ which completes the proof.
\end{Proof}

Next we would like to give characterisations of \ac{acs}/\ac{sluacs}/\ac{uacs} spaces that do not explicitly
involve the dual space. As mentioned before, a Banach space $X$ is \ac{acs} \ifif $x$ and $y$ are smooth 
points of the unit ball of the two-dimensional subspace $\lin\set{x,y}$ whenever $x,y\in S_X$ are such that 
$\norm{x+y}=2$.\par
It is possible to reformulate and refine this statement in the following way.
\begin{proposition}\label{prop:char acs without dual}
For any Banach space $X$ the following assertions are equivalent:
\begin{enumerate}[\upshape(i)]
\item $X$ is \ac{acs}.
\item For all $x,y\in S_X$ with $\norm{x+y}=2$ we have
\begin{equation*}
\lim_{t\to 0^+}\frac{\norm{x+ty}+\norm{x-ty}-2}{t}=0.
\end{equation*}
\item For all $x,y\in S_X$ with $\norm{x+y}=2$ we have
\begin{equation*}
\lim_{t\to 0^+}\frac{\norm{x-ty}-1}{t}=-1.
\end{equation*}
\item For all $x,y\in S_X$ with $\norm{x+y}=2$ there is some $1\leq p<\infty$ such that
\begin{equation*}
\lim_{t\to 0^+}\frac{\norm{x+ty}^p+\norm{x-ty}^p-2}{t^p}=0.
\end{equation*}
\item For all $x,y\in S_X$ with $\norm{x+y}=2$ there is some $1\leq p<\infty$ such that
\begin{equation*}
\lim_{t\to 0^+}\frac{(1+t)^p+\norm{x-ty}^p-2}{t^p}=0.
\end{equation*}
\end{enumerate}
\end{proposition}
The analogous characterisation for \ac{sluacs} spaces reads as follows.
\begin{proposition}\label{prop:char sluacs without dual}
For any Banach space $X$ the following assertions are equivalent:
\begin{enumerate}[\upshape(i)]
\item $X$ is \ac{sluacs}.
\item For every $\eps>0$ and every $y\in S_X$ there is some $\delta>0$ such that for all $t\in [0,\delta]$ 
and each $x\in S_X$ with $\norm{x+y}\geq2(1-t)$ we have
\begin{equation*}
\norm{x+ty}+\norm{x-ty}\leq 2+\eps t.
\end{equation*}
\item For every $\eps>0$ and every $y\in S_X$ there is some $\delta>0$ such that for all $t\in [0,\delta]$
and each $x\in S_X$ with $\norm{x+y}\geq2-t\delta$ we have
\begin{equation*}
\norm{x-ty}\leq1+t(\eps-1).
\end{equation*}
\item For every $y\in S_X$ there is some $1\leq p<\infty$ such that for every $\eps>0$ there exists $\delta>0$
such that for all $t\in [0,\delta]$ and each $x\in S_X$ with $\norm{x+y}\geq2(1-t)$ we have
\begin{equation*}
\norm{x+ty}^p+\norm{x-ty}^p\leq 2+\eps t^p.
\end{equation*}
\item For every $y\in S_X$ there is some $1\leq p<\infty$ such that for every $\eps>0$ there exists $\delta>0$
such that for all $t\in [0,\delta]$ and each $x\in S_X$ with $\norm{x+y}\geq2-t\delta$ we have
\begin{equation*}
(1+t)^p+\norm{x-ty}^p\leq 2+\eps t^p.
\end{equation*}
\end{enumerate}
\end{proposition}
Finally, we have the following characterisation for \ac{uacs} spaces.
\begin{proposition}\label{prop:char uacs without dual}
For any Banach space $X$ the following assertions are equivalent:
\begin{enumerate}[\upshape(i)]
\item $X$ is \ac{uacs}.
\item For every $\eps>0$ there exists some $\delta>0$ such that for every $t\in [0,\delta]$ and all $x,y\in S_X$
with $\norm{x+y}\geq 2(1-t)$ we have
\begin{equation*}
\norm{x+ty}+\norm{x-ty}\leq 2+\eps t.
\end{equation*}
\item For every $\eps>0$ there exists some $\delta>0$ such that for every $t\in [0,\delta]$ and all $x,y\in S_X$
with $\norm{x+y}\geq 2-\delta t$ we have
\begin{equation*}
\norm{x-ty}\leq 1+t(\eps-1).
\end{equation*}
\item There exists some $1\leq p<\infty$ such that for every $\eps>0$ there is some $\delta>0$ such that for all
$t\in [0,\delta]$ and all $x,y\in S_X$ with $\norm{x+y}\geq 2(1-t)$ we have
\begin{equation*}
\norm{x+ty}^p+\norm{x-ty}^p\leq 2+\eps t^p.
\end{equation*}
\item There exists some $1\leq p<\infty$ such that for every $\eps>0$ there is some $\delta>0$ such that for all
$t\in [0,\delta]$ and all $x,y\in S_X$ with $\norm{x+y}\geq 2-t\delta$ we have
\begin{equation*}
(1+t)^p+\norm{x-ty}^p\leq 2+\eps t^p.
\end{equation*}
\end{enumerate}
\end{proposition}

\begin{Proof}
We will only explicitly prove the characterisation for \ac{uacs} spaces. First we show $\mathrm{(i)} \Rightarrow \mathrm{(ii)}$.
So suppose $X$ is \ac{uacs} and fix $\eps>0$. Then there exists some $\tilde{\delta}>0$ such that for all $x,y\in S_X$ and 
$x^*\in S_{X^*}$ we have
\begin{equation*}
\norm{x+y}\geq 2(1-\tilde{\delta}) \ \mathrm{and} \ x^*(x)\geq 1-\tilde{\delta} \ \Rightarrow \ x^*(y)\geq 1-\eps.
\end{equation*}
Now if we put $\delta=\tilde{\delta}/2$ and take $t\in [0,\delta]$ and $x,y\in S_X$ such that $\norm{x+y}\geq 2(1-t)$ then we can
find a functional $x^*\in S_{X^*}$ such that $x^*(x-ty)=\norm{x-ty}$ and conclude that
\begin{equation*}
x^*(x)=\norm{x-ty}+tx^*(y)\geq 1-t-t=1-2t\geq 1-\tilde{\delta}.
\end{equation*}
By the choice of $\tilde{\delta}$ this implies $x^*(y)\geq 1-\eps$ and hence
\begin{equation*}
\norm{x+ty}+\norm{x-ty}=\norm{x+ty}+x^*(x-ty)\leq 1+t+1-tx^*(y)\leq 2+t\eps.
\end{equation*}
Now let us prove $\mathrm{(ii)} \Rightarrow \mathrm{(iii)}$. For a given $\eps>0$ choose $\delta>0$ to the value $\eps/2$
according to (ii). We may assume $\delta\leq \min\set{1,\eps/2}$.\par
Then if $t\in [0,\delta]$ and $x,y\in S_X$ such that $\norm{x+y}\geq 2-\delta t$ we in particular have $\norm{x+y}\geq 2(1-t)$
and hence 
\begin{equation*}
\norm{x+ty}+\norm{x-ty}\leq 2+t\frac{\eps}{2}.
\end{equation*}
But on the other hand 
\begin{equation*}
\norm{x+ty}\geq \norm{x+y}-(1-t)\norm{y}\geq 2-\delta t-1+t=1-\delta t+t\geq 1-\frac{\eps}{2}t+t.
\end{equation*}
It follows that $\norm{x-ty}\leq 1+t(\eps-1)$.\par
Next we prove that $\mathrm{(iii)} \Rightarrow \mathrm{(i)}$. Fix sequences $(x_n)_{n\in \N}$ and $(y_n)_{n\in \N}$ in $S_X$ 
such that $\norm{x_n+y_n}\to 2$ and a sequence $(x_n^*)_{n\in \N}$ of norm-one functionals with $x_n^*(x_n)\to 1$. Also, for every $n\in \N$
we fix $y_n^*\in S_{X^*}$ such that $y_n^*(y_n)=1$.\par
For given $\eps>0$ we choose $\delta>0$ according to (iii). For sufficiently large $n$ we have $\norm{x_n+y_n}\geq 2-\delta^2$ and 
$x_n^*(x_n)\geq 1-\eps\delta$ and hence
\begin{align*}
&(y_n^*-x_n^*)(\delta y_n)=x_n^*(x_n-\delta y_n)-x_n^*(x_n)+\delta\leq \norm{x_n-\delta y_n}+\delta-x_n^*(x_n) \\
&\leq\norm{x_n-\delta y_n}+\delta-1+\eps\delta\leq 1+\delta (\eps-1)+\delta -1+\eps\delta=2\delta\eps,
\end{align*}
where the last inequality holds because of $\norm{x_n+y_n}\geq 2-\delta^2$ and the choice of $\delta$.\par
It follows that $x_n^*(y_n)\geq y_n^*(y_n)-2\eps=1-2\eps$ for sufficiently large $n$.\par
The implications $\mathrm{(ii)} \Rightarrow \mathrm{(iv)}$ and $\mathrm{(iii)} \Rightarrow \mathrm{(v)}$ are clear. To prove
$\mathrm{(iv)} \Rightarrow \mathrm{(ii)}$ recall the inequalities
\begin{align*}
&(a+b)^p\leq 2^{p-1}(a^p+b^p) \ \ \forall a,b\geq 0, \forall p\in [1,\infty[ \\
&(a+b)^{\alpha}\leq a^{\alpha}+b^{\alpha} \ \ \forall a,b\geq 0, \forall \alpha\in ]0,1].
\end{align*}
They imply that for all $x,y\in S_X$, every $t>0$ and each $1\leq p<\infty$ one has
\begin{align*}
&\frac{\norm{x+ty}+\norm{x-ty}-2}{t}\leq \frac{\paren*{2^{p-1}\paren*{\norm{x+ty}^p+\norm{x-ty}^p}}^{1/p}-2}{t} \\
&\leq \paren*{\frac{2^{p-1}\paren*{\norm{x+ty}^p+\norm{x-ty}^p}-2^p}{t^p}}^{1/p} \\
&=2^{1-1/p}\paren*{\frac{\norm{x+ty}^p+\norm{x-ty}^p-2}{t^p}}^{1/p},
\end{align*}
which shows $\mathrm{(iv)} \Rightarrow \mathrm{(ii)}$. If we replace $\norm{x+ty}$ by $1+t$ in the above calculation,
we also obtain a proof for $\mathrm{(v)} \Rightarrow \mathrm{(iii)}$.
\end{Proof}

If we define the modulus $\rho_{\mathrm{uacs}}^X$ by 
\begin{equation*}
\rho_{\mathrm{uacs}}^X(\tau)=\sup\set*{1/2(\norm{x+\tau y}+\norm{x-\tau y})-1:(x,y)\in S_X(\tau)},
\end{equation*}
where $\tau>0$ and $S_X(\tau)=\set*{(x,y)\in S_X\times S_X:\norm{x+y}\geq 2(1-\tau)}$ then because of the equivalence 
of (i) and (ii) in Proposition \ref{prop:char uacs without dual} $X$ is \ac{uacs} \ifif 
$\lim_{\tau\to 0}\rho_{\mathrm{uacs}}^X(\tau)/\tau=0$ and obviously $\rho_{\mathrm{uacs}}^X(\tau)\leq \rho_X(\tau)$.\par
Let us also define 
\begin{equation*}
\tilde{\delta}_{\mathrm{uacs}}^X(\eps)=\inf\set*{\max\set*{1-\frac{1}{2}\norm{x+y},1-x^*(x)}:x,y\in S_X, x^*\in A_{\eps}(y)},
\end{equation*}
where $0<\eps\leq 2$ and  $A_{\eps}(y)=\set*{x^*\in S_{X^*}: x^*(y)\leq 1-\eps}$.\par
From the very definition of the \ac{uacs} spaces it follows that $X$ is \ac{uacs} \ifif  
$\tilde{\delta}_{\mathrm{uacs}}^X(\eps)>0$ for every $0<\eps\leq 2$.\par
Examining the proof of the implication $\mathrm{(i)} \Rightarrow \mathrm{(ii)}$ in Proposition \ref{prop:char uacs without dual} we 
see that the following holds.
\begin{lemma}\label{lemma:delta tilde rho}
If $X$ is a Banach space and $0<\eps\leq 2$ such that $\tilde{\delta}_{\mathrm{uacs}}^X(\eps)>0$ then for every $\tau>0$ with 
$2\tau<\tilde{\delta}_{\mathrm{uacs}}^X(\eps)$ we have $2\rho_{\mathrm{uacs}}^X(\tau)\leq\tau\eps$.
\end{lemma}
The reverse connection between $\rho_{\mathrm{uacs}}^X$ and $\delta_{\mathrm{uacs}}^X$ is given by the following lemma.
\begin{lemma}\label{lemma:delta rho}
Let $X$ be any Banach space and $\tau>0$ as well as $0<\eps\leq 2$. Then the inequality
\begin{equation*}
\delta_{\mathrm{uacs}}^X(\eps)\geq\frac{\eps\tau-2\rho_{\mathrm{uacs}}^X(\tau)}{2(\tau+1)}.
\end{equation*}
holds.
\end{lemma}

\begin{Proof}
We may assume $\eps\tau-2\rho_{\mathrm{uacs}}^X(\tau)>0$, because otherwise the inequality is trivially satisfied.
Let us put $R=(\eps\tau-2\rho_{\mathrm{uacs}}^X(\tau))(2(\tau+1))^{-1}$ and take $x,y\in S_X$ and $x^*\in S_{X^*}$ 
such that $x^*(x)=1$ and $\norm{x+y}>2(1-R)$.\par
Then we can find $z^*\in S_{X^*}$ with $z^*(x+y)>2(1-R)$ and hence $z^*(x)>1-2R$ and $z^*(y)>1-2R$.\par
It follows that 
\begin{align*}
&(z^*-x^*)(\tau y)=z^*(x+\tau y)+x^*(x-\tau y)-x^*(x)-z^*(x) \\
&\leq \norm{x+\tau y}+\norm{x-\tau y}-1-z^*(x)\leq 2\rho_{\mathrm{uacs}}^X(\tau)+1-z^*(x) \\
&\leq 2(\rho_{\mathrm{uacs}}^X(\tau)+R).
\end{align*}
Hence 
\begin{equation*}
x^*(y)\geq z^*(y)-\frac{2}{\tau}\big(\rho_{\mathrm{uacs}}^X(\tau)+R\big)> 1-2R-\frac{2}{\tau}\big(\rho_{\mathrm{uacs}}^X(\tau)+R\big)=1-\eps
\end{equation*}
and we are done.
\end{Proof}

Now we turn to the proof of the superreflexivity of \ac{uacs} spaces without using James's result on uniformly 
non-square Banach spaces. A key ingredient to James's proof is the following lemma of his, which may be 
found in \cite{beauzamy}*{p.51}.
\begin{lemma}\label{lemma:James reflexivity}
A Banach space $X$ is {\em not} reflexive \ifif for every $0<\theta<1$ there is a sequence $(x_k)_{k\in \N}$
in $B_X$ and a sequence $(x_n^*)_{n\in \N}$ in $B_{X^*}$ such that for every $n\in \N$ we have
\begin{equation*}
x_n^*(x_k)=
\begin{cases}
\theta & \mathrm{if} \ n\leq k \\
0 & \mathrm{if} \ n>k.
\end{cases}
\end{equation*}
\end{lemma}
Even armed with this lemma it is still difficult to prove the superreflexivity of uniformly 
non-square Banach spaces, but it easily yields the result for \ac{uacs} spaces. We can even 
prove a stronger result: it is a well known fact that a Banach space $X$ is reflexive if it 
satisfies $\liminf_{t\to 0^+}\rho_X(t)/t<1/2$ (cf.\,\cite{smithsullivan}*{Theorem 2}).\footnote{Note 
that the definition of $\rho_X$ given there differs from our definition by a factor $1/2$.} We will
see that the same holds if we replace $\rho_X$ by $\rho_{\mathrm{uacs}}^{X}$, even a bit more is true.
\begin{proposition}\label{prop:uacs reflexive}
If there is some $0<t$ such that $\rho_{\mathrm{uacs}}^{X}(t)<t/2$, then $X$ is superreflexive
(actually, it is uniformly non-square).
\end{proposition}

\begin{Proof}
Put $\theta=2\rho_{\mathrm{uacs}}^{X}(t)/t<1$ and choose $\eps>0$ such that $\theta+\eps<1$.
Also, put $\eta=\min\set*{t\eps/5,\eps/5}$.\par
If $x,y\in S_X$ such that $\norm{x+y}\geq 2(1-\eta)$ and $x^*\in S_{X^*}$ with $x^*(x)\geq 1-\eta$
fix $y^*\in S_{X^*}$ such that $y^*(x+y)\geq 2(1-\eta)$. Then $y^*(x)\geq 1-2\eta$ and $y^*(y)\geq 1-2\eta$
and hence
\begin{align*}
&(y^*-x^*)(ty)=y^*(x+ty)+x^*(x-ty)-x^*(x)-y^*(x) \\
&\leq \norm{x+ty}+\norm{x-ty}-2+3\eta\leq 2\rho_{\mathrm{uacs}}^{X}(t)+3\eta=t\theta+3\eta\leq(\theta+\frac{3}{5}\eps)t.
\end{align*}
Consequently, $x^*(y)\geq y^*(y)-\theta-\frac{3}{5}\eps\geq 1-2\eta-\theta-\frac{3}{5}\eps\geq 
1-\frac{2}{5}\eps-\theta-\frac{3}{5}\eps=1-(\theta+\eps)$.\par
Next we fix $0<\tau<1/2$ such that $\tau(1+(1-2\tau)^{-1})\leq\eta$ and put
$\beta=1-(1-\tau)(1-2\tau)(1-\theta-\eps)$. Then $0<\beta<1$.
\begin{claim*}
If $x,y\in B_X$ such that $\norm{x+y}\geq 2(1-\tau)$ and $x^*\in B_{X^*}$ such that $x^*(x)\geq 1-\tau$ then
$x^*(y)\geq 1-\beta$.
\end{claim*}
To see this, take $x,y$ and $x^*$ as above and observe $\norm{x}, \norm{y}\geq 1-2\tau$. Hence
\begin{align*}
&\norm*{\frac{x}{\norm{x}}+\frac{y}{\norm{y}}}\geq \frac{\norm{x+y}}{\norm{x}}-\abs*{\frac{1}{\norm{x}}-\frac{1}{\norm{y}}}\norm{y} \\
&\geq \norm{x+y}-\abs*{\frac{1}{\norm{x}}-\frac{1}{\norm{y}}}\geq 2(1-\tau)-\frac{2\tau}{1-2\tau}\geq 2(1-\eta)
\end{align*}
and moreover, since $\norm{x^*},\norm{x}\leq 1$,
\begin{equation*}
\frac{x^*}{\norm{x^*}}\paren*{\frac{x}{\norm{x}}}\geq 1-\tau\geq 1-\eta.
\end{equation*}
Thus by our previous considerations we must have
\begin{equation*}
x^*(y)\geq \norm{x^*}\norm{y}(1-\theta-\eps)\geq (1-\tau)(1-2\tau)(1-\theta-\eps)=1-\beta.
\end{equation*}
From the above claim together with the fact that $\beta<1$ it could be easily deduced that $X$ is uniformly 
non-square and hence superreflexive, but if we just want to prove the superreflexivity an application of 
Lemma \ref{lemma:James reflexivity} is enough. For if $X$ was not reflexive then by said Lemma  we could 
find sequences $(x_k)_{k\in \N}$ in $B_X$ and $(x_n^*)_{n\in \N}$ in $B_{X^*}$ such that $x_n^*(x_k)=0$ 
for $n>k$ and $x_n^*(x_k)=1-\tau$ for $n\leq k$.\par
We only need the first two members of the sequences to derive a contradiction, namely we have
$\norm{x_1+x_2}\geq x_1^*(x_1)+x_1^*(x_2)=2(1-\tau)$ and $x_2^*(x_2)=1-\tau$ but $x_2^*(x_1)=0<1-\beta$ 
contradicting our just established claim.\par
Thus $X$ must be reflexive and to prove the superreflexivity it only remains to show that for every 
Banach space $Y$ which is finitely representable in $X$ there exists $0<t^{\prime}$ such that 
$\rho_{\mathrm{uacs}}^{Y}(t^{\prime})<t^{\prime}/2$ which we will do in the next Lemma.
\end{Proof}

\begin{lemma}\label{lemma:superprop}
If there is some $0<t$ such that $\rho_{\mathrm{uacs}}^{X}(t)<t/2$ and $Y$ is finitely representable 
in $X$ then there is $0<t^{\prime}$ such that $\rho_{\mathrm{uacs}}^{Y}(t^{\prime})<t^{\prime}/2$.
\end{lemma}

\begin{Proof}
Let $\theta, \eps, \eta, \tau$ and $\beta$ be as in the previous proof. Put $\nu=\tau/4$.
\begin{claim*}
If $x,y\in B_X$ such that $\norm{x+y}\geq 2(1-\nu)$ then $\norm{x+\nu y}+\norm{x-\nu y}\leq 2+\nu\beta$.
\end{claim*}
To establish this, take $x,y\in B_X$ as above and also fix $x^*\in S_{X^*}$ such that $x^*(x-\nu y)=\norm{x-\nu y}$.
Observe as before that $\norm{x}, \norm{y}\geq 1-\tau/2$. Hence we have
\begin{equation*}
x^*(x)=\norm{x-\nu y}+x^*(\nu y)\geq \norm{x}-\nu \norm{y}+\nu x^*(y)\geq\norm{x}-2\nu\geq 1-\tau.
\end{equation*}
The claim we established in the previous proof now gives us $x^*(y)\geq 1-\beta$. It follows that
\begin{equation*}
\norm{x+\nu y}+\norm{x-\nu y}=\norm{x+\nu y}+x^*(x-\nu y)\leq 2+\nu(1-x^*(y))\leq 2+\nu\beta.
\end{equation*}
Next fix $\beta<\alpha<1$ and $0<\tilde{\eta}<\nu$ such that $(\beta\nu+3\tilde{\eta})(\nu-\tilde{\eta})^{-1}<\alpha$.
Put $t^{\prime}=\nu-\tilde{\eta}$. Finally, choose $\tilde{\eps}>0$ such that $(1-t^{\prime})(1+\tilde{\eps})^{-1}>1-\nu$
and $(1+\tilde{\eps})(2+\nu\beta)\leq 2+\nu\beta+\tilde{\eta}$.\par
Now take $y_1, y_2\in S_Y$ with $\norm{y_1+y_2}\geq 2(1-t^{\prime})$ and put $F=\lin\set*{y_1,y_2}$. Since $Y$ is finitely
representable in $X$ there is a subspace $E\ssq X$ and an isomorphism $T:F \rightarrow E$ such that $\norm{T}=1$ and
$\norm{T^{-1}}\leq 1+\tilde{\eps}$. Let $x_i=Ty_i$ for $i=1,2$.\par
It easily follows that $\norm{x_1+x_2}\geq 2(1-t^{\prime})(1+\tilde{\eps})^{-1}>2(1-\nu)$, whence 
$\norm{x_1+\nu x_2}+\norm{x_1-\nu x_2}\leq 2+\nu\beta$ which implies $\norm{y_1+\nu y_2}+\norm{y_1-\nu y_2}
\leq (1+\tilde{\eps})(2+\nu\beta)$. Thus we have
\begin{align*}
&\norm{y_1+t^{\prime}y_2}+\norm{y_1-t^{\prime}y_2}\leq \norm{y_1+\nu y_2}+\norm{y_1-\nu y_2}+2\abs*{\nu-t^{\prime}} \\
&\leq(1+\tilde{\eps})(2+\nu\beta)+2\tilde{\eta}\leq 2+\nu\beta+3\tilde{\eta}\leq 2+\alpha(\nu-\tilde{\eta})=2+\alpha t^{\prime}.
\end{align*}
So we have proved $2\rho_{\mathrm{uacs}}^{Y}(t^{\prime})/t^{\prime}\leq\alpha<1$.
\end{Proof}

We remark that the uniform non-squareness of a space $X$ satisfying $2\rho_{\mathrm{uacs}}^{X}(t)<t$ for some $0<t$ 
could also be deduced from our Lemma \ref{lemma:delta rho} and \cite{gao1}*{Theorem 2}, where it is observed that 
$\delta_{\mathrm{uacs}}^{X}(1)>0$ is sufficient to ensure that $X$ is uniformly non-square.\par
Now let us have a look at the quantitative connection between the moduli $\delta_{\mathrm{uacs}}^X$ and $\tilde{\delta}_{\mathrm{uacs}}^X$.
\begin{lemma}\label{lemma:delta deltatilde}
If $X$ is \ac{uacs} then 
\begin{equation*}
\tilde{\delta}_{\mathrm{uacs}}^X(\eps)\geq \delta_{\mathrm{ucas}}^X\paren*{\delta_{\mathrm{uacs}}^X(\eps)}
\end{equation*}
for every $0<\eps\leq 2$.
\end{lemma}

\begin{Proof}
Here we can adopt Sirotkin's idea from the proof of Proposition \ref{prop:char-uacs} in \cite{sirotkin}.
Put $\delta=\delta_{\mathrm{uacs}}^X\paren*{\delta_{\mathrm{uacs}}^X(\eps)}$ and take $x,y\in S_X$ and $x^*\in S_{X^*}$
such that $\norm{x+y}>2(1-\delta)$ and $x^*(x)>1-\delta$.\par
Since $X$ is reflexive, there is some $z\in S_X$ with $x^*(z)=1$. It follows that $\norm{x+z}\geq x^*(x+z)>2(1-\delta)$.\par
Now fix $y^*\in S_{X^*}$ such that $y^*(x)=1$. Then by the definition of $\delta$ we must have $y^*(z)>1-\delta_{\mathrm{uacs}}^X(\eps)$
and $y^*(y)>1-\delta_{\mathrm{uacs}}^X(\eps)$ and hence $\norm{y+z}>2(1-\delta_{\mathrm{uacs}}^X(\eps))$.\par
Because of $x^*(z)=1$ this implies $x^*(y)>1-\eps$ and the proof is finished.
\end{Proof}

It is claimed in \cite{dhompongsa0}*{Lemma 3.10} that the modulus of $U$-convexity, which coincides with our modulus 
$\delta_{\mathrm{ucas}}^{X}$, is continuous on $]0,2[$, but it seems that the proof given there only works in the case $\eps<1$
(this is not a major drawback since one is usually interested in small values of $\eps$). We wish to point out that for values 
between $0$ and $1$ even more is true, namely $\delta_{\mathrm{uacs}}^{X}$ is uniformly continuous on $[a,1[$ for every $0<a<1$.
\begin{lemma}\label{lemma:uacs mod unif cont}
For every Banach space $X$ and all $0<\eps,\eps^{\prime}<1$ we have
\begin{equation*}
\abs*{\delta_{\mathrm{uacs}}^{X}(\eps)-\delta_{\mathrm{uacs}}^{X}(\eps^{\prime})}\leq\frac{\abs*{\eps-\eps^{\prime}}}{\min\set*{\eps,\eps^{\prime}}}.
\end{equation*}
In particular, $\delta_{\mathrm{uacs}}^{X}$ is uniformly continuous on $[a,1[$ for all $0<a<1$.
\end{lemma}

\begin{Proof}
Let $0<\eps<1$ and $0<\beta<1-\eps$. Put $\tau=\beta/(\eps+\beta)$ and take $x,y\in S_X$  and $x^*\in S_{X^*}$
such that $x^*(x)=1$ and $x^*(y)\leq 1-\eps$. Let $z=(y-\tau x)/\norm{y-\tau x}$. Note that, since 
$\norm{y-\tau x}\geq 1-\tau$ and $\eps+\tau<1$, we have
\begin{equation*}
x^*(z)\geq\frac{1-\eps-\tau}{1-\tau}=1-\eps\paren*{1+\frac{\tau}{1-\tau}}=1-(\eps+\beta)
\end{equation*}
and hence
\begin{equation*}
1-\norm*{\frac{x+z}{2}}\geq\delta_{\mathrm{uacs}}^{X}(\eps+\beta).
\end{equation*}
Furthermore, we have
\begin{equation*}
\norm*{y-z}\leq \frac{\norm*{\paren{\norm{y-\tau x}-1}y+\tau x}}{1-\tau}\leq\frac{2\tau}{1-\tau}=\frac{2\beta}{\eps}.
\end{equation*}
It follows that
\begin{equation*}
1-\norm*{\frac{x+y}{2}}\geq\delta_{\mathrm{uacs}}^{X}(\eps+\beta)-\frac{\beta}{\eps}.
\end{equation*}
Thus we have
\begin{equation*}
\delta_{\mathrm{uacs}}^{X}(\eps+\beta)\geq\delta_{\mathrm{uacs}}^{X}(\eps)\geq\delta_{\mathrm{uacs}}^{X}(\eps+\beta)-\frac{\beta}{\eps}
\end{equation*}
for all $0<\eps<1$ and every $0<\beta<1-\eps$, which finishes the proof.
\end{Proof}

Let us mention yet another characterisation of \ac{uacs} spaces: if $1<p<\infty$ is fixed then a Banach space $X$ is \ac{uacs} 
\ifif for all bounded sequences $(x_n)_{n\in \N}$ and $(y_n)_{n\in \N}$ in $X$ and every bounded sequence $(x_n^*)_{n\in \N}$ in $S_{X^*}$ 
the two conditions $2^{p-1}(\norm{x_n}^p+\norm{y_n}^p)-\norm{x_n+y_n}^p\to 0$ and $x_n^*(x_n)-\norm{x_n^*}\norm{x_n}\to 0$
imply $x_n^*(y_n)-\norm{x_n^*}\norm{y_n}\to 0$. We had already mentioned the special case $p=2$ earlier. The proof is
completely analogous to the one for the well-known corresponding characterisation of \ac{UR} spaces. Similar 
characterisations hold for \ac{acs}, \ac{luacs}, \ac{wuacs} and \ac{sluacs} spaces.\par
Next we will deal with some duality results. In \cite{lau}*{Theorem 2.4} a proof of the fact that a Banach space 
$X$ is a $U$-space \ifif its dual $X^*$ is a $U$-space is proposed and in \cite{dutta}*{Theorem 2.6} the stronger 
statement that the moduli of u-convexity of $X$ and $X^*$ coincide is claimed. Both proofs make use of the following
claim from \cite{lau}*{Remark after Definition 2.2}:
\begin{claim*}
$X$ is a $U$-space \ifif for every $\eps>0$ there is some $\delta>0$ such that whenever $x,y\in S_X$ and 
$x^*,y^*\in S_{X^*}$ with $x^*(x)=1=y^*(y)$ and $\norm{x+y}>2(1-\delta)$ then $\norm{x^*+y^*}>2(1-\eps)$.
\end{claim*}
\noindent A $U$-space certainly has the above property. However, the converse need not be true, not even in a two-dimensional space.\par
To see this, first note that if $X$ is finite-dimensional then by an easy compactness argument the condition of the claim 
is equivalent to the following one: whenever $x,y\in S_X$ and $x^*,y^*\in S_{X^*}$ with $x^*(x)=1=y^*(y)$ and $\norm{x+y}=2$ we 
also have $\norm{x^*+y^*}=2$.\par
Therefore, if $X$ is finite-dimensional it fulfils the condition of the claim if for each $x,y\in S_X$ with $\norm{x+y}=2$ at
least one of the two points $x$ and $y$ is a smooth point of the unit ball. But as we have mentioned before, a two-dimensional
space is \ac{acs} (equivalently a $U$-space) \ifif whenever $x,y\in S_X$ with $\norm{x+y}=2$ then {\em both} points $x$ and $y$
are smooth points of the unit ball.\par
Taking all this into account, we see that the space $\R^2$ endowed with the norm whose unit ball is sketched below will be an example
of a space which fulfils the condition of the claim but is not a $U$-space.

\begin{figure}[H]
\begin{center}
\begin{tikzpicture}
  \draw[help lines] (0,-2.5)--(0,2.5);
  \draw[help lines] (-2.5,0)--(2.5,0);
  \draw (-2,0)[rounded corners]--(0,2)--(2,0);
  \draw (2,0)[rounded corners]--(0,-2)--(-2,0);
\end{tikzpicture}
\end{center}
\CAP\label{fig:5}
\end{figure}

\noindent However, it is possible to modify the proof from \cite{lau}*{Theorem 2.4} to show that the desired self-duality result is true
nonetheless.
\begin{proposition}\label{prop:uacs-dual}
Let $X$ be a Banach space whose dual $X^*$ is \ac{uacs}. Then we have
\begin{equation}\label{eq:3}
\delta_{\mathrm{uacs}}^X(\eps)\geq \delta_{\mathrm{uacs}}^{X^*}\paren*{\delta_{\mathrm{uacs}}^{X^*}(\eps)} \ \forall \eps\in ]0,2].
\end{equation}
In particular, $X$ is also \ac{uacs}.
\end{proposition}

\begin{Proof}
Take any $\eps\in ]0,2]$ and put $\delta=\delta_{\mathrm{uacs}}^{X^*}(\eps)$ and 
$\tilde{\delta}=\delta_{\mathrm{uacs}}^{X^*}(\delta)$.\par
Now if $x,y\in S_X$ and $x^*\in S_{X^*}$ with $x^*(x)=1$ and $\norm*{x+y}>2(1-\tilde{\delta})$
choose $y^*, z^*\in S_{X^*}$ such that $y^*(y)=1$ and $z^*(x+y)=\norm*{x+y}$.\par
Then we must have $z^*(x)>1-2\tilde{\delta}$ and $z^*(y)>1-2\tilde{\delta}$. It follows that 
$(z^*+x^*)(x)>2-2\tilde{\delta}$ and $(z^*+y^*)(y)>2-2\tilde{\delta}$ and hence
\begin{equation}\label{eq:4}
\norm*{\frac{z^*+x^*}{2}}>1-\tilde{\delta} \ \mathrm{and} \ \norm*{\frac{z^*+y^*}{2}}>1-\tilde{\delta}.
\end{equation}
Next we pick any $z^{**}\in S_{X^{**}}$ with $z^{**}(z^*)=1$. Then from \eqref{eq:4} and the 
definition of $\tilde{\delta}$ we get that $z^{**}(x^*)>1-\delta$ and $z^{**}(y^*)>1-\delta$.\par
It follows that $\norm*{x^*+y^*}>2(1-\delta)$ and because of $y^*(y)=1$ and the definition of 
$\delta$ this implies $x^*(y)>1-\eps$ and thus we have shown 
$\delta_{\mathrm{uacs}}^X(\eps)\geq\tilde{\delta}=\delta_{\mathrm{uacs}}^{X^*}\paren*{\delta_{\mathrm{uacs}}^{X^*}(\eps)}$.
\end{Proof}

Taking into account that \ac{uacs} spaces are reflexive we finally get that being \ac{uacs} is a self-dual property.
\begin{corollary}\label{cor:uacs-dual}
A Banach space $X$ is \ac{uacs} \ifif $X^*$ is \ac{uacs}.
\end{corollary}
The author does not know whether the equality $\delta_{\mathrm{uacs}}^X=\delta_{\mathrm{uacs}}^{X^*}$ that 
was claimed in \cite{dutta}*{Theorem 2.6} is actually true.\par

Alternatively, we could also derive the self-duality from the following lemma (cf. the proof of 
\cite{fabian}*{Lemma 9.9}). The modulus $\tilde{\rho}_{\mathrm{uacs}}^X$ is defined exactly as 
$\rho_{\mathrm{uacs}}^X$ except that one replaces $S_X$ by $B_X$. The argument that $X$ is \ac{uacs} 
\ifif $\lim_{\tau \to 0}\tilde{\rho}_{\mathrm{uacs}}^X(\tau)/\tau=0$ is analogous to the one for 
$\rho_{\mathrm{uacs}}^X$.
\begin{lemma}\label{lemma:dual rho delta}
If $X$ is any Banach space then for every $\tau>0$ and every $0<\eps\leq 2$ the following inequalities hold:
\begin{enumerate}[\upshape(i)]
\item $\delta_{\mathrm{uacs}}^X(\eps)+\rho_{\mathrm{uacs}}^{X^*}(\tau)\geq\tau\frac{\eps}{2}$
\item $\delta_{\mathrm{uacs}}^{X^*}(\eps)+\tilde{\rho}_{\mathrm{uacs}}^X(\tau)\geq\tau\frac{\eps}{2}$
\end{enumerate}
\end{lemma}

\begin{Proof}
We only prove the slightly more difficult inequality (ii). To this end, fix $x^*,y^*\in S_{X^*}$ and $x^{**}\in S_{X^{**}}$ 
such that $x^{**}(x^*)=1$ and $x^{**}(y^*)\leq 1-\eps$.\par
If $\norm{x^*+y^*}\leq 2(1-\tau)$ then we certainly have $2-\norm{x^*+y^*}\geq \tau\eps-2\tilde{\rho}_{\mathrm{uacs}}^X(\tau)$.\par
If $\norm{x^*+y^*}>2(1-\tau)$ then take an arbitrary $0<\alpha<\norm{x^*+y^*}-2(1-\tau)$. By Goldstine's theorem there is some
$x\in B_X$ such that 
\begin{equation*}
\abs*{x^{**}(x^*)-x^*(x)}\leq \frac{\alpha}{2} \ \mathrm{and} \ \abs*{x^{**}(y^*)-y^*(x)}\leq \frac{\alpha}{2}.
\end{equation*}
Now choose $y\in S_X$ such that $(x^*+y^*)(y)>\norm{x^*+y^*}-\alpha/2$. It follows that $(x^*+y^*)(y)>2(1-\tau)+\alpha/2$
and hence $x^*(y), y^*(y)>1-2\tau+\alpha/2$.\par
Thus we have
\begin{equation*}
\norm{x+y}\geq x^*(x+y)\geq x^{**}(x^*)-\frac{\alpha}{2}+1-2\tau+\frac{\alpha}{2}=2(1-\tau)
\end{equation*}
and hence
\begin{align*}
&2\tilde{\rho}_{\mathrm{uacs}}^X(\tau)\geq\norm{y+\tau x}+\norm{y-\tau x}-2\geq x^*(y+\tau x)+y^*(y-\tau x)-2 \\
&=(x^*+y^*)(y)+\tau(x^*(x)-y^*(x))-2 \\
&\geq \norm{x^*+y^*}-\frac{\alpha}{2}+\tau\paren*{x^{**}(x^*)-x^{**}(y^*)-\alpha}-2 \\
&\geq \norm{x^*+y^*}-\frac{\alpha}{2}+\tau(\eps-\alpha)-2.
\end{align*}
For $\alpha \to 0$ we get $2-\norm{x^*+y^*}\geq \tau\eps-2\tilde{\rho}_{\mathrm{uacs}}^X(\tau)$ and we are done.
\end{Proof}

There are also some duality result on \ac{acs}, \ac{luacs}, \ac{sluacs} and \ac{wuacs} spaces which we will treat 
in the following. The proof of the first statement is very easy and will therefore be omitted.
\begin{proposition}\label{prop:dual acs}
A Banach space $X$ is \ac{acs} \ifif for all $x^*, y^*\in S_{X^*}$ and all $x,y\in S_X$ the implication
\begin{equation*}
(x^*+y^*)(x)=2 \ \mathrm{and} \ x^*(y)=1 \ \Rightarrow \ y^*(y)=1
\end{equation*}
holds. In particular, if $X^*$ is \ac{acs} then so is $X$ and the converse is true if $X$ is reflexive.
\end{proposition}
We will say that a dual space $X^*$ is $\mathrm{luacs}^*$ resp. $\mathrm{wuacs}^*$ if it fulfils the definition
of an \ac{luacs} resp. \ac{wuacs} space with for all weak*-continuous functionals on $X^*$. With this terminology
the following is valid.
\begin{proposition}\label{prop:dual}
For any Banach space $X$ we have the following implications.
  \begin{enumerate}[\upshape(i)]
  \item $X^* \ \mathrm{luacs}^* \ \iff \ X \ \mathrm{\ac{luacs}}$
  \item $X^* \ \mathrm{wuacs}^* \ \iff \ X \ \mathrm{\ac{sluacs}}$
  \item $X^* \ \mathrm{\ac{sluacs}} \ \iff \ X \ \mathrm{\ac{wuacs}}$
  \end{enumerate}
\end{proposition}

\begin{Proof}
We only prove (iii). Let us first assume that $X^*$ is \ac{sluacs} and take sequences 
$(x_n)_{n\in \N}$, $(y_n)_{n\in \N}$ in $S_X$ and a functional $x^*\in S_{X^*}$
such that $\norm*{x_n+y_n}\to 2$ and $x^*(x_n)\to 1$.\par
Choose a sequence $(x_n^*)_{n\in \N}$ in $S_{X^*}$ with $x_n^*(x_n+y_n)=\norm*{x_n+y_n}$
for every $n$. It follows that $x_n^*(x_n)\to 1$ and $x_n^*(y_n)\to 1$.\par
From $x^*(x_n)\to 1$ and $x_n^*(x_n)\to 1$ we get $\norm*{x_n^*+x^*}\to 2$.
Together with $x_n^*(y_n)\to 1$ and the fact that $X^*$ is \ac{sluacs} this 
implies $x^*(y_n)\to 1$ and we are done.\par
Now assume $X$ is \ac{wuacs} and fix a sequence $(x_n^*)_{n\in \N}$ in $S_{X^*}$ and $x^*\in S_{X^*}$
such that $\norm{x_n^*+x^*}\to 2$ as well as a sequence $(x_n^{**})_{n\in \N}$ in $S_{X^{**}}$ with
$x_n^{**}(x_n^*)\to 1$.\par
Because of $\norm{x_n^*+x^*}\to 2$ we can find a sequence $(x_n)_{n\in \N}$ in $S_X$ such that $x_n^*(x_n)\to 1$ 
and $x^*(x_n)\to 1$.\par
By Goldstine's theorem we can also find a sequence $(y_n)_{n\in \N}$ in $B_X$ which satisfies 
\begin{equation*}
\abs*{x_n^*(y_n)-x_n^{**}(x_n^*)}\leq \frac{1}{n} \ \mathrm{and} \ \abs*{x^*(y_n)-x_n^{**}(x^*)}\leq \frac{1}{n} \ \ \forall n\in \N.
\end{equation*}
So we have $x_n^*(x_n+y_n)\to 2$ and hence $\norm{x_n+y_n}\to 2$. Since $X$ is \ac{wuacs} and $x^*(x_n)\to 1$ we must also have
$x^*(y_n)\to 1$ and consequently $x_n^{**}(x^*)\to 1$.
\end{Proof}

If $X$ is reflexive then by (i) and (ii) of the preceding proposition $X^*$ is \ac{luacs} (resp. \ac{wuacs}) \ifif $X$ is \ac{luacs}
(resp. \ac{sluacs}). Next we would like to give necessary and sufficient conditions for a dual space to be \ac{acs} resp. \ac{luacs} 
resp. \ac{wuacs} that do not explicitly involve the bidual space. We start with the \ac{acs} case. The characterisation is inspired by
\cite{yorke}*{Proposition 3}.\par
\begin{proposition}\label{prop:char X* acs}
Let $X$ be any Banach space. The dual space $X^*$ is \ac{acs} \ifif for all sequences $(x_n)_{n\in \N}$ and $(y_n)_{n\in \N}$
in $B_X$ and all functionals $x^*, y^*\in S_{X^*}$ the implication
\begin{equation*}
x^*(x_n+y_n)\to 2 \ \mathrm{and} \ y^*(x_n)\to 1 \ \Rightarrow \ y^*(y_n)\to 1
\end{equation*}
holds.
\end{proposition}

\begin{Proof}
To prove the necessity, assume that $X^*$ is \ac{acs} and take sequences $(x_n)_{n\in \N}, (y_n)_{n\in \N}$ and functionals $x^*, y^*$ as above.
It follows that $\norm{x^*+y^*}=2$. By the weak*-compactness of $B_{X^{**}}$ we can find for an arbitrary subsequence $(y_{n_k})_{k\in \N}$ 
a subnet $(y_{n_{\phi(i)}})_{i\in I}$ that weak*-converges to some $y^{**}\in B_{X^{**}}$. It follows that $y^{**}(x^*)=1$ and since $X^*$ is 
\ac{acs} we must also have $y^{**}(y^*)=1$. Thus $y^*(y_{n_{\phi(i)}})\to 1$ and the proof of the necessity is finished.\par
Now assume that $X^*$ fulfils the above condition and take $x^*, y^*\in S_{X^*}$ and $x^{**}\in S_{X^{**}}$ such that 
$\norm{x^*+y^*}=2$ and $x^{**}(x^*)=1$. Then we can find a sequence $(x_n)_{n\in \N}$ in $B_X$ such that $x^*(x_n)\to 1$
and $y^*(x_n)\to 1$.\par
By Goldstine's theorem there is a sequence $(y_n)_{n\in \N}$ in $B_X$ such that $x^*(y_n)\to x^{**}(x^*)=1$ and $y^*(y_n)\to x^{**}(y^*)$.\par
Thus we have $x^*(x_n+y_n)\to 2$ and $y^*(x_n)\to 1$ and hence by our assumption we get $y^*(y_n)\to 1$, so $x^{**}(y^*)=1$.
\end{Proof}

The characterisations for the dual space to be \ac{luacs} resp. \ac{wuacs} are a bit more complicated. They read as follows.
\begin{proposition}
Let $X$ be a Banach space.
\begin{enumerate}[\upshape(i)]
\item $X^*$ is \ac{luacs} \ifif for every $x^*\in S_{X^*}$ and all sequences $(x_n^*)_{n\in \N}$ and $(x_k)_{k\in \N}$ in $S_{X^*}$ and
$B_X$ respectively, the implication
\begin{equation*}
\norm{x^*+x_n^*}\to 2 \ \mathrm{and} \ x_n^*(x_k)\xrightarrow[k\geq n]{k,n\to \infty} 1 \ \Rightarrow \ x^*(x_k)\to 1
\end{equation*}
holds.
\item $X^*$ is \ac{wuacs} \ifif for all sequences $(x_n^*)_{n\in \N}, (y_n^*)_{n\in \N}$ in $S_{X^*}$ and $(x_k)_{k\in \N}$ in $B_X$ 
the implication
\begin{equation*}
\norm{x_n^*+y_n^*}\to 2 \ \mathrm{and} \ x_n^*(x_k)\xrightarrow[k\geq n]{k,n\to \infty} 1 \ \Rightarrow \ \lim_{n\to \infty}\sup_{k\geq n}y_n^*(x_k)=1.
\end{equation*}
holds.
\end{enumerate}
\end{proposition}

\begin{Proof}
To prove (ii) we first assume that $X^*$ is \ac{wuacs} and fix sequences $(x_n^*)_{n\in \N}, (y_n^*)_{n\in \N}$ in $S_{X^*}$
and $(x_k)_{k\in \N}$ in $B_X$ as above. Since $B_{X^{**}}$ is weak*-compact there is a subnet $(x_{\phi(i)})_{i\in I}$ that
is weak*-convergent to some $x^{**}\in B_{X^{**}}$. We will show that $x^{**}(x_n^*)\to 1$.\par
Given any $\eps>0$ by our assumption on $(x_n^*)_{n\in \N}$ and $(x_k)_{k\in \N}$ we can find an $N\in \N$ such that
\begin{equation*}
\abs*{x_n^*(x_k)-1}\leq\eps \ \ \forall k\geq n\geq N.
\end{equation*}
For every $n\geq N$ it is possible to find an index $i\in I$ with $\phi(i)\geq n$ and $\abs*{x_n^*(x_{\phi(i)})-x^{**}(x_n^*)}\leq\eps$.
It follows that $\abs*{x^{**}(y_n^*)-1}\leq 2\eps$ and the convergence is proved.\par
So we have $\norm{x_n^*+y_n^*}\to 2$ and $x^{**}(x_n^*)\to 1$. Since $X^*$ is \ac{wuacs} this implies $x^{**}(y_n^*)\to 1$. Thus for
any $\delta>0$ there is some $n_0\in \N$ such that $\abs*{x^{**}(y_n^*)-1}\leq\delta$ for all $n\geq n_0$ and for any such $n$ we find $j\in I$
with $\phi(j)\geq n$ and $\abs*{y_n^*(x_{\phi(i)})-x^{**}(y_n^*)}\leq\delta$. Hence $\abs*{y_n^*(x_{\phi(i)})-1}\leq2\delta$ and we have
shown $\sup_{k\geq n}y_n^*(x_k)\geq 1-2\delta$ for all $n\geq n_0$.\par
Now let us prove the converse. We take sequences $(x_n^*)_{n\in \N}, (y_n^*)_{n\in \N}$ in $S_{X^*}$ such that $\norm{x_n^*+y_n^*}\to 2$ 
and a functional $x^{**}\in S_{X^{**}}$ with $x^{**}(x_n^*)\to 1$.\par
By means of Goldstine's theorem we find a sequence $(x_k)_{k\in \N}$ in $B_X$ that satisfies
\begin{equation*}
\abs*{x_n^*(x_k)-x^{**}(x_n^*)}\leq \frac{1}{k} \ \mathrm{and} \ \abs*{y_n^*(x_k)-x^{**}(y_n^*)}\leq \frac{1}{k} \ \ \forall n\leq k.
\end{equation*}
It is then easy to see that $(x_n^*(x_k))_{k\geq n}$ tends to $1$ and hence our assumption gives us $\lim_{n\to \infty}\sup_{k\geq n}y_n^*(x_k)=1$.\par
Thus for any $\eps>0$ there exists $N\in \N$ with $\sup_{k\geq n}y_n^*(x_k)>1-\eps$ and $1/n\leq \eps$ for each $n\geq N$.\par
If we fix $n\geq N$ we find $k\geq n$ with $y_n^*(x_k)\geq 1-\eps$ and because of $\abs*{x^{**}(y_n^*)-y_n^*(x_k)}\leq 1/k\leq \eps$ it follows
that $x^{**}(y_n^*)\geq 1-2\eps$ and the proof is finished. Part (i) is proved similarly.
\end{Proof}

One can also give some more characterisations of \ac{acs}, \ac{luacs} and \ac{sluacs} spaces by apparently stronger properties.
\begin{proposition}\label{prop:very acs}
For a Banach space $X$, the following assertions are equivalent:
  \begin{enumerate}[\upshape(i)]
  \item $X$ is \ac{acs}
  \item For all sequences $(x_n^*)_{n\in \N}, (y_n^*)_{n\in \N}$ in $B_{X^*}$ and all $x,y\in S_X$ the implication
  \begin{equation*}
  (x_n^*+y_n^*)(x)\to 2 \ \mathrm{and} \ y_n^*(y)\to 1 \ \Rightarrow \ x_n^*(y)\to 1
  \end{equation*}
  holds.
  \item For every sequence $(x_n^{*})_{n\in \N}$ in $S_{X^{*}}$ and all $x,y\in S_X$ the implication
  \begin{equation*}
  \norm{x+y}=2 \ \mathrm{and} \ x_n^{*}(x)\to 1 \ \Rightarrow \ x_n^{*}(y)\to 1
  \end{equation*}
  holds.
  \end{enumerate}
\end{proposition}

\begin{Proof}
$\mathrm{(i)} \Rightarrow \mathrm{(ii)}$ follows from Proposition \ref{prop:dual acs} together with the fact that $B_{X^*}$ is weak*-compact, 
the implication $\mathrm{(iii)} \Rightarrow \mathrm{(i)}$ is trivial and $\mathrm{(ii)} \Rightarrow \mathrm{(iii)}$ is also quite easy to see.
\end{Proof}

By means of Goldstine's theorem one can also prove the following cha\-racterisation of \ac{luacs} spaces (we omit the details).
\begin{proposition}\label{prop:very luacs}
A Banach space $X$ is \ac{luacs} if and only if for every sequence $(x_n^{**})_{n\in \N}$ in $S_{X^{**}}$, every $x\in S_X$ and 
each $x^*\in S_{X^*}$ the implication
\begin{equation*}
\norm{x_n^{**}+x}\to 2 \ \mathrm{and} \ x_n^{**}(x^*)\to 1 \ \Rightarrow \ x^*(x)=1.
\end{equation*}
holds.
\end{proposition}

Let us denote by $X^{(k)}$ the $k$-th dual of $X$. Then $X$ resp. $X^*$ naturally embeds into $X^{(2k)}$ resp. $X^{(2k+1)}$ for each $k$.
For \ac{sluacs} spaces we have the following stronger result.
\begin{proposition}\label{prop:very sluacs}
A Banach space $X$ is \ac{sluacs} \ifif for every $k\in \N$, for every sequence $(z_n)_{n\in \N}$ in $B_{X^{(2k)}}$, every $x\in S_X$ and each sequence 
$(z_n^*)_{n\in \N}$ in $B_{X^{(2k+1)}}$ the implication
\begin{equation*}
\norm{z_n+x}\to 2 \ \mathrm{and} \ z_n^*(z_n)\to 1 \ \Rightarrow \ z_n^*(x)\to 1
\end{equation*}
holds.
\end{proposition}

\begin{Proof}
The sufficiency is obvious. To prove the necessity, we first take sequences $(x_n^{**})_{n\in \N}$ in $B_{X^{**}}$
and $(x_n^{***})_{n\in \N}$ in $B_{X^{***}}$ as well as an element $x\in S_X$ such that $\norm{x_n^{**}+x}\to 2$ and
$x_n^{***}(x_n^{**})\to 1$. Then we can find a sequence $(y_n^*)_{n\in \N}$ in $S_{X^*}$ such that $x_n^{**}(y_n^*)\to 1$ 
and $y_n^*(x)\to 1$.\par
By Goldstine's theorem (applied to $X^*$) there is a sequence $(x_n^*)_{n\in \N}$ in $B_{X^*}$ such that $x_n^{***}(x_n^{**})-x_n^{**}(x_n^*)\to 0$ 
and $x_n^{***}(x)-x_n^*(x)\to 0$. Hence $x_n^{**}(x_n^*)\to 1$.\par
Again by Goldstine's theorem (now applied to $X$) there exists a sequence $(x_n)_{n\in \N}$ in $B_X$ such that $x_n^{**}(x_n^*)-x_n^*(x_n)\to 0$
and $x_n^{**}(y_n^*)-y_n^*(x_n)\to 0$. It follows that $x_n^*(x_n)\to 1$ and $y_n^*(x_n)\to 1$.\par
Taking into account that $y_n^*(x)\to 1$ we get $\norm{x_n+x}\to 2$. Since $X$ is \ac{sluacs} it follows $x_n^*(x)\to 1$ and hence $x_n^{***}(x)\to 1$.\par
Thus we have proved our claim for $k=1$. Continuing by induction with the above argument we can show it for all $k\in \N$.
\end{Proof}

If we use the preceding proposition and the technique from the proof of Proposition \ref{prop:char uacs without dual} we see that the following holds.
\begin{proposition}
For a Banach space $X$ the following assertions are equivalent:
\begin{enumerate}[\upshape(i)]
\item $X$ is \ac{sluacs}.
\item For every $k\in \N$, every $\eps>0$ and every $y\in S_X$ there is some $\delta>0$ such that for all $t\in [0,\delta]$ 
and each $z\in S_{X^{(2k)}}$ with $\norm{z+y}\geq2(1-t)$ we have
\begin{equation*}
\norm{z+ty}+\norm{z-ty}\leq 2+\eps t.
\end{equation*}
\item For every $k\in \N$, every $\eps>0$ and every $y\in S_X$ there is some $\delta>0$ such that for all $t\in [0,\delta]$
and each $z\in S_{X^{(2k)}}$ with $\norm{z+y}\geq2-t\delta$ we have
\begin{equation*}
\norm{z-ty}\leq1+t(\eps-1).
\end{equation*}
\end{enumerate}
\end{proposition}
Finally, let us consider quotient spaces. If $U$ is a closed subpace of $X$ then $(X/U)^*$ is 
isometrically isomorphic to $U^\perp$ (the annihilator of $U$ in $X^*$). Using this together 
with Corollary \ref{cor:uacs-dual} and the obvious fact that closed subspaces of \ac{uacs} 
spaces are again \ac{uacs}, one immediately gets that quotients of \ac{uacs} spaces are \ac{uacs}
as well. An analogous argument using part (iii) of Proposition \ref{prop:dual} works for \ac{wuacs} 
spaces, so in the summary we have the following proposition.\par
\begin{proposition}\label{prop:quot uacs}
Let $U$ be a closed subspace of the Banach space $X$. If $X$ is \ac{uacs} (resp. \ac{wuacs}) then
$X/U$ is also \ac{uacs} (resp. \ac{wuacs}).
\end{proposition}
As for quotients of \ac{acs}, \ac{luacs} and \ac{sluacs} spaces we have the following result
which is an analogue of \cite{klee}*{Proposition 3.2}.
\begin{proposition}\label{prop:quot acs}
If $U$ is a reflexive subspace of the Banach space $X$ then the properties \ac{acs}, \ac{luacs}
and \ac{sluacs} pass from $X$ to $X/U$.
\end{proposition}

\begin{Proof}
Let $\omega: X \rightarrow X/U$ be the canonical quotient map. As was observed in the proof of 
\cite{klee}*{Proposition 3.2} the reflexivity of $U$ implies $\omega(B_X)=B_{X/U}$.\par
Now suppose that $X$ is \ac{sluacs} and take a sequence $(z_n)_{n\in \N}$ in $S_{X/U}$ and 
an element $z\in S_{X/U}$ such that $\norm{z_n+z}\to 2$. Further, take a sequence 
$(\psi_n)_{n\in \N}$ in $S_{(X/U)^*}$ with $\psi_n(z_n)\to 1$.\par
Since $\omega(B_X)=B_{X/U}$ we can find a sequence $(x_n)_{n\in \N}$ in $S_X$ and a point 
$x\in S_X$ such that $z_n=\omega(x_n)$ for every $n$ and $z=\omega(x)$.\par
It easily follows from $\norm{z_n+z}\to 2$ that we also have $\norm{x_n+x}\to 2$.\par
We put $x_n^*:=\psi_n\circ\omega\in S_{U^\perp}$ for every $n$ and observe that $x_n^*(x_n)=
\psi_n(z_n)\to 1$. Since $X$ is \ac{sluacs} this implies $x_n^*(x)=\psi_n(z)\to 1$.\par
The proofs for \ac{acs} and \ac{luacs} spaces are analogous.
\end{Proof}

Using again the relation $(X/U)^*\cong U^\perp$ for every closed subspace $U$ of $X$ we can derive 
the following from Propositions \ref{prop:dual acs} and \ref{prop:dual}.

\begin{proposition}\label{prop:dual, stronger version}
If $U$ is a closed subspace of the Banach space $X$ the following implications hold.
  \begin{enumerate}[\upshape(i)]
  \item $X^* \ \mathrm{\ac{acs}} \ \Rightarrow \ X/U \ \mathrm{\ac{acs}}$
  \item $X^* \ \mathrm{\ac{luacs}} \ \Rightarrow \ X/U \ \mathrm{\ac{luacs}}$
  \item $X^* \ \mathrm{\ac{wuacs}} \ \Rightarrow \ X/U \ \mathrm{\ac{sluacs}}$
  \end{enumerate}
\end{proposition}
It is known (cf. \cite{day5}*{p.145}) that for any Banach space $X$ the dual $X^*$ is \ac{R} (resp. \ac{S}) 
\ifif every quotient space of $X$ is \ac{S} (resp. \ac{R}) \ifif every two-dimensional quotient space of $X$ 
is \ac{S} (resp. \ac{R}). By an analogous argument we can get the following result.
\begin{proposition}\label{prop:dual acs 2-dim quotients}
For a Banach space $X$ the following assertions are equivalent.
  \begin{enumerate}[\upshape(i)]
  \item $X^*$ is \ac{acs}.
  \item $X/U$ is \ac{acs} for every closed subspace $U$ of $X$.
  \item $X/U$ is \ac{acs} for every closed subspace $U$ of $X$ with $\dim{X/U}=2$.
  \end{enumerate}
\end{proposition}

\begin{Proof}
$\mathrm{(i)} \Rightarrow \mathrm{(ii)}$ holds according to Proposition \ref{prop:dual, stronger version} 
and $\mathrm{(ii)} \Rightarrow \mathrm{(iii)}$ is trivial, so it only remains to prove 
$\mathrm{(iii)} \Rightarrow \mathrm{(i)}$. Obviously it suffices to show that every two-dimensional subspace
of $X^*$ is \ac{acs}, so let us take such a subspace $V=\lin{\set{x^*,y^*}}$. Then $V=U^\perp=(X/U)^*$, where 
$U=\ker{x^*}\cap\ker{y^*}$. The quotient space $X/U$ is two-dimensional and hence by our assumption it is 
\ac{acs}. Since $X/U$ is in particular reflexive it follows from Proposition \ref{prop:dual} that 
$(X/U)^*=V$ is also \ac{acs}.
\end{Proof}

By \cite{klee}*{Proposition 3.4} there is an equivalent norm $\Norm{\,.\,}$ on $\ell^1$ such that $(\ell^1,\Norm{\,.\,})$ 
is \ac{R} and every separable Banach space is isometrically isomorphic to a quotient space of $(\ell^1,\Norm{\,.\,})$, 
so in particular $\ell^1$ is a quotient of $(\ell^1,\Norm{\,.\,})$. Thus quotients of \ac{acs} spaces are in general not \ac{acs}
and it also follows (in view of Proposition \ref{prop:dual acs 2-dim quotients}) that the fact that $X$ is \ac{acs} is not
sufficient to ensure that $X^*$ is \ac{acs}.\par
There is also an analogue of Proposition \ref{prop:dual acs 2-dim quotients} for \ac{uacs} spaces which reads as follows.
(The corresponding result for \ac{UR} spaces was proved by Day (cf. \cite{day-factorspaces}*{Theorem 5.5}).)
\begin{proposition}\label{prop:dual uacs 2-dim quotients}
For a Banach space $X$ let $\operatorname{\mathcal{S}}(X)$ denote the set of all closed subspaces of $X$ and
$\operatorname{\mathcal{S}_2}(X)$ the set of all closed subspaces $U$ of $X$ such that $\dim{X/U}\leq 2$. Then the 
following assertions are equivalent:
  \begin{enumerate}[\upshape(i)]
  \item $X$ is \ac{uacs}.
  \item $\inf\set*{\delta_{\mathrm{uacs}}^{X/U}(\eps): U\in \operatorname{\mathcal{S}}(X)}>0 \ \ \forall \eps\in ]0,2]$.
  \item $\inf\set*{\delta_{\mathrm{uacs}}^{X/U}(\eps): U\in \operatorname{\mathcal{S}_2}(X)}>0 \ \ \forall \eps\in ]0,2]$.
  \end{enumerate}
\end{proposition}

\begin{Proof}
$\mathrm{(i)} \Rightarrow \mathrm{(ii)}$ Let $X$ be \ac{uacs}. If $U\in \operatorname{\mathcal{S}}(X)$ then $(X/U)^*\cong U^\perp$, hence
$\delta_{\mathrm{uacs}}^{(X/U)^*}(\eps)\geq \delta_{\mathrm{uacs}}^{X^*}(\eps)\geq \delta_{\mathrm{uacs}}^X\paren*{\delta_{\mathrm{uacs}}^X(\eps)}$
by Proposition \ref{prop:uacs-dual} and the reflexivity of $X$.\par
Using again Proposition \ref{prop:uacs-dual} (now applied to $X/U$) and the monotonicity of the \ac{uacs} modulus we obtain 
\begin{equation*}
\delta_{\mathrm{uacs}}^{X/U}(\eps)\geq \delta_{\mathrm{uacs}}^X\paren*{\delta_{\mathrm{uacs}}^X\paren*{\delta_{\mathrm{uacs}}^X\paren*{\delta_{\mathrm{uacs}}^X(\eps)}}}>0,
\end{equation*}
which finishes our argument.\par
Since $\mathrm{(ii)} \Rightarrow \mathrm{(iii)}$ is obvious it only remains to prove $\mathrm{(iii)} \Rightarrow \mathrm{(i)}$. Denote the infimum in (iii)
by $\delta(\eps)$ and take sequence $(x_n^*)_{n\in \N}, (y_n^*)_{n\in \N}$ in $S_{X^*}$  such that $\norm{x_n^*+y_n^*}\to 2$ and a sequence $(x_n^{**})_{n\in \N}$
in $S_{X^{**}}$ with $x_n^{**}(x_n^*)\to 1$.\par
We put $V_n=\lin\set*{x_n^*,y_n^*}$ and $U_n=\ker{x_n^*}\cap\ker{y_n^*}$ for every $n$. Then $V_n=U_n^\perp=(X/U_n)^*$. Again by Proposition \ref{prop:uacs-dual}
(and reflexivity of $X/U_n$) we get that $\delta_{\mathrm{uacs}}^{V_n}(\eps)\geq \delta_{\mathrm{uacs}}^{X/U_n}\paren*{\delta_{\mathrm{ucas}}^{X/U_n}(\eps)}\geq 
\delta\paren*{\delta(\eps)}$.\par
Let $\varphi_n$ denote the restriction of $x_n^{**}$ to $V_n$ and fix any $\eps_0>0$. Because of $\norm{x_n^*+y_n^*}\to 2$ we have 
$1-2^{-1}\norm{x_n^*+y_n^*}<\delta\paren*{\delta(\eps_0)}\leq \delta_{\mathrm{uacs}}^{V_n}(\eps_0)$ for sufficiently large $n$.\par
Since $\varphi_n(x_n^*)=1$ this implies that we eventually have $\varphi_n(y_n^*)=x_n^{**}(y_n^*)\geq 1-\eps_0$.\par
Thus we have shown that $X^*$ is \ac{uacs} and by Proposition \ref{prop:uacs-dual} $X$ is \ac{uacs} as well.
\end{Proof}

In the next section we will study absolute sums of \ac{uacs} spaces and their relatives, but first
we have to introduce two more definitions that will be needed, namely a kind of symmetrised versions
of the notions of \ac{luacs} and \ac{sluacs} spaces.
\begin{definition}\label{def:luacs+sluacs+}
A Banach space $X$ is called 
  \begin{enumerate}[(i)]
  \item a $\mathrm{luacs}^+$ space if for every $x\in S_X$, every sequence $(x_n)_{n\in \N}$ in $S_X$
  with $\norm*{x_n+x}\to 2$ and all $x^*\in S_{X^*}$ we have
  \begin{equation*}
  x^*(x_n)\to 1 \ \iff \ x^*(x)=1,
  \end{equation*}
  \item a $\mathrm{sluacs}^+$ space if for every $x\in S_X$, every sequence $(x_n)_{n\in \N}$ in $S_X$
  with $\norm*{x_n+x}\to 2$ and all sequences $(x_n^*)_{n\in \N}$ in $S_{X^*}$ we have
  \begin{equation*}
  x_n^*(x_n)\to 1 \ \iff \ x_n^*(x)\to 1.
  \end{equation*}
  \end{enumerate}
\end{definition}
\noindent If we include these two properties in our implication chart we get the following.

\begin{figure}[H]
\begin{center}
\begin{tikzpicture}
  \node (UR) at (-2,0) {UR};
  \node (WUR) at (0,1) {WUR};
  \node (LUR) at (0,-1) {LUR};
  \node (WLUR) at (2,0) {WLUR};
  \node (R) at (4,0) {R};
  \node (uacs) at (-2,-1) {uacs};
  \node (wuacs) at (0,0) {wuacs};
  \node (sluacs+) at (0,-2) {$\mathrm{sluacs}^+$};
  \node (luacs+) at (2,-1) {$\mathrm{luacs}^+$};
  \node (acs) at (4,-2) {acs};
  \node (sluacs) at (0,-3) {sluacs};
  \node (luacs) at (2,-2) {luacs};
  \draw[->] (UR)--(LUR);
  \draw[->] (UR)--(WUR);
  \draw[->] (LUR)--(WLUR);
  \draw[->] (WUR)--(WLUR);
  \draw[->] (WLUR)--(R);
  \draw[->] (uacs)--(sluacs+);
  \draw[->] (uacs)--(wuacs);
  \draw[->] (sluacs+)--(luacs+);
  \draw[->] (wuacs)--(luacs+);
  \draw[->] (sluacs+)--(sluacs);
  \draw[->] (luacs+)--(luacs);
  \draw[->] (sluacs)--(luacs);
  \draw[->] (luacs)--(acs);
  \draw[->] (UR)--(uacs);
  \draw[->] (LUR)--(sluacs+);
  \draw[->] (WUR)--(wuacs);
  \draw[->] (WLUR)--(luacs+);
  \draw[->] (R)--(acs);
\end{tikzpicture}
\end{center}
\CAP\label{fig:6}
\end{figure}

Let us mention that Proposition \ref{prop:quot acs} also holds for $\mathrm{luacs}^+$ 
and $\mathrm{sluacs}^+$ spaces (with the same argument). Also, Propositions \ref{prop:very luacs}
resp. \ref{prop:very sluacs} hold accordingly for $\mathrm{luacs}^+$ resp. $\mathrm{sluacs}^+$ spaces.\par
In analogy to Proposition \ref{prop:char uacs without dual} one can prove that for any Banach space $X$ 
the following conditions are equivalent:
\begin{enumerate}[(i)]
\item For all sequences $(x_n)_{n\in \N}$ in $S_X$, $(x_n^*)_{n\in \N}$ in $S_{X^*}$ and every $x\in S_X$
with $\norm{x_n+x}\to 2$ and $x_n^*(x)\to 1$ one has $x_n^*(x_n)\to 1$.
\item For every $x\in S_X$ and every $\eps>0$ there exists a $\delta>0$ such that 
\begin{equation*}
\norm{x+ty}+\norm{x-ty}\leq 2+\eps t
\end{equation*}
whenever $t\in [0,\delta]$ and $y\in S_X$ with $\norm{x+y}\geq 2(1-t)$.
\end{enumerate}
In particular, every \ac{FS} space fulfils (i) and hence a space which is \ac{FS} and \ac{sluacs}
is $\mathrm{sluacs}^+$. In the context of \ac{FS} spaces we also have the following proposition.\par
\begin{proposition}\label{prop:X FS X* acs implies X luacs}
If $X$ is \ac{FS} and $X^*$ is \ac{acs} then $X$ is $\text{luacs}^+$. In particular, every
reflexive \ac{FS} space is $\text{luacs}^+$.
\end{proposition}

\begin{Proof}
By our previous considerations we only have to show that $X$ is \ac{luacs}.\par
Take a sequence $(x_n)_{n\in \N}$ in $S_X$ and a point $x\in S_X$ with $\norm{x_n+x}\to 2$ as 
well as a functional $x^*\in S_{X^*}$ with $x^*(x_n)\to 1$. Choose a sequence $(y_n^*)_{n\in \N}$
in $S_{X^*}$ such that $y_n^*(x_n+x)=\norm{x_n+x}$ for every $n\in \N$. It follows that $y_n^*(x_n)\to 1$ 
and $y_n^*(x)\to 1$.\par
Because of $\norm{y_n^*+x^*}\geq y_n^*(x_n)+x^*(x_n)$ for every $n$ it follows that $\norm{y_n^*+x^*}\to 2$.
If $y^*\in S_{X^*}$ is the Fr\'echet-derivative of $\|\,.\,\|$ at $x$ then $y_n^*(x)\to 1$ implies 
$\norm{y_n^*-y^*}\to 0$. Hence we get $\norm{x^*+y^*}=2$ and $y^*(x)=1$.\par
Since $X^*$ is \ac{acs} we can conclude that $x^*(x)=1$.
\end{Proof}

We conclude this section with a simple lemma that will be frequently used in the sequel.
It is the generalization of \cite{abramovich}*{Lemma 2.1} to sequences, while the proof 
remains virtually the same.
\begin{lemma}\label{aux lemma}
Let $(x_n)_{n\in \N}$ and $(y_n)_{n\in \N}$ be sequences in the (real or complex) normed 
space $X$ such that $\norm*{x_n+y_n}-\norm*{x_n}-\norm*{y_n}\to 0$.\par
Then for any two bounded sequences $(\alpha_n)_{n\in \N}$, $(\beta_n)_{n\in \N}$ 
of non-negative real numbers we also have 
$\norm*{\alpha_nx_n+\beta_ny_n}-\alpha_n\norm*{x_n}-\beta_n\norm*{y_n}\to 0$.
\end{lemma}

\begin{Proof}
Let $n\in \N$ be arbitrary. If $\alpha_n\geq \beta_n$ then
\begin{align*}
&\norm*{\alpha_nx_n+\beta_ny_n}\geq \alpha_n\norm*{x_n+y_n}-\paren*{\alpha_n-\beta_n}\norm*{y_n}\\
&=\alpha_n\paren*{\norm*{x_n+y_n}-\norm*{x_n}-\norm*{y_n}}+\alpha_n\norm*{x_n}+\beta_n\norm*{y_n}
\end{align*}
and hence
\begin{equation*}
\norm*{\alpha_nx_n+\beta_ny_n}-\alpha_n\norm*{x_n}-\beta_n\norm*{y_n}\geq
\alpha_n\paren*{\norm*{x_n+y_n}-\norm*{x_n}-\norm*{y_n}}.
\end{equation*}
Analogously one can show that
\begin{equation*}
\norm*{\alpha_nx_n+\beta_ny_n}-\alpha_n\norm*{x_n}-\beta_n\norm*{y_n}\geq
\beta_n\paren*{\norm*{x_n+y_n}-\norm*{x_n}-\norm*{y_n}}
\end{equation*}
if $\alpha_n<\beta_n$. Since $(\alpha_n)_{n\in \N}$ and $(\beta_n)_{n\in \N}$ are bounded
we obtain the desired conclusion.
\end{Proof}

\section{Absolute sums}\label{sec:abs sums}
We begin by recalling some preliminaries on absolute sums. Let $I$ be a non-empty set, $E$ a subspace
of $\R^I$ with $e_i\in E$ for all $i\in I$ and $\norm*{\,.\,}_E$ a complete norm on $E$ (here $e_i$
denotes the characteristic function of $\set*{i}$).\par
\noindent The norm $\norm*{\,.\,}_E$ is called {\em absolute} if the following holds
\begin{align*}
&(a_i)_{i\in I}\in E, \ (b_i)_{i\in I}\in \R^I \ \mathrm{and} \ \abs*{a_i}=\abs*{b_i} \ \forall i\in I \\
&\Rightarrow \ (b_i)_{i\in I}\in E \ \mathrm{and} \ \norm*{(a_i)_{i\in I}}_E=\norm*{(b_i)_{i\in I}}_E.
\end{align*}
The norm is called {\em normalised} if $\norm*{e_i}=1$ for every $i\in I$.\par
Standard examples of subspaces of $\R^I$ with absolute normalised norm are the spaces $\ell^p(I)$
for $1\leq p\leq \infty$.\par
We have the following important lemma on absolute normalised norms, whose proof can be found for example 
in \cite{lee}*{Remark 2.1}.
\begin{lemma}\label{lemma:abs norms}
Let $(E,\norm*{\,.\,}_E)$ be a subspace of $\R^I$ with an absolute normalised norm. Then the following is true.
\begin{align*}
&(a_i)_{i\in I}\in E, \ (b_i)_{i\in I}\in \R^I \ \mathrm{and} \ \abs*{b_i}\leq\abs*{a_i} \ \forall i\in I \\
&\Rightarrow \ (b_i)_{i\in I}\in E \ \mathrm{and} \ \norm*{(b_i)_{i\in I}}_E\leq\norm*{(a_i)_{i\in I}}_E.
\end{align*}
Furthermore, the inclusions  $\ell^1(I)\ssq E\ssq\ell^\infty(I)$ hold and the respective inclusion mappings
are contractive.
\end{lemma}
For a given subspace $(E,\norm*{\,.\,}_E)$ of $\R^I$ endowed with an absolute normalised norm we put
\begin{equation*}
E^\prime:=\set*{(a_i)_{i\in I}\in \R^I: \sup_{(b_i)_{i\in I}\in B_E}\sum_{i\in I}\abs*{a_ib_i}<\infty}.
\end{equation*}
It is easy to check that  $E^\prime$ is a subspace of $\R^I$ and that
\begin{equation*}
\norm*{(a_i)_{i\in I}}_{E^\prime}:=\smashoperator[l]{\sup_{(b_i)_{i\in I}\in B_E}}\sum_{i\in I}\abs*{a_ib_i} 
\ \ \forall (a_i)_{i\in I}\in E^\prime
\end{equation*}
defines an absolute normalised norm on $E^\prime$.\par
The map $T: E^\prime \rightarrow E^*$ defined by
\begin{equation*}
T((a_i)_{i\in I})((b_i)_{i\in I}):=\sum_{i\in I}a_ib_i \ \ \forall (a_i)_{i\in I}\in E^\prime, \forall (b_i)_{i\in I}\in E
\end{equation*}
is easily seen to be an isometric embedding. Moreover, if $\lin\set*{e_i:i\in I}$ is dense in $E$ then $T$ is onto,
so in this case we can identify $E^*$ and $E^\prime$.\par
Now if $(X_i)_{i\in I}$ is a family of (real or complex) Banach spaces we put
\begin{equation*}
\Bigl[\bigoplus_{i\in I}X_i\Bigr]_E:=\set*{(x_i)_{i\in I}\in \prod_{i\in I}X_i: (\norm*{x_i})_{i\in I}\in E}.
\end{equation*}
It is not hard to see that this defines a subspace of the product space $\prod_{i\in I}X_i$ which becomes a 
Banach space when endowed with the norm 
\begin{equation*}
\norm*{(x_i)_{i\in I}}_E:=\norm*{(\norm*{x_i})_{i\in I}}_E \ \forall (x_i)_{i\in I}\in \Bigl[\bigoplus_{i\in I}X_i\Bigr]_E.
\end{equation*}
We call this Banach space the absolute sum of the family $(X_i)_{i\in I}$ with respect to $E$. Again, the map
\begin{align*}
&S: \Bigl[\bigoplus_{i\in I}X_i^*\Bigr]_{E^\prime} \rightarrow \Bigl[\bigoplus_{i\in I}X_i\Bigr]_E^* \\
&S((x_i^*)_{i\in I})((x_i)_{i\in I}):=\sum_{i\in I}x_i^*(x_i)
\end{align*}
is an isometric embedding and it is onto if $\lin\set*{e_i:i\in I}$ is dense in $E$.\par
We also mention the following well-known fact, which will be needed later.
\begin{lemma}\label{lemma:abs norms l1}
If $E$ is a subspace of $\R^I$ endowed with an absolute normalised norm and $\lin\set*{e_i:i\in I}$ is dense in $E$ 
then $E$ contains no isomorphic copy of $\ell^1$ \ifif $\lin\set*{e_i:i\in I}$ is dense in $E^\prime$.\Todo{Add a reference.}
\end{lemma}
If not otherwise stated, we shall henceforth assume $E$ to be a subspace of $\R^I$ with an absolute
normalised norm such that $\lin\set*{e_i:i\in I}$ is dense in $E$.\par
Now let us first have a look at absolute sums of \ac{acs} spaces.
\begin{proposition}\label{prop:sum acs}
If $(X_i)_{i\in I}$ is a family of \ac{acs} spaces and $E$ is \ac{acs} then 
$\bigl[\bigoplus_{i\in I}X\bigr]_E$ is also \ac{acs}.
\end{proposition}

\begin{Proof}
Let $x=(x_i)_{i\in I}$ and $y=(y_i)_{i\in I}$ be elements of the unit sphere of
$\bigl[\bigoplus_{i\in I}X_i\bigr]_E$ and $x^*=(x_i^*)_{i\in I}$ an element of the 
dual unit sphere such that $\norm{x+y}_E=2$ and $x^*(x)=1$. We then have
\begin{equation*}
1=x^*(x)=\sum_{i\in I}x_i^*(x_i)\leq\sum_{i\in I}\norm{x_i^*}\norm{x_i}\leq\norm{x^*}_{E^\prime}\norm{x}_E=1
\end{equation*}
and hence
\begin{equation}\label{eq:5}
x_i^*(x_i)=\norm{x_i^*}\norm{x_i} \ \forall i\in I \ \mathrm{and} \ \sum_{i\in I}\norm{x_i^*}\norm{x_i}=1.
\end{equation}
Moreover, by Lemma \ref{lemma:abs norms} we have
\begin{align*}
&2=\norm{x+y}_E=\norm*{(\norm{x_i+y_i})_{i\in I}}_E\leq\norm*{(\norm{x_i}+\norm{y_i})_{i\in I}}_E \\
&\leq\norm{x}_E+\norm{y}_E=2
\end{align*}
and thus
\begin{equation}\label{eq:6}
\norm*{(\norm{x_i}+\norm{y_i})_{i\in I}}_E=2.
\end{equation}
Since $E$ is \ac{acs} \eqref{eq:6} and the second part of \eqref{eq:5} imply that
\begin{equation}\label{eq:7}
\sum_{i\in I}\norm{x_i^*}\norm{y_i}=1.
\end{equation}
Another application of Lemma \ref{lemma:abs norms} shows
\begin{equation}\label{eq:8}
\norm*{(\norm{x_i+y_i}+\norm{x_i}+\norm{y_i})_{i\in I}}_E=4.
\end{equation}
Again, since $E$ is \ac{acs} we get from \eqref{eq:8}, \eqref{eq:7} and the second part of \eqref{eq:5} that
\begin{equation*}
\sum_{i\in I}\norm{x_i^*}\norm{x_i+y_i}=2
\end{equation*}
which together with \eqref{eq:5} and \eqref{eq:7} implies
\begin{equation}\label{eq:9}
\norm{x_i^*}\paren*{\norm{x_i}+\norm{y_i}-\norm{x_i+y_i}}=0 \ \ \forall i\in I.
\end{equation}
Next we claim that
\begin{equation}\label{eq:10}
x_i^*(y_i)=\norm{x_i^*}\norm{y_i} \ \ \forall i\in I.
\end{equation}
To see this, fix any $i_0\in I$ with $x_{i_0}^*\neq0$ and $y_{i_0}\neq0$. Define $a_i=\norm{x_i^*}$ for all
$i\in I\sm\set{i_0}$ and $a_{i_0}=0$. Then $(a_i)_{i\in I}\in B_{E^\prime}$, because of Lemma \ref{lemma:abs norms}.\par
If $x_{i_0}=0$ it would follow that $\sum_{i\in I}a_i\norm{x_i}=\sum_{i\in I}\norm{x_i^*}\norm{x_i}=1$
and hence (because of \eqref{eq:6} and since $E$ is \ac{acs}) we would also have $\sum_{i\in I}a_i\norm*{y_i}=1$.
But by \eqref{eq:7} this would imply $\norm{y_{i_0}}\norm{x_{i_0}^*}=\sum_{i\in I}\norm{y_i}\paren*{\norm{x_i^*}-a_i}=0$,
a contradiction.\par
Thus $x_{i_0}\neq0$. From \eqref{eq:9} and Lemma \ref{aux lemma} we get that 
\begin{equation*}
\norm*{\frac{x_{i_0}}{\norm{x_{i_0}}}+\frac{y_{i_0}}{\norm{y_{i_0}}}}=2.
\end{equation*}
Taking into account the first part of \eqref{eq:5} and the fact that $X_{i_0}$ is \ac{acs} we get 
$x_{i_0}^*(y_{i_0})=\norm{x_{i_0}^*}\norm{y_{i_0}}$, as desired.\par
Now from \eqref{eq:10} and \eqref{eq:7} it follows that $x^*(y)=1$ and we are done.
\end{Proof}

We remark that the special case of finitely many summands in the above proposition has already been treated in
\cite{dhompongsa} (for two summand) and \cite{mitani} (for finitely many summands) in the context of $u$-spaces 
and the so called $\psi$-direct sums.\par
Before we can get on, we have to introduce another technical definition.
\begin{definition}
The space $E$ is said to have the property $(P)$ if for every sequence $(a_n)_{n\in \N}$ in $S_E$
and every $a\in S_E$ we have
\begin{equation*}
\norm{a_n+a}_E\to 2 \ \Rightarrow \ a_n\to a \ \mathrm{pointwise}.
\end{equation*}
\end{definition}
\noindent If $E$ is \ac{WLUR} then it obviously has property $(P)$. The converse is true if $E$ contains
no isomorphic copy of $\ell^1$ by Lemma \ref{lemma:abs norms l1}.\par
With this notion we can formulate the following proposition.
\begin{proposition}\label{prop:sum sluacs (P)}
If $(X_i)_{i\in I}$ is a family of \ac{sluacs} (resp. \ac{luacs}) spaces and $E$ is \ac{sluacs} (resp. \ac{luacs})
and has the property $(P)$ then $\bigl[\bigoplus_{i\in I}X_i\bigr]_E$ is \ac{sluacs} (resp. \ac{luacs}) as well.
\end{proposition}

\begin{Proof}
We only prove the \ac{sluacs} case. The argument for \ac{luacs} spaces is analogous.\par
So let $(x_n)_{n\in \N}$ be a sequence in the unit sphere of $\bigl[\bigoplus_{i\in I}X_i\bigr]_E$ 
and $x=(x_i)_{i\in I}$ another element of norm one such that $\norm{x_n+x}_E\to 2$ and let
$(x_n^*)_{n\in \N}$ be a sequence in the dual unit sphere such that $x_n^*(x_n)\to 1$.\par
Write $x_n=(x_{n,i})_{i\in I}$ and $x_n^*=(x_{n,i}^*)_{i\in I}$ for each $n$. We then have
\begin{equation*}
x_n^*(x_n)=\sum_{i\in I}x_{n,i}^*(x_{n,i})\leq\sum_{i\in I}\norm{x_{n,i}^*}\norm{x_{n,i}}\leq\norm{x_n^*}_{E^\prime}\norm{x_n}_E=1
\end{equation*}\
which gives us
\begin{equation}\label{eq:11}
\lim_{n\to \infty}\sum_{i\in I}\norm{x_{n,i}^*}\norm{x_{n,i}}=1
\end{equation}
and
\begin{equation}\label{eq:12}
\lim_{n\to \infty}\paren*{x_{n,i}^*(x_{n,i})-\norm{x_{n,i}^*}\norm{x_{n,i}}}=0 \ \ \forall i\in I.
\end{equation}
Applying Lemma \ref{lemma:abs norms} we also get
\begin{equation*}
\norm*{x_n+x}_E\leq\norm*{(\norm{x_{n,i}}+\norm{x_i})_{i\in I}}_E\leq\norm{x_n}_E+\norm{x}_E=2
\end{equation*}
and hence
\begin{equation}\label{eq:13}
\lim_{n\to \infty}\norm*{(\norm{x_{n,i}}+\norm{x_i})_{i\in I}}_E=2.
\end{equation}
Since $E$ has property $(P)$ this implies
\begin{equation}\label{eq:14}
\lim_{n\to \infty}\norm*{x_{n,i}}=\norm*{x_i} \ \ \forall i\in I.
\end{equation}
Because $E$ is \ac{sluacs} we get from \eqref{eq:11} and \eqref{eq:13} that
\begin{equation}\label{eq:15}
\lim_{n\to \infty}\sum_{i\in I}\norm{x_{n,i}^*}\norm{x_i}=1.
\end{equation}
If we apply Lemma \ref{lemma:abs norms} again we arrive at
\begin{equation}\label{eq:16}
\lim_{n\to \infty}\norm*{(\norm{x_{n,i}+x_i}+\norm{x_{n,i}}+\norm{x_i})_{i\in I}}_E=4.
\end{equation}
We further have
\begin{align*}
&\norm{x_n+x}_E+1\geq\norm*{(\norm{x_{n,i}+x_i}+\norm{x_{n,i}})_{i\in I}}_E \\
&\geq\norm*{(\norm{x_{n,i}+x_i}+\norm{x_{n,i}}+\norm{x_i})_{i\in I}}_E-1
\end{align*}
and thus it follows from \eqref{eq:16} that
\begin{equation}\label{eq:17}
\lim_{n\to \infty}\norm*{(\norm{x_{n,i}+x_i}+\norm{x_{n,i}})_{i\in I}}_E=3.
\end{equation}
Analogously one can shown
\begin{equation}\label{eq:18}
\lim_{n\to \infty}\norm*{(\norm{x_{n,i}+x_i}+\norm{x_i})_{i\in I}}_E=3.
\end{equation}
But because of Lemma \ref{lemma:abs norms} we also have
\begin{align*}
&\norm*{(\norm{x_{n,i}+x_i}+\norm{x_{n,i}})_{i\in I}}_E+3 \\
&\geq\norm*{(\norm{x_{n,i}+x_i}+\norm{x_{n,i}}+3\norm{x_i})_{i\in I}}_E \\
&\geq2\norm*{(\norm{x_{n,i}+x_i}+\norm{x_i})_{i\in I}}_E
\end{align*}
and thus \eqref{eq:18} and \eqref{eq:17} imply
\begin{equation}\label{eq:19}
\lim_{n\to \infty}\norm*{(\norm{x_{n,i}+x_i}+\norm{x_{n,i}}+3\norm{x_i})_{i\in I}}_E=6.
\end{equation}
Since $E$ has property $(P)$ it follows from \eqref{eq:17} and \eqref{eq:19} (and some standard
normalisation arguments) that
\begin{equation*}
\lim_{n\to \infty}\paren*{\norm{x_{n,i}+x_i}+\norm{x_{n,i}}}=3\norm{x_i} \ \ \forall i\in I
\end{equation*}
which together with  \eqref{eq:14} gives us
\begin{equation}\label{eq:20}
\lim_{n\to \infty}\norm{x_{n,i}+x_i}=2\norm{x_i} \ \ \forall i\in I.
\end{equation}
Because each $X_i$ is \ac{sluacs} it follows from \eqref{eq:14}, \eqref{eq:20} and \eqref{eq:12}
(and again some standard normalisation arguments) that
\begin{equation}\label{eq:21}
\lim_{n\to \infty}\paren*{x_{n,i}^*(x_i)-\norm{x_{n,i}^*}\norm{x_i}}=0 \ \ \forall i\in I.
\end{equation}
Now take any $\eps>0$. Then there is a finite subset $J\ssq I$ such that
\begin{equation}\label{eq:22}
\norm*{\sum_{i\in J}\norm{x_i}e_i-(\norm{x_i})_{i\in I}}_E\leq\eps.
\end{equation}
By \eqref{eq:21} we can find an index $n_0\in \N$ such that
\begin{equation}\label{eq:23}
\abs*{\sum_{i\in J}\paren*{x_{n,i}^*(x_i)-\norm{x_{n,i}^*}\norm{x_i}}}\leq\eps \ \ \forall n\geq n_0.
\end{equation}
Then for all $n\geq n_0$ we have
\begin{align*}
&\abs*{x_n^*(x)-\sum_{i\in I}\norm{x_{n,i}^*}\norm{x_i}}
=\abs*{\sum_{i\in I}\paren*{x_{n,i}^*(x_i)-\norm{x_{n,i}^*}\norm{x_i}}} \\
&\leq\abs*{\sum_{i\in J}\paren*{x_{n,i}^*(x_i)-\norm{x_{n,i}^*}\norm{x_i}}}
+\abs*{\sum_{i\in I\sm J}\paren*{x_{n,i}^*(x_i)-\norm{x_{n,i}^*}\norm{x_i}}} \\
&\stackrel{\eqref{eq:23}}{\leq}\eps+2\sum_{i\in I\sm J}\norm{x_{n,i}^*}\norm{x_i}\leq
\eps+2\norm*{\sum_{i\in J}\norm{x_i}e_i-(\norm{x_i})_{i\in I}}_E\stackrel{\eqref{eq:22}}{\leq}3\eps.
\end{align*}
Thus we have shown $x_n^*(x)-\sum_{i\in I}\norm{x_{n,i}^*}\norm{x_i}\to 0$ which together
with \eqref{eq:15} leads to $x_n^*(x)\to 1$ finishing the proof.
\end{Proof}

In our next result we shall see that instead of supposing that $E$ possesses the property $(P)$ we can
also assume that $E$ is $\mathrm{sluacs}^+$ (resp. $\mathrm{luacs}^+$) to come to the same conclusion.
\begin{proposition}\label{prop:sum sluacs}
If $(X_i)_{i\in I}$ is a family of \ac{sluacs} (resp. \ac{luacs}) spaces and $E$ is $\text{sluacs}^+$
(resp. $\text{luacs}^+$) then $\bigl[\bigoplus_{i\in I}X_i\bigr]_E$ is also \ac{sluacs} (resp. \ac{luacs}).
\end{proposition}

\begin{Proof}
Again we only show the \ac{sluacs} case, the \ac{luacs} case being analogous.\par
So fix a sequence $(x_n)_{n\in \N}$, a point $x$ and a sequence $(x_n^*)_{n\in \N}$ 
of functionals just like in the proof of the preceding proposition.\par
As in this very proof we can show
\begin{equation}\label{eq:24}
\lim_{n\to \infty}\sum_{i\in I}\norm{x_{n,i}^*}\norm{x_{n,i}}=1
\end{equation}
and
\begin{equation}\label{eq:25}
\lim_{n\to \infty}\paren*{x_{n,i}^*(x_{n,i})-\norm{x_{n,i}^*}\norm{x_{n,i}}}=0 \ \ \forall i\in I
\end{equation}
as well as
\begin{equation}\label{eq:26}
\lim_{n\to \infty}\norm*{(\norm{x_{n,i}}+\norm{x_i})_{i\in I}}_E=2
\end{equation}
and
\begin{equation}\label{eq:27}
\lim_{n\to \infty}\sum_{i\in I}\norm{x_{n,i}^*}\norm{x_i}=1.
\end{equation}
Also as in the proof of Proposition \ref{prop:sum sluacs (P)} we can see
\begin{equation}\label{eq:28}
\lim_{n\to \infty}\norm*{(\norm{x_{n,i}+x_i}+\norm{x_{n,i}})_{i\in I}}_E=3
\end{equation}
and
\begin{equation}\label{eq:29}
\lim_{n\to \infty}\norm*{(\norm{x_{n,i}+x_i}+\norm{x_{n,i}}+3\norm{x_i})_{i\in I}}_E=6.
\end{equation}
Since $E$ is $\mathrm{sluacs}^+$ it follows from \eqref{eq:28}, \eqref{eq:29} and 
\eqref{eq:27} (with the usual normalisation arguments) that
\begin{equation*}
\lim_{n\to \infty}\sum_{i\in I}\norm{x_{n,i}^*}\paren*{\norm{x_{n,i}+x_i}+\norm{x_{n,i}}}=3.
\end{equation*}
Together with \eqref{eq:24} we get
\begin{equation*}
\lim_{n\to \infty}\sum_{i\in I}\norm{x_{n,i}^*}\paren*{\norm{x_{n,i}+x_i}-\norm{x_{n,i}}-\norm{x_i}}=0
\end{equation*}
and hence
\begin{equation}\label{eq:30}
\lim_{n\to \infty}\norm{x_{n,i}^*}\paren*{\norm{x_{n,i}+x_i}-\norm{x_{n,i}}-\norm{x_i}}=0 \ \ \forall i\in I.
\end{equation}
Next we show that
\begin{equation}\label{eq:31}
\lim_{n\to \infty}\paren*{x_{n,i}^*(x_i)-\norm{x_{n,i}^*}\norm{x_i}}=0 \ \ \forall i\in I.
\end{equation}
To see this we fix $i_0\in I$ with $x_{i_0}\neq0$. If $\norm{x_{n,i_0}^*}\to 0$ the statement 
is clear. Otherwise there is some $\eps>0$ such that $\norm{x_{n,i_0}^*}\geq\eps$ for infinitely
many $n$. Without loss of generality we may assume that this inequality holds for every $n\in \N$.\par
For each $n\in \N$ we put $a_{n,i}=\norm{x_{n,i}^*}$ for $i\in I\sm \set{i_0}$ and $a_{n,i_0}=0$.
Then $(a_{n,i})_{i\in I}\in B_{E^\prime}$ for every $n$.\par
It is $\abs{\sum_{i\in I}\paren{a_{n,i}-\norm{x_{n,i}^*}}\norm{x_{n,i}}}=
\norm{x_{n,i_0}^*}\norm{x_{n,i_0}}\leq\norm{x_{n,i_0}}$.\par
So if $\norm{x_{n,i_0}}\to 0$ then by \eqref{eq:24} we would also have 
$\lim_{n\to \infty}\sum_{i\in I}a_{n,i}\norm{x_{n,i}}=1$.\par
But since $E$ is a $\mathrm{sluacs}^+$  space this together with \eqref{eq:26} would also imply 
$\lim_{n\to \infty}\sum_{i\in I}a_{n,i}\norm{x_i}=1$, which in turn implies (because of \eqref{eq:27}) 
$\norm{x_{n,i_0}^*}\norm{x_{i_0}}=\abs{\sum_{i\in I}\paren{a_{n,i}-\norm{x_{n,i}^*}}\norm{x_i}}\to 0$, where
on the other hand $\norm{x_{n,i_0}^*}\norm{x_{i_0}}\geq\eps\norm{x_{i_0}}>0$ for all $n\in \N$, a contradiction.\par
So we must have $\norm{x_{n,i_0}}\not\to 0$ and hence there is some $\delta>0$ such that $\norm{x_{n,i_0}}\geq\delta$
for infinitely many (say for all) $n\in \N$.\par
Now since $(\norm{x_{n,i_0}^*})_{n\in \N}$ is bounded away from zero \eqref{eq:30} gives us that 
$\lim_{n\to \infty}\paren*{\norm{x_{n,i_0}+x_{i_0}}-\norm{x_{n,i_0}}-\norm{x_{i_0}}}=0$.\par
Because $(\norm{x_{n,i_0}})_{n\in \N}$ is bounded away from zero as well this together with Lemma \ref{aux lemma}
tells us that 
\begin{equation*}
\lim_{n\to \infty}\norm*{\frac{x_{n,i_0}}{\norm{x_{n,i_0}}}+\frac{x_{i_0}}{\norm{x_{i_0}}}}=2.
\end{equation*}
Using \eqref{eq:25} and the fact that $X_{i_0}$ is \ac{sluacs} we now get the desired conclusion.\par
Now that we have established \eqref{eq:31}, the rest of the proof can be carried out exactly as in 
Proposition \ref{prop:sum sluacs (P)}.
\end{Proof}

The next two propositions deal with sums of $\mathrm{luacs}^+$ and $\mathrm{sluacs}^+$ spaces. 
\begin{proposition}\label{prop:sum luacs+ (P)}
If $(X_i)_{i\in I}$ is a family of $\text{luacs}^+$ spaces and $E$ is $\text{luacs}^+$ and has
the property $(P)$ then $\bigl[\bigoplus_{i\in I}X_i\bigr]_E$ is also a $\text{luacs}^+$ space.
\end{proposition}

\begin{Proof}
By Proposition \ref{prop:sum sluacs} (or Proposition \ref{prop:sum sluacs (P)}) we already know that
the space $\bigl[\bigoplus_{i\in I}X_i\bigr]_E$ is \ac{luacs}.\par
Now take a sequence $(x_n)_{n\in \N}$ and an element $x=(x_i)_{i\in I}$ in the unit sphere of 
$\bigl[\bigoplus_{i\in I}X_i\bigr]_E$ such that $\norm{x_n+x}\to 2$ and a functional $x^*=(x_i^*)_{i\in I}$
of norm one with $x^*(x)=1$. Write $x_n=(x_{n,i})_{i\in I}$ for al $n\in \N$.\par
As in the proof of Proposition \ref{prop:sum acs} it follows from $x^*(x)=1$ that
\begin{equation}\label{eq:32}
x_i^*(x_i)=\norm{x_i^*}\norm{x_i} \ \forall i\in I \ \mathrm{and} \ \sum_{i\in I}\norm{x_i^*}\norm{x_i}=1
\end{equation}
and as in the proof of Proposition \ref{prop:sum sluacs (P)} one can show that
\begin{equation}\label{eq:33}
\lim_{n\to \infty}\norm*{(\norm{x_{n,i}}+\norm{x_i})_{i\in I}}_E=2.
\end{equation}
Since $E$ is $\mathrm{luacs}^+$ it follows from \eqref{eq:33} and the second part of \eqref{eq:32} that
we also have
\begin{equation}\label{eq:34}
\lim_{n\to \infty}\sum_{i\in I}\norm{x_i^*}\norm{x_{n,i}}=1.
\end{equation}
Because $E$ has property $(P)$ it also follows from \eqref{eq:33} that
\begin{equation}\label{eq:35}
\lim_{n\to \infty}\norm{x_{n,i}}=\norm{x_i} \ \ \forall i\in I.
\end{equation}
Exactly as in the proof of Proposition \ref{prop:sum sluacs (P)} we can see 
\begin{equation}\label{eq:36}
\lim_{n\to \infty}\norm{x_{n,i}+x_i}=2\norm{x_i} \ \ \forall i\in I.
\end{equation}
Since each $X_i$ is $\mathrm{luacs}^+$ we infer from \eqref{eq:36}, \eqref{eq:35} and the first part
of \eqref{eq:32} that 
\begin{equation}\label{eq:37}
\lim_{n\to \infty}x_i^*(x_{n,i})=\norm{x_i^*}\norm{x_i} \ \ \forall i\in I.
\end{equation}
Now take an arbitrary $\eps>0$ and fix a finite subset $J\ssq I$ such that
\begin{equation}\label{eq:38}
\norm*{\sum_{i\in J}\norm{x_i}e_i-(\norm{x_i})_{i\in I}}_E\leq\eps.
\end{equation}
From \eqref{eq:32}, \eqref{eq:34} and \eqref{eq:35} it follows that
\begin{equation*}
\lim_{n\to \infty}\sum_{i\in I\sm J}\norm{x_i^*}\norm{x_{n,i}}=\sum_{i\in I\sm J}\norm{x_i^*}\norm{x_i}
\end{equation*}
and by \eqref{eq:37} we also have 
\begin{equation*}
\lim_{n\to \infty}\sum_{i\in J}x_i^*(x_{n,i})=\sum_{i\in J}\norm{x_i^*}\norm{x_i}.
\end{equation*}
Hence there is some $n_0\in \N$ such that
\begin{align}\label{eq:39}
&\abs*{\sum_{i\in J}\paren*{x_i^*(x_{n,i})-\norm{x_i^*}\norm{x_i}}}\leq\eps \ \ \mathrm{and} \\ 
\label{eq:40}
&\abs*{\sum_{i\in I\sm J}\norm{x_i^*}\paren*{\norm{x_{n,i}}-\norm{x_i}}}\leq\eps \ \ \forall n\geq n_0. 
\end{align}
But then we have for every $n\geq n_0$
\begin{align*}
&\abs*{x^*(x_n)-1}\stackrel{\eqref{eq:32}}{=}\abs*{\sum_{i\in I}\paren*{x_i^*(x_{n,i})-\norm{x_i^*}\norm{x_i}}} \\
&\stackrel{\eqref{eq:39}}{\leq}\eps+\abs*{\sum_{i\in I\sm J}\paren*{x_i^*(x_{n,i})-\norm{x_i^*}\norm{x_i}}} \\
&\ \ \leq\eps+\sum_{i\in I\sm J}\norm{x_i^*}\paren*{\norm{x_{n,i}}+\norm{x_i}} \\ 
&\stackrel{\eqref{eq:40}}{\leq}2\eps+2\sum_{i\in I\sm J}\norm{x_i^*}\norm{x_i}\stackrel{\eqref{eq:38}}{\leq}4\eps.
\end{align*}
Thus we have $x^*(x_n)\to 1$ and the proof is finished.
\end{Proof}

\begin{proposition}\label{prop:sum sluacs+}
If $(X_i)_{i\in I}$ is a family of $\text{sluacs}^+$ (resp. $\text{luacs}^+$) spaces and $E$ is $\text{sluacs}^+$ 
then $\bigl[\bigoplus_{i\in I}X_i\bigr]_E$ is $\text{sluacs}^+$ (resp. $\text{luacs}^+$) as well.
\end{proposition}

\begin{Proof}
Suppose all the $X_i$ and $E$ are $\mathrm{sluacs}^+$. Then by Proposition \ref{prop:sum sluacs} 
$\bigl[\bigoplus_{i\in I}X_i\bigr]_E$ is \ac{sluacs}.\par
Now take sequences $(x_n)_{n\in \N}$ and $(x_n^*)_{n\in \N}$ in the unit sphere and in the dual 
unit sphere of  $\bigl[\bigoplus_{i\in I}X_i\bigr]_E$ respectively, as well as another  element 
$x=(x_i)_{i\in I}$ in $\bigl[\bigoplus_{i\in I}X_i\bigr]_E$ of norm one such that $\norm{x_n+x}_E\to 2$
and $x_n^*(x)\to 1$.\par
As usual we write $x_n=(x_{n,i})_{i\in I}$ and $x_n^*=(x_{n,i}^*)_{i\in I}$ for every $n\in \N$.\par
Much as we have done before we can show that
\begin{equation}\label{eq:41}
\lim_{n\to \infty}\paren*{x_{n,i}^*(x_i)-\norm{x_{n,i}^*}\norm{x_i}}=0 \ \forall i\in I \ \mathrm{and} \ 
\lim_ {n\to \infty}\sum_{i\in I}\norm{x_{n,i}^*}\norm{x_i}=1
\end{equation}
as well as
\begin{equation}\label{eq:42}
\lim_{n\to \infty}\norm*{(\norm{x_{n,i}}+\norm{x_i})_{i\in I}}_E=2.
\end{equation}
It follows from \eqref{eq:42}, the second part of \eqref{eq:41}, and the fact that $E$ is $\mathrm{sluacs}^+$ that
\begin{equation}\label{eq:43}
\lim_{n\to \infty}\sum_{i\in I}\norm{x_{n,i}^*}\norm{x_{n,i}}=1.
\end{equation}
As in the proof of Proposition \ref{prop:sum sluacs} we see that
\begin{equation}\label{eq:44}
\lim_{n\to \infty}\norm{x_{n,i}^*}\paren*{\norm{x_{n,i}+x_i}-\norm{x_{n,i}}-\norm{x_i}}=0 \ \ \forall i\in I.
\end{equation}
Now using an argument analogous to that in the proof of Proposition \ref{prop:sum sluacs} shows
\begin{equation}\label{eq:45}
\lim_{n\to \infty}\paren*{x_{n,i}^*(x_{n,i})-\norm{x_{n,i}^*}\norm{x_{n,i}}}=0 \ \ \forall i\in I.
\end{equation}
Put $b_J=(\norm{x_i})_{i\in I}-\sum_{i\in J}\norm{x_i}e_i$ and $c_{n,J}=\sum_{i\in J}\norm{x_{n,i}^*}e_i$ for
every $n\in \N$ and every finite subset $J\ssq I$. Then for every $n$ and $J$ we have
\begin{align}
\nonumber&\abs*{c_{n,J}((\norm{x_i})_{i\in I})-1}=\abs*{\sum_{i\in J}\norm{x_{n,i}^*}\norm{x_i}-1} \\
\nonumber&\leq\abs*{\sum_{i\in I\sm J}\norm{x_{n,i}^*}\norm{x_i}}+\abs*{\sum_{i\in I}\norm{x_{n,i}^*}\norm{x_i}-1} \\
\label{eq:46}&\leq\norm{b_J}_E+\abs*{\sum_{i\in I}\norm{x_{n,i}^*}\norm{x_i}-1}.
\end{align}
Now take any $\eps>0$. Because $E$ is $\mathrm{sluacs}^+$ there is some $\delta>0$ such that 
\begin{align}
\nonumber&a\in S_E, \ g\in B_{E^*} \ \mathrm{with} \ \norm{a+(\norm{x_i})_{i\in I}}_E\geq2-\delta \\
\label{eq:47}&\mathrm{and} \ g((\norm{x_i})_{i\in I})\geq1-\delta \ \Rightarrow \ g(a)\geq1-\eps.
\end{align}
Fix a finite subset $J_0\ssq I$ such that $\norm{b_{J_0}}_E\leq\delta/2$ and also fix an index $n_0$ such that
$\abs*{\sum_{i\in I}\norm{x_{n,i}^*}\norm{x_i}-1}\leq\delta/2$ and $\norm*{(\norm{x_{n,i}}+\norm{x_i})_{i\in I}}_E\geq2-\delta$
for all $n\geq n_0$ (which is possible because of \eqref{eq:41} and \eqref{eq:42}).\par
Then \eqref{eq:46} and \eqref{eq:47} give us
\begin{equation}\label{eq:48}
c_{n,J_0}((\norm{x_{n,i}})_{i\in I})=\sum_{i\in J_0}\norm{x_{n,i}^*}\norm{x_{n,i}}\geq1-\eps \ \ \forall n\geq n_0.
\end{equation}
By \eqref{eq:45} we may also assume that
\begin{equation}\label{eq:49}
\abs*{\sum_{i\in J_0}\paren*{x_{n,i}^*(x_{n,i})-\norm{x_{n,i}^*}\norm{x_{n,i}}}}\leq\eps \ \ \forall n\geq n_0.
\end{equation}
Then for every $n\geq n_0$ we have 
\begin{align*}
&\abs*{x_n^*(x_n)-\sum_{i\in I}\norm{x_{n,i}^*}\norm{x_{n,i}}}=
\abs*{\sum_{i\in I}\paren*{x_{n,i}^*(x_{n,i})-\norm{x_{n,i}^*}\norm{x_{n,i}}}} \\
&\stackrel{\eqref{eq:49}}{\leq}\eps+\abs*{\sum_{i\in I\sm J_0}\paren*{x_{n,i}^*(x_{n,i})-\norm{x_{n,i}^*}\norm{x_{n,i}}}} \\
&\ \ \leq\eps+2\sum_{i\in I\sm J_0}\norm{x_{n,i}^*}\norm{x_{n,i}}\stackrel{\eqref{eq:48}}{\leq}3\eps.
\end{align*}
Thus $x_n^*(x_n)-\sum_{i\in I}\norm{x_{n,i}^*}\norm{x_{n,i}}\to 0$ which together with \eqref{eq:43} implies
$x_n^*(x_n)\to 1$.\par
The proof for the $\mathrm{luacs}^+$ case can be done in a very similar fashion.
\end{Proof}
In our next result we consider sums of \ac{wuacs} and $\mathrm{luacs}^+$ spaces for the case that $E$
does not contain $\ell^1$.
\begin{proposition}\label{prop:sum wuacs}
If $(X_i)_{i\in I}$ is a family of \ac{wuacs} (resp. $\text{luacs}^+$) spaces and if $E$ is 
\ac{wuacs} (resp. $\text{luacs}^+$) and does not contain an isomorphic copy of $\ell^1$ then 
$\bigl[\bigoplus_{i\in I}X_i\bigr]_E$ is also \ac{wuacs} (resp. $\text{luacs}^+$).
\end{proposition}

\begin{Proof}
Let us suppose that $E$ and all the $X_i$ are \ac{wuacs} and fix to sequences $(x_n)_{n\in \N}$
and $(y_n)_{n\in \N}$ in the unit sphere of $\bigl[\bigoplus_{i\in I}X_i\bigr]_E$ as well as a
norm one functional $x^*=(x_i^*)_{i\in I}$ on $\bigl[\bigoplus_{i\in I}X_i\bigr]_E$ such that 
$\norm{x_n+y_n}_E\to 2$ and $x^*(x_n)\to 1$. Write $x_n=(x_{n,i})_{i\in I}$ and $y_n=(y_{n,i})_{i\in I}$
for each $n$.\par
As we have often done before we deduce
\begin{equation}\label{eq:50}
\lim_{n\to \infty}\paren*{x_i^*(x_{n,i})-\norm{x_i^*}\norm{x_{n,i}}}=0 \ \forall i\in I \ \mathrm{and} \ 
\lim_{n\to \infty}\sum_{i\in I}\norm{x_i^*}\norm{x_{n,i}}=1
\end{equation}
and
\begin{equation}\label{eq:51}
\lim_{n\to \infty}\norm*{(\norm{x_{n,i}}+\norm{y_{n,i}})_{i\in I}}_E=2
\end{equation}
as well as
\begin{equation}\label{eq:52}
\lim_{n\to \infty}\norm*{(\norm{x_{n,i}+y_{n,i}}+\norm{x_{n,i}}+\norm{y_{n,i}})_{i\in I}}_E=4.
\end{equation}
Since $E$ is \ac{wuacs} \eqref{eq:51} and the second part of \eqref{eq:50} imply
\begin{equation}\label{eq:53}
\lim_{n\to \infty}\sum_{i\in I}\norm{x_i^*}\norm{y_{n,i}}=1.
\end{equation}
Applying again the fact that $E$ is \ac{wuacs} together with \eqref{eq:52}, \eqref{eq:53} and 
the second part of \eqref{eq:50} gives us 
\begin{equation*}
\lim_{n\to \infty}\sum_{i\in I}\norm{x_i^*}\paren*{\norm{x_{n,i}}+\norm{y_{n,i}}-\norm{x_{n,i}+y_{n,i}}}=0
\end{equation*}
and hence
\begin{equation}\label{eq:54}
\lim_{n\to \infty}\norm{x_i^*}\paren*{\norm{x_{n,i}}+\norm{y_{n,i}}-\norm{x_{n,i}+y_{n,i}}}=0 \ \ \forall i\in I.
\end{equation}
Now we can show
\begin{equation}\label{eq:55}
\lim_{n\to \infty}\paren*{x_i^*(y_{n,i})-\norm{x_i^*}\norm{y_{n,i}}}=0 \ \ \forall i\in I.
\end{equation}
The argument for this is similiar to what we have done before but we state it here for the sake
of completeness. Fix $i_0\in I$ with $x_{i_0}^*\neq0$ and $y_{n,i_0}\not\to 0$. Then there is 
$\tau>0$ such that $\norm{y_{n,i_0}}\geq\tau$ for infinitely many (without loss of generality for 
all) $n\in \N$.\par
Put $a_{i_0}=0$ and $a_i=\norm{x_i^*}$ for every $i\in I\sm \set*{i_0}$. If $\norm{x_{n,i_0}}\to 0$
then because of the second part of \eqref{eq:50} it would follow that 
$\lim_{n\to \infty}\sum_{i\in I}a_i\norm{x_{n,i}}=1$.\par
Since $E$ is \ac{wuacs} this together with \eqref{eq:51} would imply that we also have
$\lim_{n\to \infty}\sum_{i\in I}a_i\norm{y_{n,i}}=1$ which because \eqref{eq:53} would give us 
$\norm{x_{i_0}^*}\norm{y_{n,i_0}}\to 0$, a contradiction.\par
Hence there must be some $\delta>0$ such that $\norm{x_{n,i_0}}\geq\delta$ for infinitely many
(say for every) $n\in \N$.\par
Now since the sequences $(\norm{x_{n,i_0}})_{n\in \N}$ and $(\norm{y_{n,i_0}})_{n\in \N}$ are 
bounded away from zero it follows from \eqref{eq:50}, \eqref{eq:54} and Lemma \ref{aux lemma} that 
\begin{equation*}
\lim_{n\to \infty}\norm*{\frac{x_{n,i_0}}{\norm{x_{n,i_0}}}+\frac{y_{n,i_0}}{\norm{y_{n,i_0}}}}=2 \ 
\mathrm{and} \ \lim_{n\to \infty}\frac{x_{i_0}^*}{\norm{x_{i_0}^*}}\paren*{\frac{x_{n,i_0}}{\norm{x_{n,i_0}}}}=1.
\end{equation*}
Since $X_{i_0}$ is \ac{wuacs} this implies our desired conclusion.\par
Now we fix any $\eps>0$. Because $\ell^1\not\ssq E$ by Lemma \ref{lemma:abs norms l1} there must be some
finite set $J\ssq I$ such that
\begin{equation}\label{eq:56}
\norm*{(\norm{x_i^*})_{i\in I}-\sum_{i\in J}\norm{x_i^*}e_i}_{E^\prime}\leq\eps.
\end{equation}
By \eqref{eq:55} we can find some $n_0\in \N$ such that
\begin{equation}\label{eq:57}
\abs*{\sum_{i\in J}\paren*{x_i^*(y_{n,i})-\norm{x_i^*}\norm{y_{n,i}}}}\leq\eps \ \ \forall n\geq n_0.
\end{equation}
We then have for every $n\geq n_0$
\begin{align*}
&\abs*{x^*(y_n)-\sum_{i\in I}\norm{x_i^*}\norm{y_{n,i}}}\stackrel{\eqref{eq:57}}{\leq}
\eps+\abs*{\sum_{i\in I\sm J}\paren*{x_i^*(y_{n,i})-\norm{x_i^*}\norm{y_{n,i}}}} \\
&\leq\eps+2\sum_{i\in I\sm J}\norm{x_i^*}\norm{y_{n,i}}\stackrel{\eqref{eq:56}}{\leq}3\eps.
\end{align*}
So we have $x^*(y_n)-\sum_{i\in I}\norm{x_i^*}\norm{y_{n,i}^*}\to 0$. From \eqref{eq:53} it now 
follows that $x^*(y_n)\to 1$.\par
The $\mathrm{luacs}^+$ case is proved analogously.
\end{Proof}

Note that the above Proposition especially applies to the case that $E$ is \ac{WUR} because a
\ac{WUR} space cannot contain an isomorphic copy of $\ell^1$ (cf. \cite{zizler}*{Remark 4}).
Frankly, the author does not know whether a \ac{wuacs} space can contain an isomorphic copy of $\ell^1$
at all, but at least it cannot contain particularly `good' copies of $\ell^1$ in the following sense 
(introduced in \cite{dowling1}).

\begin{definition}\label{def:asymp l1}
A Banach space $X$ is said to contain an {\em asymptotically isometric copy of $\ell^1$} if
there is a sequence $(x_n)_{n\in \N}$ in $B_X$ and a decreasing sequence $(\eps_n)_{n\in \N}$
in $[0,1[$ with $\eps_n\to 0$ such that for each $m\in \N$ and all scalars 
$a_1,\dots,a_m$ we have
\begin{equation*}
\sum_{i=1}^m(1-\eps_i)\abs{a_i}\leq\norm*{\sum_{i=1}^ma_ix_i}\leq\sum_{i=1}^m\abs{a_i}.
\end{equation*}
Likewise, $X$ is said to contain an {\em asymptotically isomorphic copy of $c_0$} if there are two
such sequences $(x_n)_{n\in \N}$ and $(\eps_n)_{n\in \N}$ which fulfil
\begin{equation*}
\max_{i=1,\dots,m}(1-\eps_i)\abs{a_i}\leq\norm*{\sum_{i=1}^ma_ix_i}\leq\max_{i=1,\dots,m}\abs{a_i}
\end{equation*}
for each $m\in \N$ and all scalars $a_1,\dots,a_m$.
\end{definition}
We then have the following observation.
\begin{proposition}\label{prop:asymp l1 wuacs}
If the Banach space $X$ is \ac{wuacs} then it does not contain an asymptotically isometric 
copy of $\ell^1$.
\end{proposition}

\begin{Proof}
Suppose that $X$ contains an asymptotically isometric copy of $\ell^1$. Then fix two sequences 
$(x_n)_{n\in \N}$ and $(\eps_n)_{n\in \N}$ as in the above definition.\par
We can find $\alpha>1$ such that $\alpha\eps_n<1$ for every $n\in \N$. Put 
$\tilde{x}_n=(1-\alpha\eps_n)^{-1}x_n$ for each $n$. Then for every finite sequence 
$(a_i)_{i=1}^m$ of scalars we have
\begin{equation}\label{eq:58}
\norm*{\sum_{i=1}^ma_i\tilde{x}_i}=\norm*{\sum_{i=1}^m\frac{a_i}{1-\alpha\eps_i}x_i}\geq
\sum_{i=1}^m\frac{1-\eps_i}{1-\alpha\eps_i}\abs{a_i}\geq\sum_{i=1}^m\abs{a_i}.
\end{equation}
In other words, the operator $T: \ell^1 \rightarrow X$ defined by 
$T((a_n)_{n\in \N})=\sum_{n=1}^{\infty}a_n\tilde{x}_n$ is an isomorphism onto its range 
$U=\ran{T}$ with $\norm{T^{-1}}\leq1$.\par
Define $(b_n)_{n\in \N}\in \ell^{\infty}=(\ell^1)^*$ by $b_n=1$ if $n$ is even and $b_n=0$ 
if $n$ is odd. Then $u^*=(T^{-1})^*((b_n)_{n\in \N})\in B_{U^*}$. Take a Hahn-Banach extension
$x^*$ of $u^*$ to $X$.\par
Note that because of \eqref{eq:58} we have in particular $\norm{\tilde{x}_n}\geq1$ for every 
$n$ and on the other hand $\norm{\tilde{x}_n}\leq(1-\alpha\eps_n)^{-1}$ and $\eps_n\to 0$, 
hence $\norm{\tilde{x}_n}\to 1$. Again because of \eqref{eq:58} we have 
$\norm{\tilde{x}_n+\tilde{x}_{n+1}}\geq2$ for every $n$. It follows that 
$\norm{\tilde{x}_n+\tilde{x}_{n+1}}\to 2$ and thus in particular 
$\norm{\tilde{x}_{2n}+\tilde{x}_{2n+1}}\to 2$.\par
But we also have $x^*(\tilde{x}_{2n})=u^*(\tilde{x}_{2n})=b_{2n}=1$ and likewise 
$x^*(\tilde{x}_{2n+1})=b_{2n+1}=0$ for every $n$ and hence $X$ cannot be a \ac{wuacs} space.
\end{Proof}

If the space $X$ contains an asymptotically isometric copy of $c_0$ then by \cite{dowling1}*{Theorem 2}
$X^*$ contains an asymptotically isometric copy of $\ell^1$ and thus we get the following corollary.
\begin{corollary}\label{cor:asymp c0 wuacs}
If $X$ is a Banach space whose dual $X^*$ is \ac{wuacs} then $X$ does not contain an asymptotically
isometric copy of $c_0$.
\end{corollary}
We also remark that since $\ell^p(I)$ is \ac{UR} for every $1<p<\infty$ we can obtain the following
corollary from our above results.
\begin{corollary}\label{cor:p-sum}
If $(X_i)_{i\in I}$ is a family of Banach space such that each $X_i$ is \ac{acs} resp. \ac{luacs}
resp. $\text{luacs}^+$ resp. \ac{sluacs} resp. $\text{sluacs}^+$ resp. \ac{wuacs} then
$\bigl[\bigoplus_{i\in I}X_i\bigr]_p$ is also \ac{acs} resp. \ac{luacs} resp. $\text{luacs}^+$ 
resp. \ac{sluacs} resp. $\text{sluacs}^+$ resp. \ac{wuacs} for every $1<p<\infty$.
\end{corollary}
Now we turn to sums of \ac{uacs} spaces. We first consider sums of finitely many spaces.
In fact, this has been done before in \cite{dhompongsa} (for two summands) and in \cite{mitani} 
(for finitely many summands) in the context of $U$-spaces and the so called $\psi$-direct sums.
However, we include a sketch of our own slightly different proof here, for the sake of completeness.
\begin{proposition}\label{prop:finite sum uacs}
If $I$ is a finite set, $(X_i)_{i\in I}$ a family of \ac{uacs} Banach spaces and $\norm*{\,.\,}_E$
is an absolute normalized norm on $\R^I$ such that $E:=(\R^I,\norm*{\,.\,}_E)$ is \ac{acs} then
$\bigl[\bigoplus_{i\in I}X_i\bigr]_E$ is also a \ac{uacs} space.
\end{proposition}

\begin{Proof}
First note that since $E$ is finite-dimensional it is actually \ac{uacs}. Now if we take to sequences 
$(x_n)_{n\in \N}$ and $(y_n)_{n\in \N}$ in the unit sphere of $\bigl[\bigoplus_{i\in I}X_i\bigr]_E$
and a sequence $(x_n^*)_{n\in \N}$ in the dual unit sphere such that $\norm{x_n+y_n}_E\to 2$ and
$x_n^*(x_n)\to 1$ then we can show just as we have done before that
\begin{equation}\label{eq:59}
\lim_{n\to \infty}\sum_{i\in I}\norm{x_{n,i}^*}\norm{x_{n,i}}=1
\end{equation}
and
\begin{equation}\label{eq:60}
\lim_{n\to \infty}\paren*{x_{n,i}^*(x_{n,i})-\norm{x_{n,i}^*}\norm{x_{n,i}}}=0 \ \ \forall i\in I
\end{equation}
as well as
\begin{equation}\label{eq:61}
\lim_{n\to \infty}\norm*{(\norm{x_{n,i}}+\norm{y_{n,i}})_{i\in I}}_E=2
\end{equation}
and
\begin{equation}\label{eq:62}
\lim_{n\to \infty}\norm*{(\norm{x_{n,i}+y_{n,i}}+\norm{x_{n,i}}+\norm{y_{n,i}})_{i\in I}}_E=4.
\end{equation}
Since $E$ is \ac{uacs} it follows from \eqref{eq:59} and \eqref{eq:61} that 
\begin{equation}\label{eq:63}
\lim_{n\to \infty}\sum_{i\in I}\norm{x_{n,i}^*}\norm{y_{n,i}}=1.
\end{equation}
Again, since $E$ is \ac{uacs} it follows from \eqref{eq:59}, \eqref{eq:63} and \eqref{eq:62} that
\begin{equation}\label{eq:64}
\lim_{n\to \infty}\norm{x_{n,i}^*}\paren*{\norm{x_{n,i}}+\norm{y_{n,i}}-\norm{x_{n,i}+y_{n,i}}}=0 \ \ \forall i\in I.
\end{equation}
Now using \eqref{eq:64}, Lemma \ref{aux lemma}, \eqref{eq:59}, \eqref{eq:63}, \eqref{eq:60}, the fact that each
$X_i$ is \ac{uacs} and an argument similiar the one used in the proof of Proposition \ref{prop:sum wuacs} we
can infer that
\begin{equation*}
\lim_{n\to \infty}\paren*{x_{n,i}^*(y_{n,i})-\norm{x_{n,i}^*}\norm{y_{n,i}}}=0 \ \ \forall i\in I.
\end{equation*}
Since $I$ is finite it follows that $x_n^*(y_n)-\sum_{i\in I}\norm{x_{n,i}^*}\norm{y_{n,i}}\to 0$ which 
together with \eqref{eq:63} gives us $x_n^*(y_n)\to 1$ and the proof is over.
\end{Proof}

Before we can come to the study of absolute sums of infinitely many \ac{uacs} spaces we have to 
introduce one more definition.
\begin{definition}\label{def:u+}
The space $E$ is said to have the property $(u^+)$ if for every $\eps>0$ there is some $\delta>0$ such
that for all $(a_i)_{i\in I}, (b_i)_{i\in I}\in S_E$ and each $(c_i)_{i\in I}\in S_{E^\prime}=S_{E^*}$
we have
\begin{equation*}
\sum_{i\in I}a_ic_i=1 \ \mathrm{and} \ \norm*{(a_i+b_i)_{i\in I}}_E\geq2(1-\delta) \ \Rightarrow \ 
\sum_{i\in I}\abs{c_i}\abs{a_i-b_i}\leq\eps.
\end{equation*}
\end{definition}
\noindent Clearly, if $E$ is \ac{UR} then it has property $(u^+)$ and the property $(u^+)$ in turn 
implies that $E$ is \ac{uacs}.\par
Now we can formulate and prove the following theorem, which is an analogue of Day's results on
sums of \ac{UR} spaces from \cite{day2}*{Theorem 3} (for the $\ell^p$-case) and \cite{day3}*{Theorem 3}
(for the general case). Also, its proof is just a slight modification of Day's technique.
\begin{theorem}\label{thm:sum uacs}
If $(X_i)_{i\in I}$ is a family of Banach spaces such that for every $0<\eps\leq2$ we have 
$\delta(\eps):=\inf_{i\in I}\delta_{\mathrm{uacs}}^{X_i}(\eps)>0$ and if the space $E$ has  
the property $(u^+)$ then $\bigl[\bigoplus_{i\in I}X_i\bigr]_E$ is also \ac{uacs}.
\end{theorem}

\begin{Proof}
As in \cite{day2} and \cite{day3} the proof is divided into two steps. In the first step we 
show that for every $0<\eps\leq2$ there is some $\eta>0$ such that for any two elements 
$x=(x_i)_{i\in I}$ and $y=(y_i)_{i\in I}$ of the unit sphere of $\bigl[\bigoplus_{i\in I}X_i\bigr]_E$
with $\norm{x_i}=\norm{y_i}$ for every $i\in I$ and each functional $x^*=(x_i^*)_{i\in I}$ with 
$\norm{x^*}_{E^\prime}=x^*(x)=1$ and $x^*(y)<1-\eps$ we have $\norm{x+y}_E\leq2(1-\eta)$.\par
So let $0<\eps\leq2$ be arbitrary. Since $E$ is \ac{uacs} there exists some $\eta>0$ such that
\begin{align}
\nonumber
&a, b\in B_E, l\in B_{E^*}, l(a)=1 \ \mathrm{and} \ l(b)<1-\frac{\eps}{4}\delta\paren*{\frac{\eps}{2}} \\
\label{eq:65}
&\ \Rightarrow \ \norm{a+b}_E\leq2(1-\eta).
\end{align}
We claim that this $\eta$ fulfils our requirement. To show this, fix $x, y$ and $x^*$ as above and put
$\beta_i=\norm{x_i}=\norm{y_i}$, $\nu_i=\norm{x_i^*}$ and $\gamma_i=\nu_i\beta_i-x_i^*(y_i)$ for 
each $i\in I$. Then we have 
\begin{equation}\label{eq:66}
0\leq\gamma_i\leq2\beta_i\nu_i \ \ \forall i\in I.
\end{equation}
From $x^*(x)=1=\norm{x^*}_{E^\prime}=\norm{x}_E$ we get
\begin{equation}\label{eq:67}
\sum_{i\in I}\nu_i\beta_i=1 \ \mathrm{and} \ x_i^*(x_i)=\nu_i\beta_i \ \forall i\in I.
\end{equation}
Next we define
\begin{equation}\label{eq:68}
\alpha_i=
\begin{cases}
\frac{1}{2}\delta\paren*{\frac{\gamma_i}{\nu_i\beta_i}} & \mathrm{if} \ \gamma_i>0 \\
0 & \mathrm{if} \ \gamma_i=0.
\end{cases}
\end{equation}
From the definition of the $\delta_{\mathrm{uacs}}^{X_i}$ and the second part of \eqref{eq:67}
it easily follows that
\begin{equation}\label{eq:69}
\norm{x_i+y_i}\leq2(1-\alpha_i)\beta_i \ \ \forall i\in I.
\end{equation}
By \eqref{eq:66} and the first part of \eqref{eq:67} we have $\sum_{i\in I}\gamma_i\leq2$ and 
further it is
\begin{equation*}
\eps<1-x^*(y)=x^*(x-y)=\sum_{i\in I}x_i^*(x_i-y_i)\leq\sum_{i\in I}\gamma_i
\end{equation*}
thus
\begin{equation}\label{eq:70}
\eps<\sum_{i\in I}\gamma_i\leq2.
\end{equation}
Now put $A=\set*{i\in I:2\gamma_i>\eps\nu_i\beta_i}$ and $B=I\sm A$. Then we get
\begin{equation}\label{eq:71}
\sum_{i\in B}\gamma_i\leq\frac{\eps}{2}\sum_{i\in B}\nu_i\beta_i
\leq\frac{\eps}{2}\sum_{i\in I}\nu_i\beta_i\stackrel{\eqref{eq:67}}{=}\frac{\eps}{2}.
\end{equation}
From \eqref{eq:70} and \eqref{eq:71} it follows that
\begin{equation}\label{eq:72}
\sum_{i\in A}\gamma_i=\sum_{i\in I}\gamma_i-\sum_{i\in B}\gamma_i>\frac{\eps}{2}.
\end{equation}
Using \eqref{eq:66} and \eqref{eq:72} we now get
\begin{equation}\label{eq:73}
\sum_{i\in A}\nu_i\beta_i>\frac{\eps}{4}.
\end{equation}
Write $t=(\beta_i\chi_B(i))_{i\in I}$ and $t^{\prime}=(\beta_i\chi_A(i))_{i\in I}$, where $\chi_B$ 
and $\chi_A$ denote the characteristic function of $B$ and $A$ respectively. Then $t,t^{\prime}\in B_E$
(by Lemma \ref{lemma:abs norms}) and $t+t^{\prime}=(\beta_i)_{i\in I}$. We also put 
$t^{\prime\prime}=(1-\delta(\eps/2))t^{\prime}$. Again by Lemma \ref{lemma:abs norms} we have 
$\norm{t+t^{\prime\prime}}_E\leq\norm{t+t^{\prime}}_E=1$. Further, $l=(\nu_i)_{i\in I}$ defines an 
element of $S_{E^*}$ such that $l(t+t^{\prime})=\sum_{i\in I}\nu_i\beta_i=1$ (by \eqref{eq:67}) and 
$l(t+t^{\prime\prime})=1-\delta(\eps/2)l(t^{\prime})=1-\delta(\eps/2)\sum_{i\in A}\nu_i\beta_i$ and
hence (by \eqref{eq:73}) 
\begin{equation*}
l(t+t^{\prime\prime})<1-\frac{\eps}{4}\delta\paren*{\frac{\eps}{2}}.
\end{equation*}
Thus we can apply \eqref{eq:65} to deduce
\begin{equation}\label{eq:74}
\frac{1}{2}\norm*{2t+t^{\prime}+t^{\prime\prime}}_E=
\norm*{t+\paren*{1-\frac{1}{2}\delta\paren*{\frac{\eps}{2}}}t^{\prime}}_E\leq1-\eta.
\end{equation}
Since $\delta$ is obviously an increasing function we also have
\begin{equation}\label{eq:75}
\alpha_i\geq\frac{1}{2}\delta\paren*{\frac{\eps}{2}} \ \ \forall i\in A.
\end{equation}
Now we can conclude (with the aid of Lemma \ref{lemma:abs norms})
\begin{align*}
&\norm{x+y}_E=\norm*{(\norm{x_i+y_i})_{i\in I}}_E
\stackrel{\eqref{eq:69}}{\leq}2\norm*{((1-\alpha_i)\beta_i)_{i\in I}}_E \\
&\stackrel{\eqref{eq:75}}{\leq}2\norm*{\paren*{\paren*{1-\frac{1}{2}\delta\paren*{\frac{\eps}{2}}}\beta_i\chi_A(i)
+\beta_i\chi_B(i)}_{i\in I}}_E \\
&\ \ =2\norm*{\paren*{1-\frac{1}{2}\delta\paren*{\frac{\eps}{2}}}t^{\prime}+t}_E\stackrel{\eqref{eq:74}}{\leq}2(1-\eta),
\end{align*}
finishing the first step of the proof. Note that so far we have only used the fact that $E$ is \ac{uacs}
and not the property $(u^+)$.\par
Now for the second step we fix $0<\eps\leq2$  and choose an $\eta>0$ to the value $\eps/2$ according 
to step one. Then we take $0<\nu<2\eta/3$. Since $E$ is \ac{uacs} we can find $\tau>0$ such that
\begin{align}
\nonumber
& a, b\in B_E, l\in B_{E^*}, l(a)\geq1-\tau \ \mathrm{and} \ \norm{a+b}_E\geq2(1-\tau) \\
\label{eq:76}
& \ \Rightarrow \ l(b)\geq1-\nu.
\end{align}
Next we fix $0<\alpha<\min\set*{\eps/2,2\tau,\nu}$. Now we can find a number $\tilde{\tau}>0$ to the value 
$\alpha$ according to the definition of the property $(u^+)$ (Definition \ref{def:u+}). Finally, we take
$0<\xi<\min\set*{\tau,\tilde{\tau}}$.\par
Now suppose $x=(x_i)_{i\in I}$ and $y=(y_i)_{i\in I}$ are elements of the unit sphere of 
$\bigl[\bigoplus_{i\in I}X_i\bigr]_E$ and $x^*=(x_i^*)_{i\in I}$ is an element of the dual unit 
sphere such that $\norm{x+y}_E\geq2(1-\xi)$ and $x^*(x)=1$. We will show that $x^*(y)>1-\eps$.\par
To do so, we define
\begin{equation}\label{eq:77}
z_i=
\begin{cases}
\frac{\norm{x_i}}{\norm{y_i}}y_i & \mathrm{if} \ y_i\neq0 \\
x_i & \mathrm{if} \ y_i=0.
\end{cases}
\end{equation}
Then we have 
\begin{equation}\label{eq:78}
\norm{z_i}=\norm{x_i} \ \mathrm{and} \ \norm{z_i-y_i}=\abs*{\norm{x_i}-\norm{y_i}} \ \forall i\in I.
\end{equation}
As before we can see that $\sum_{i\in I}\norm{x_i^*}\norm{x_i}=1$ and further we have 
$2(1-\tilde{\tau})\leq2(1-\xi)\leq\norm{x+y}_E\leq\norm*{(\norm{x_i}+\norm{y_i})_{i\in I}}_E$.\par
Thus we get from the choice of $\tilde{{\tau}}$ that
\begin{equation}\label{eq:79}
\sum_{i\in I}\norm{x_i^*}\norm{z_i-y_i}
\stackrel{\eqref{eq:78}}{=}\sum_{i\in I}\norm{x_i^*}\abs*{\norm{x_i}-\norm{y_i}}\leq\alpha.
\end{equation}
Further, we have
\begin{equation*}
\norm*{(\norm{x_i}+\norm{y_i}+\norm{x_i+y_i})_{i\in I}}_E\geq2\norm{x+y}_E\geq4(1-\xi)\geq4(1-\tau)
\end{equation*}
and
\begin{align*}
& \sum_{i\in I}\norm{x_i^*}\paren*{\norm{x_i}+\norm{y_i}}=1+\sum_{i\in I}\norm{x_i^*}\norm{y_i} \\
& \geq1+\sum_{i\in I}\norm{x_i^*}\norm{x_i}-\sum_{i\in I}\norm{x_i^*}\abs*{\norm{x_i}-\norm{y_i}} \\
& =2-\sum_{i\in I}\norm{x_i^*}\abs*{\norm{x_i}-\norm{y_i}}\stackrel{\eqref{eq:79}}{\geq}2-\alpha\geq2(1-\tau).
\end{align*}
Hence we can conclude from \eqref{eq:76} that
\begin{equation}\label{eq:80}
\sum_{i\in I}\norm{x_i^*}\norm{x_i+y_i}\geq2(1-\nu).
\end{equation}
Using \eqref{eq:79} and \eqref{eq:80} we get
\begin{align*}
& \norm{x+z}_E\geq\sum_{i\in I}\norm{x_i^*}\norm{x_i+z_i} \\
& \geq \sum_{i\in I}\norm{x_i^*}\norm{x_i+y_i}-\sum_{i\in I}\norm{x_i^*}\norm{y_i-z_i} \\
& \geq2(1-\nu)-\alpha>2(1-\eta)
\end{align*}
and thus the choice of $\eta$ implies $x^*(z)\geq1-\eps/2$. But from \eqref{eq:79} it also follows 
that $\abs{x^*(y)-x^*(z)}\leq\alpha$ and hence $x^*(y)\geq1-\eps/2-\alpha>1-\eps$.
\end{Proof}

Because of the uniform rotundity of $\ell^p(I)$ for $1<p<\infty$ we have the following corollary.
\begin{corollary}\label{cor:p-sum uacs}
If $(X_i)_{i\in I}$ is a family of Banach spaces such that for every $0<\eps\leq2$ 
we have $\inf_{i\in I}\delta_{\mathrm{uacs}}^{X_i}(\eps)>0$ then 
$\bigl[\bigoplus_{i\in I}X_i\bigr]_p$ is also \ac{uacs} for every $1<p<\infty$.
\end{corollary}
We can also get a more general corollary for a \ac{US} space $E$.
\begin{corollary}\label{cor:smooth sum uacs}
If $(X_i)_{i\in I}$ is a family of Banach spaces such that for every $0<\eps\leq2$
we have $\delta(\eps):=\inf_{i\in I}\delta_{\mathrm{uacs}}^{X_i}(\eps)>0$ and if $E$
is \ac{US} then $\bigl[\bigoplus_{i\in I}X_i\bigr]_E$ is also a \ac{uacs} space.
\end{corollary}

\begin{Proof}
Since $E$ is \ac{US} it is reflexive and hence it cannot contain an isomorpic copy of $\ell^1$.
Thus by Lemma \ref{lemma:abs norms l1} $\lin\set*{e_i:i\in I}$ is dense in $E^{\prime}$.\par
Further, since $E$ is \ac{US} the dual space $E^*=E^{\prime}$ is \ac{UR}, as already mentioned 
in the introduction. Because the spaces $X_i$ are \ac{uacs} they are also reflexive and hence 
Proposition \ref{prop:uacs-dual} and the monotonicity of the functions $\delta_{\mathrm{uacs}}^{X_i}$ gives us 
$\inf_{i\in I}\delta_{\mathrm{uacs}}^{X_i^*}(\eps)\geq\delta(\delta(\eps))>0$ for every $0<\eps\leq2$.\par
So by Theorem \ref{thm:sum uacs} the space $\bigl[\bigoplus_{i\in I}X_i^*\bigr]_{E^{\prime}}=
\bigl[\bigoplus_{i\in I}X_i\bigr]_E^*$ is \ac{uacs} and hence $\bigl[\bigoplus_{i\in I}X_i\bigr]_E$ is also
\ac{uacs} by Proposition \ref{prop:uacs-dual}.
\end{Proof}

Finally, we summarise all the results on absolute sums we have obtained in this 
section in the following table.

\begin{table}[H]
\begin{center}
\caption{Summary of the results}\label{tab:summary}
\ \\
\begin{tabular}{ccc}
        $E$                             &               $X_i$                           & $\bigl[\bigoplus_{i\in I}X_i\bigr]_E$ \\ \addlinespace[1pt] \toprule
\ac{acs}                                &             \ac{acs}                          &            \ac{acs}                   \\ \midrule
\ac{luacs} + $(P)$                      &             \ac{luacs}                        &            \ac{luacs}                 \\ \midrule
$\mathrm{luacs}^+$                      &             \ac{luacs}                        &            \ac{luacs}                 \\ \midrule
$\mathrm{luacs}^+$ + $(P)$              &         $\mathrm{luacs}^+$                    &        $\mathrm{luacs}^+$             \\ \midrule
$\mathrm{luacs}^+$ + $\ell^1\not\ssq E$ &         $\mathrm{luacs}^+$                    &        $\mathrm{luacs}^+$             \\ \midrule
\ac{sluacs} + $(P)$                     &             \ac{sluacs}                       &            \ac{sluacs}                \\ \midrule
$\mathrm{sluacs}^+$                     &             \ac{sluacs}                       &            \ac{sluacs}                \\ \midrule
$\mathrm{sluacs}^+$                     &         $\mathrm{luacs}^+$                    &        $\mathrm{luacs}^+$             \\ \midrule 
$\mathrm{sluacs}^+$                     &         $\mathrm{sluacs}^+$                   &        $\mathrm{sluacs}^+$            \\ \midrule
\ac{wuacs} + $\ell^1\not\ssq E$         &             \ac{wuacs}                        &            \ac{wuacs}                 \\ \midrule
\ac{acs} + $I$ finite                   &             \ac{uacs}                         &            \ac{uacs}                  \\ \midrule
$(u^+)$                                 & $\inf_{i\in I}\delta_{\mathrm{uacs}}^{X_i}>0$ &            \ac{uacs}                  \\ \midrule
\ac{US}                                 & $\inf_{i\in I}\delta_{\mathrm{uacs}}^{X_i}>0$ &            \ac{uacs}                  \\ \bottomrule
\end{tabular}
\end{center}
\end{table}

\section{Midpoint versions}\label{sec:midpoint versions}
Recall that a Banach space $X$ is said to be {\em \ac{MLUR}} if for any two sequence $(x_n)_{n\in \N}$
and $(y_n)_{n\in \N}$ in $S_X$ and every $x\in S_X$ we have
\begin{equation*}
\norm*{x-\frac{x_n+y_n}{2}}\to 0 \ \Rightarrow \ \norm{x_n-y_n}\to 0.
\end{equation*}
This notion was originally introduced in \cite{anderson}.\par
Recall also that $X$ is called {\em \ac{WMLUR}} if it satisfies the above condition 
with $\norm{x_n-y_n}\to 0$ replaced by $x_n-y_n\xrightarrow{\sigma} 0$, where the  
symbol $\xrightarrow{\sigma}$ denotes the convergence in the weak topology of $X$.\par
We now introduce in an analogous way midpoint versions of \ac{luacs} and \ac{sluacs} spaces.

\begin{definition}
Let $X$ be a Banach space.
  \begin{enumerate}[(i)]
  \item The space $X$ is said to be {\em \ac{mluacs}} if for any two sequences $(x_n)_{n\in \N}$ 
  and $(y_n)_{n\in \N}$ in $S_X$, every $x\in S_X$ and every $x^*\in S_{X^*}$ we have that
  \begin{equation*}
  \norm*{x-\frac{x_n+y_n}{2}}\to 0 \ \mathrm{and} \ x^*(x_n)\to 1 \ \Rightarrow \ x^*(y_n)\to 1.
  \end{equation*}
  \item The space $X$ is called {\em \ac{msluacs}} if for any two sequences $(x_n)_{n\in \N}$ and 
  $(y_n)_{n\in \N}$ in $S_X$, every $x\in S_X$ and every sequence $(x_n^*)_{n\in \N}$ in $S_{X^*}$ we 
  have that
  \begin{equation*}
  \norm*{x-\frac{x_n+y_n}{2}}\to 0 \ \mathrm{and} \ x_n^*(x_n)\to 1 \ \Rightarrow \ x_n^*(y_n)\to 1.
  \end{equation*}
  \end{enumerate}
\end{definition}

We then get the following implication chart.
\begin{figure}[H]
\begin{center}
  \begin{tikzpicture}
  \node (LUR) at (-2,0) {LUR};
  \node (MLUR) at (0,1) {MLUR};
  \node (WLUR) at (0,-1) {WLUR};
  \node (WMLUR) at (2,0) {WMLUR};
  \node (R) at (4,0) {R};
  \node (sluacs) at (-2,-1) {sluacs};
  \node (msluacs) at (0,0) {msluacs};
  \node (luacs) at (0,-2) {luacs};
  \node (mluacs) at (2,-1) {mluacs};
  \node (acs) at (4,-1) {acs};
  \draw[->] (LUR)--(MLUR);
  \draw[->] (LUR)--(WLUR);
  \draw[->] (WLUR)--(WMLUR);
  \draw[->] (MLUR)--(WMLUR);
  \draw[->] (WMLUR)--(R);
  \draw[->] (sluacs)--(msluacs);
  \draw[->] (sluacs)--(luacs);
  \draw[->] (luacs)--(mluacs);
  \draw[->] (msluacs)--(mluacs);
  \draw[->] (mluacs)--(acs);
  \draw[->] (LUR)--(sluacs);
  \draw[->] (MLUR)--(msluacs);
  \draw[->] (WLUR)--(luacs);
  \draw[->] (WMLUR)--(mluacs);
  \draw[->] (R)--(acs);
  \end{tikzpicture}
\end{center}
\CAP\label{fig:7}
\end{figure}

No other implications are valid in general, as is shown by the examples in Section \ref{sec:examples}.\par
Note that in the definition of \ac{msluacs} spaces we can replace the condition $x_n^*(x_n)\to 1$ by
$x_n^*(x_n)=1$ for every $n\in \N$ and obtain an equivalent definition, by the same argument used in
the proof of Proposition \ref{prop:char-uacs-sluacs}. Also, it is well known (and not to hard to see)
that a Banach space $X$ is \ac{MLUR} (resp. \ac{WMLUR}) \ifif for every sequence $(x_n)_{n\in \N}$ in 
$X$ and each element $x\in X$ the condition $\norm{x\pm x_n}\to \norm{x}$ implies $\norm{x_n}\to 0$ 
(resp. $x_n\xrightarrow{\sigma} 0$). In much the same way one can prove that $X$ is \ac{msluacs}
\ifif for every sequence $(x_n)_{n\in \N}$ in $X$, each $x\in X$ and every bounded sequence $(x_n^*)_{n\in \N}$
in $X^*$ the two conditions $\norm{x\pm x_n}\to \norm{x}$ and $x_n^*(x+x_n)-\norm{x_n^*}\norm{x}\to 0$ imply
$x_n^*(x_n)\to 0$ and that an analogous characterisation holds for \ac{mluacs} spaces.\par
It was noted in \cite{smith2}*{p.663} that by using the principle of local reflexivity one can easily check that
$X$ is \ac{WMLUR} \ifif every point $x\in S_X$ is an extreme point of $B_{X^{**}}$\footnote{We consider $X$ canonically 
embedded into its second dual.}(in particular, \ac{WMLUR} and \ac{R} coincide in reflexive spaces). In analogy to this 
result we can prove the following characterisation of \ac{mluacs} spaces, which especially yields that \ac{mluacs} and 
\ac{acs} coincide in reflexive spaces.\par
\begin{proposition}\label{prop:second dual char mluacs}
A Banach space $X$ is \ac{mluacs} \ifif the following holds: for any two elements $x^{**},y^{**}\in S_{X^{**}}$
with $x^{**}+y^{**}\in 2S_X$ and every $x^{*}\in S_{X^*}$ with $x^{**}(x^{*})=1$ we also have $y^{**}(x^{*})=1$.
\end{proposition}

\begin{Proof}
The sufficiency is straightforwardly proved using the weak*-compact\-ness of the bidual unit ball.\par
To prove the necessity, fix $x^{**},y^{**}\in S_{X^{**}}$ such that $x^{**}+y^{**}\in 2S_X$ and $x^{*}\in S_{X^*}$ 
with $x^{**}(x^*)=1$. Put $F=\lin\set*{x^{**},y^{**}}$. By the principle of local reflexivity (cf. \cite{albiac}*{Theorem 11.2.4}) 
we can find for each $n\in \N$ a finite-dimensional subspace $E_n\ssq X$ and an isomorphism $T_n:F \rightarrow E_n$ such that 
$\norm{T_n}\leq 1+2^{-n}$, $\norm{T_n^{-1}}\leq1+2^{-n}$, $T_nz=z$ for every $z\in X\cap F$ and $x^{*}(T_nz^{**})=z^{**}(x^{*})$ 
for all $z^{**}\in F$.\par
If we put $x_n=T_nx^{**}$ and $y_n=T_ny^{**}$ for every $n\in \N$ then we have $x_n+y_n=2(x^{**}+y^{**})$ and 
$x^{*}(x_n)=1$ as well as $\norm{x_n}, \norm{y_n}\to 1$. Since $X$ is \ac{mluacs} it follows that $x^{*}(y_n)\to 1$.
But $x^{*}(y_n)=y^{**}(x^{*})$ for every $n$, thus $y^{**}(x^{*})=1$.
\end{Proof}

We can also prove a characterisation of \ac{msluacs} spaces that is analogous to the results on \ac{acs}, \ac{sluacs} and \ac{uacs}
spaces given in the Propositions \ref{prop:char acs without dual}, \ref{prop:char sluacs without dual} and \ref{prop:char uacs without dual}
(the proof is completely analogous as well).
\begin{proposition}\label{prop:char msluacs without dual}
For a Banach space $X$ the following assertions are equivalent:
\begin{enumerate}[\upshape(i)]
\item $X$ is \ac{msluacs}.
\item For every $z\in S_X$ and every $\eps>0$ there exists some $\delta>0$ such that for every $t\in [0,\delta]$ and all $x,y\in S_X$
with $\norm{x+y-2z}\leq 2t$ we have
\begin{equation*}
\norm{x+ty}+\norm{x-ty}\leq 2+\eps t.
\end{equation*}
\item For every $z\in S_X$ and every $\eps>0$ there exists some $\delta>0$ such that for every $t\in [0,\delta]$ and all $x,y\in S_X$
with $\norm{x+y-2z}\leq \delta t$ we have
\begin{equation*}
\norm{x-ty}\leq 1+t(\eps-1).
\end{equation*}
\item For every $z\in S_X$ there exists some $1\leq p<\infty$ such that for every $\eps>0$ there is some $\delta>0$ such that for all
$t\in [0,\delta]$ and all $x,y\in S_X$ with $\norm{x+y-2z}\leq 2t$ we have
\begin{equation*}
\norm{x+ty}^p+\norm{x-ty}^p\leq 2+\eps t^p.
\end{equation*}
\item For every $z\in S_X$ there exists some $1\leq p<\infty$ such that for every $\eps>0$ there is some $\delta>0$ such that for all
$t\in [0,\delta]$ and all $x,y\in S_X$ with $\norm{x+y-2z}\leq t\delta$ we have
\begin{equation*}
(1+t)^p+\norm{x-ty}^p\leq 2+\eps t^p.
\end{equation*}
\end{enumerate}
\end{proposition}
Recall that the space $X$ is said to have the property $(H)$ (also known as the Kadets-Klee property) if for every 
sequence $(x_n)_{n\in \N}$ in $X$ and each $x\in X$ the implication 
\begin{equation*}
x_n\xrightarrow{\sigma} x \ \mathrm{and} \ \norm{x_n}\to \norm{x} \ \Rightarrow \ \norm{x_n-x}\to 0
\end{equation*}
holds. For example, every \ac{LUR} space has property $(H)$.\par
It was proved in \cite{kadets82} that a Banach space which is \ac{R}, has property $(H)$ and does not contain 
an isomorpic copy of $\ell^1$ is actually \ac{MLUR}. We can adopt the proof from \cite{kadets82} to show the 
analogous result for \ac{acs} spaces.\par
\begin{proposition}\label{prop:acs H msluacs}
Let $X$ be an \ac{acs} space which has property $(H)$ and does not contain an isomorphic copy of $\ell^1$.
Then $X$ is \ac{msluacs}.
\end{proposition}

\begin{Proof}
Take a sequence $(x_n)_{n\in \N}$ in $X$ and $x\in X$ with $\norm{x_n\pm x}\to \norm{x}$ and a sequence $(x_n^*)_{n\in \N}$ 
in $S_{X^*}$ such that $x_n^*(x_n+x)\to \norm{x}$. If $(x_n^*(x_n))_{n\in \N}$ was not convergent to zero then by passing 
to an appropriate subsequence we could assume $\abs*{x_n^*(x_n)}\geq\eps$ for all $n$ and some $\eps>0$.\par
Since $X$ does not contain $\ell^1$ we can, by Rosenthal's theorem (cf. \cite{albiac}*{Theorem 10.2.1}), pass to a 
further subsequence such that $(x_n)_{n\in \N}$ is weakly Cauchy. But then the double-sequence $(x_n-x_m)_{n,m\in \N}$ 
is weakly null, so $\liminf \norm{2x+x_n-x_m}\geq 2\norm{x}$.\par
On the other hand, because of 
\begin{equation*}
\norm{2x+x_n-x_m}\leq \norm{x+x_n}+\norm{x-x_m} \ \ \forall n,m\in \N
\end{equation*}
we also have $\limsup \norm{2x+x_n-x_m}\leq 2\norm{x}$ and hence $\lim \norm{2x+x_n-x_m}=2\norm{x}$.
Since $2x+x_n-x_m\xrightarrow{\sigma} 2x$ the property $(H)$ of $X$ implies that $\lim\norm{x_n-x_m}=0$,
i.\,e. $(x_n)_{n\in \N}$ is norm Cauchy. Let $y$ be the limit of $(x_n)_{n\in \N}$. It follows that
$\norm{x\pm y}=\norm{x}$, so if we put $z_1=(x+y)/\norm{x}$ and $z_2=(x-y)/\norm{x}$ then 
$\norm{z_1}=\norm{z_2}=1$ and $\norm{z_1+z_2}=2$.\par
Since $B_{X^*}$ is weak*-compact we can pass to a subnet $(x_{\varphi(i)}^*)_{i\in I}$ that is weak*-convergent
to some $x^*\in B_{X^*}$. Because of $x_{\varphi(i)}^*(x_{\varphi(i)}+x)\to \norm{x}$ this implies $x_{\varphi(i)}^*(x_{\varphi(i)})\to \norm{x}-x^*(x)$.
But we also have $x_{\varphi(i)}^*(y)\to x^*(y)$ and $\norm{x_{\varphi(i)}-y}\to 0$, thus $x_{\varphi(i)}^*(x_{\varphi(i)})\to x^*(y)$ and 
hence $x^*(z_1)=1$.\par
Because of $\abs{x_{\varphi(i)}^*(x_{\varphi(i)})}\geq\eps$ for every $i\in I$ we have $x^*(y)\neq 0$. It follows that $x^*(z_2)\neq 1$ and 
hence $X$ cannot be \ac{acs}, contradicting our hypothesis.
\end{Proof}

In the spirit of the duality results from section \ref{sec:general facts} it is also not difficult to prove
the following assertions (we omit the details).
\begin{proposition}\label{prop:dual mluacs}
Let $X$ be a Banach space such that for all functionals $x^*,y^*\in S_{X^*}$ and every $x\in S_X$ the
implication
\begin{equation*}
\norm{x^*+y^*}=2 \ \mathrm{and} \ x^*(x)=1 \ \Rightarrow \ y^*(x)=1
\end{equation*}
is valid. Then $X$ is \ac{mluacs}. In particular, $X$ is \ac{mluacs} whenever $X^*$ is \ac{acs}.
\end{proposition}

\begin{proposition}\label{prop:dual msluacs}
Let $X$ be a Banach space such that for every sequence $(x_n^*)_{n\in \N}$ in $S_{X^*}$, every $x^*\in S_{X^*}$
and each $x\in S_X$ the implication
\begin{equation*}
\norm{x_n^*+x^*}\to 2 \ \mathrm{and} \ x^*(x)=1 \ \Rightarrow \ x_n^*(x)\to 1
\end{equation*}
is valid. Then $X$ is \ac{msluacs}. In particular, $X$ is \ac{msluacs} whenever $X^*$ is $\text{luacs}^+$.
\end{proposition}
There are some results on absolute sums of \ac{msluacs} and \ac{mluacs} spaces which we will prove in the
following. The proof of the first one uses ideas from the proof of \cite{dowling2}*{Proposition 4}.
The notation is the same as in the previous section.
\begin{proposition}\label{prop:sum msluacs MLUR}
If $(X_i)_{i\in I}$ is a family of \ac{msluacs} (resp. \ac{mluacs}) Banach spaces and if $E$ is \ac{MLUR}
then $\bigl[\bigoplus_{i\in I}X_i\bigr]_E$ is also \ac{msluacs} (resp. \ac{mluacs}).
\end{proposition}

\begin{Proof}
Suppose all the $X_i$ are \ac{msluacs} and take two sequences $(x_n)_{n\in \N}$ and $(y_n)_{n\in \N}$ as well 
as an element $x=(x_i)_{i\in I}$ in the unit sphere of $\bigl[\bigoplus_{i\in I}X_i\bigr]_E$ such that $\norm{x_n+y_n-2x}_E\to 0$.
Also, fix a sequence $(x_n^*)_{n\in \N}$ of norm one functionals with $x_n^*(x_n)\to 1$. We write $x_n=(x_{n,i})_{i\in I}$,
$y_n=(y_{n,i})_{i\in I}$ and $x_n^*=(x_{n,i}^*)_{i\in I}$ for each $n\in \N$.\par
As we have done many times before, we conclude 
\begin{equation}\label{eq:81}
\lim_{n\to \infty}\paren*{x_{n,i}^*(x_{n,i})-\norm{x_{n,i}^*}\norm{x_{n,i}}}=0 \ \ \forall i\in I
\end{equation}
and
\begin{equation}\label{eq:82}
\lim_{n\to \infty}\sum_{i\in I}\norm{x_{n,i}^*}\norm{x_{n,i}}=1
\end{equation}
as well as
\begin{equation}\label{eq:83}
\lim_{n\to \infty}\norm*{(\norm{x_{n,i}}+\norm{y_{n,i}})_{i\in I}}_E=2.
\end{equation}
We also have 
\begin{equation}\label{eq:84}
\lim_{n\to \infty}\norm{x_{n,i}+y_{n,i}-2x_i}=0 \ \ \forall i\in I.
\end{equation}
Because of Lemma \ref{lemma:abs norms} we get
\begin{equation*}
\norm*{(2\norm{x_i}-\norm{x_{n,i}+y_{n,i}})_{i\in I}}_E\leq\norm*{2x-x_n-y_n}_E
\end{equation*}
and hence
\begin{equation}\label{eq:85}
\lim_{n\to \infty}\norm*{(2\norm{x_i}-\norm{x_{n,i}+y_{n,i}})_{i\in I}}_E=0.
\end{equation}
Now we put for every $n\in \N$
\begin{equation}\label{eq:86}
a_n=(a_{n,i})_{i\in I}=\paren*{2\norm{x_i}-\frac{1}{2}(\norm{x_{n,i}}+\norm{y_{n,i}})}_{i\in I}
\end{equation}
and
\begin{equation}\label{eq:87}
b_n=(b_{n,i})_{i\in I}=\paren*{\norm{x_i}-\frac{1}{2}\norm{x_{n,i}+y_{n,i}}}_{i\in I}.
\end{equation}
We then have
\begin{align}
\nonumber
&\norm{x_i}\leq\norm{x_i}-\frac{1}{2}\norm{x_{n,i}+y_{n,i}}+\frac{1}{2}(\norm{x_{n,i}}+\norm{y_{n,i}}) \\
\label{eq:88}
&=b_{n,i}+\frac{1}{2}(\norm{x_{n,i}}+\norm{y_{n,i}}) \ \ \forall i\in I.
\end{align}
If $a_{n,i}\geq0$ then 
\begin{equation*}
\abs{a_{n,i}}=a_{n,i}=2\norm{x_i}-\frac{1}{2}(\norm{x_{n,i}}+\norm{y_{n,i}})
\stackrel{\eqref{eq:88}}{\leq}2\abs{b_{n,i}}+\frac{1}{2}(\norm{x_{n,i}}+\norm{y_{n,i}})
\end{equation*}
and if $a_{n,i}<0$ then
\begin{equation*}
\abs{a_{n,i}}=-a_{n,i}=\frac{1}{2}(\norm{x_{n,i}}+\norm{y_{n,i}})-2\norm{x_i}\leq
2\abs{b_{n,i}}+\frac{1}{2}(\norm{x_{n,i}}+\norm{y_{n,i}}).
\end{equation*}
Thus we have
\begin{equation}\label{eq:89}
\abs{a_{n,i}}\leq2\abs{b_{n,i}}+\frac{1}{2}(\norm{x_{n,i}}+\norm{y_{n,i}}) \ \ \forall i\in I, \forall n\in \N.
\end{equation}
Using \eqref{eq:89} and Lemma \ref{lemma:abs norms} we deduce that
\begin{align*}
&\frac{1}{2}\norm*{(\norm{x_{n,i}}+\norm{y_{n,i}})_{i\in I}}_E+2\norm{b_n}_E\geq
\norm*{\paren*{2\abs{b_{n,i}}+\frac{1}{2}(\norm{x_{n,i}}+\norm{y_{n,i}})}_{i\in I}}_E \\
&\geq\norm{a_n}_E\geq2-\frac{1}{2}\norm{(\norm{x_{n,i}}+\norm{y_{n,i}})_{i\in I}}_E
\end{align*}
holds for every $n\in \N$ which together with \eqref{eq:83} and \eqref{eq:85} implies that
$\norm{a_n}_E\to 1$. Because of the definition of the sequence $(a_n)_{n\in \N}$ and the fact 
that $E$ is \ac{MLUR} this leads to $\norm*{2a_n-(\norm{x_{n,i}}+\norm{y_{n,i}})_{i\in I}}_E\to 0$,
in other words
\begin{equation}\label{eq:90}
\lim_{n\to \infty}\norm*{\paren*{\norm{x_i}-\frac{1}{2}(\norm{x_{n,i}}+\norm{y_{n,i}})}_{i\in I}}_E=0.
\end{equation}
Again, since $E$ is \ac{MLUR} it follows from \eqref{eq:90} that
\begin{equation}\label{eq:91}
\lim_{n\to \infty}\norm*{(\norm{x_{n,i}}-\norm{y_{n,i}})_{i\in I}}_E=0.
\end{equation}
From \eqref{eq:90} resp. \eqref{eq:91} it follows that $\norm{x_{n,i}}+\norm{y_{n,i}}\to 2\norm{x_i}$
and $\norm{x_{n,i}}-\norm{y_{n,i}}\to 0$ for every $i\in I$ and hence
\begin{equation}\label{eq:92}
\lim_{n\to \infty}\norm{x_{n,i}}=\norm{x_i}=\lim_{n\to \infty}\norm{y_{n,i}} \ \ \forall i\in I.
\end{equation}
Since each $X_i$ is \ac{msluacs} it follows from \eqref{eq:81}, \eqref{eq:84} and \eqref{eq:92} that
\begin{equation}\label{eq:93}
\lim_{n\to \infty}(x_{n,i}^*(x_i)-\norm{x_{n,i}^*}\norm{x_i})=0 \ \ \forall i\in I.
\end{equation}
From \eqref{eq:90} and \eqref{eq:91} we get $\norm*{(\norm{x_i}-\norm{x_{n,i}})_{i\in I}}_E\to 0$ which
together with \eqref{eq:82} implies
\begin{equation}\label{eq:94}
\lim_{n\to \infty}\sum_{i\in I}\norm{x_{n,i}^*}\norm{x_i}=1.
\end{equation}
From \eqref{eq:93} we can infer exactly as in the proof of Proposition \ref{prop:sum sluacs (P)} that
$x_n^*(x)-\sum_{i\in I}\norm{x_{n,i}^*}\norm{x_i}\to 0$. By \eqref{eq:94} this implies $x_n^*(x)\to 1$. 
Together with $x_n^*(x_n)\to 1$ and $\norm{x_n+y_n-2x}_E\to 0$ it follows $x_n^*(y_n)\to 1$ and 
the proof is finished. The case of \ac{mluacs} spaces is proved analogously.
\end{Proof}

As mentioned before, $\ell^p(I)$ is \ac{UR} for every $1<p<\infty$ thus we get the following corollary.
\begin{corollary}\label{cor:p-sum msluacs}
If $(X_i)_{i\in I}$ is a family of \ac{msluacs} (resp. \ac{mluacs}) Banach spaces then 
$\bigl[\bigoplus_{i\in I}X_i\bigr]_p$ is also \ac{msluacs} (resp. \ac{mluacs}) for every $1<p<\infty$.
\end{corollary}
The second result on sums of \ac{msluacs} (and \ac{mluacs}) spaces reads as follows.
\begin{proposition}\label{prop:sum msluacs (P)}
If $(X_i)_{i\in I}$ is a family of \ac{msluacs} (resp. \ac{mluacs}) spaces and if $E$ is 
\ac{sluacs} (resp. \ac{luacs}) and has the property $(P)$ then $\bigl[\bigoplus_{i\in I}X_i\bigr]_E$ 
is also \ac{msluacs} (resp. \ac{mluacs}).
\end{proposition}

\begin{Proof}
We suppose that every $X_i$ is \ac{msluacs} and that $E$ is \ac{sluacs}. Let us fix sequences $(x_n)_{n\in \N}$, 
$(y_n)_{n\in \N}$ and $(x_n^*)_{n\in \N}$ as well as an element $x$ just as in the proof of Proposition 
\ref{prop:sum msluacs MLUR}. Exactly as in this proof we can show that
\begin{equation}\label{eq:95}
\lim_{n\to \infty}\paren*{x_{n,i}^*(x_{n,i})-\norm{x_{n,i}^*}\norm{x_{n,i}}}=0 \ \ \forall i\in I,
\end{equation}
\begin{equation}\label{eq:96}
\lim_{n\to \infty}\sum_{i\in I}\norm{x_{n,i}^*}\norm{x_{n,i}}=1
\end{equation}
and
\begin{equation}\label{eq:97}
\lim_{n\to \infty}\norm{x_{n,i}+y_{n,i}-2x_i}=0 \ \ \forall i\in I.
\end{equation}
From $\norm{x_n+y_n-2x}_E\to 0$ and $\norm{x_n}_E=\norm{x}_E=\norm{y_n}_E=1$ for every $n$ we can infer that 
$\norm{x+x_n}_E\to 2$ and $\norm{x+y_n}_E\to 2$ which in turn implies
\begin{equation}\label{eq:98}
\lim_{n\to \infty}\norm*{(\norm{x_{n,i}}+\norm{x_i})_{i\in I}}_E=2
=\lim_{n\to \infty}\norm*{(\norm{y_{n,i}}+\norm{x_i})_{i\in I}}_E.
\end{equation}
Since $E$ is \ac{sluacs} it follows from \eqref{eq:96} and \eqref{eq:98} that
\begin{equation}\label{eq:99}
\lim_{n\to \infty}\sum_{i\in I}\norm{x_{n,i}^*}\norm{x_i}=1
\end{equation}
and since $E$ has the property $(P)$ we also get from \eqref{eq:98} that
\begin{equation}\label{eq:100}
\lim_{n\to \infty}\norm{x_{n,i}}=\norm{x_i}=\lim_{n\to \infty}\norm{y_{n,i}} \ \ \forall i\in I.
\end{equation}
Now we can use \eqref{eq:95}, \eqref{eq:97}, \eqref{eq:100} and the fact that each $X_i$ is \ac{msluacs}
to get
\begin{equation}\label{eq:101}
\lim_{n\to \infty}\paren*{x_{n,i}^*(x_i)-\norm{x_{n,i}^*}\norm{x_i}}=0 \ \ \forall i\in I
\end{equation}
and then based on \eqref{eq:99} and \eqref{eq:101} the rest of the proof can be carried out exactly as 
the proof of Proposition \ref{prop:sum msluacs MLUR}. Again, the \ac{mluacs} case is proved analogously.
\end{Proof}

\section{Directional versions}\label{sec:direc versions}
Recall that a Banach space $X$ is said to be {\em \ac{URED}} if for any $z\in X\sm \set*{0}$ and all sequences $(x_n)_{n\in \N}$
and $(y_n)_{n\in \N}$ in $S_X$ such that $\norm{x_n+y_n}\to 2$ and $x_n-y_n\in \lin\set*{z}$ for every $n$ one already has
$\norm{x_n-y_n}\to 0$.\par
This notion was first introduced by Garkavi in \cite{garkavi} and further studied by the authors of \cite{dayjamesswaminathan}.\par
In this spirit, we define the following directionalisation of \ac{uacs} spaces.
\begin{definition}\label{def:uacsed}
A Banach space $X$ is called {\em \ac{uacsed}} if for every $z\in X\sm \set*{0}$ and all sequences $(x_n)_{n\in \N},
(y_n)_{n\in \N}$ in $S_X$ and $(x_n^*)_{n\in \N}$ in $S_{X^*}$ such that $\norm{x_n+y_n}\to 2$, $x_n^*(x_n)\to 1$ and
$x_n-y_n\in \lin\set*{z}$ for every $n$ one also has $x_n^*(y_n)\to 1$.
\end{definition}
Obviously, the following implications hold.
\begin{figure}[H]
\begin{center}
  \begin{tikzpicture}
  \node (UR) at (-2,1) {UR};
  \node (URED) at (0,1) {URED};
  \node (R) at (2,1) {R};
  \node (uacs) at (-2,-0.5) {uacs};
  \node (uacsed) at (0,-0.5) {uacsed};
  \node (acs) at (2,-0.5) {acs};
  \draw[->] (UR)--(URED);
  \draw[->] (URED)--(R);
  \draw[->] (uacs)--(uacsed);
  \draw[->] (uacsed)--(acs);
  \draw[->] (UR)--(uacs);
  \draw[->] (URED)--(uacsed);
  \draw[->] (R)--(acs);
  \end{tikzpicture}
\end{center}
\CAP\label{fig:8}
\end{figure}
Let us now give some equivalent characterisations for a Banach space to be \ac{uacsed} in analogy to the characterisations
for a Banach space to be \ac{URED} given in \cite{dayjamesswaminathan}*{Theorem 1}.
\begin{proposition}\label{prop:char uacsed}
Let $X$ be a Banach space and $2\leq p<\infty$. Then the following assertions are equivalent.
  \begin{enumerate}[\upshape(i)]
  \item $X$ is \ac{uacsed}.
  \item For all sequences $(x_n)_{n\in \N}, (y_n)_{n\in \N}$ in $B_X$, $(x_n^*)_{n\in \N}$ in $S_{X^*}$ and every $z\in X$
  the implication 
  \begin{equation*}
  \norm{x_n+y_n}\to 2, \ x_n^*(x_n)\to 1 \ \mathrm{and} \ x_n-y_n\to z \ \Rightarrow \ x_n^*(z)\to 0
  \end{equation*}
  holds.
  \item For all sequences $(x_n)_{n\in \N}, (y_n)_{n\in \N}$ in $S_X$, $(x_n^*)_{n\in \N}$ in $S_{X^*}$ and every $z\in X$
  the implication 
  \begin{equation*}
  \norm{x_n+y_n}\to 2, \ x_n^*(x_n)\to 1 \ \mathrm{and} \ x_n-y_n\to z \ \Rightarrow \ x_n^*(z)\to 0
  \end{equation*}
  holds.
  \item For all sequences $(x_n)_{n\in \N}, (y_n)_{n\in \N}$ in $B_X$, $(x_n^*)_{n\in \N}$ in $S_{X^*}$ and every $z\in X\sm \set*{0}$
  the conditions
  \begin{equation*}
  \norm{x_n+y_n}\to 2, \ x_n^*(x_n)\to 1 \ \mathrm{and} \ x_n-y_n\in \lin\set*{z} \ \forall n\in \N
  \end{equation*}
  imply $x_n^*(y_n)\to 1$.
  \item For all sequences $(x_n)_{n\in \N}, (y_n)_{n\in \N}$ in $S_X$, $(x_n^*)_{n\in \N}$ in $S_{X^*}$ and every $z\in X$
  the implication 
  \begin{equation*}
  \norm{x_n+y_n}\to 2, \ x_n^*(x_n)=1 \ \forall n\in \N \ \mathrm{and} \ x_n-y_n\to z \ \Rightarrow \ x_n^*(z)\to 0
  \end{equation*}
  holds.
  \item For all sequences $(x_n)_{n\in \N}$ in $B_X$, $(x_n^*)_{n\in \N}$ in $S_{X^*}$ and every $z\in X$ the two conditions
  \begin{equation*}
  2^{p-1}\paren*{\norm{x_n+z}^p+\norm{x_n}^p}-\norm{2x_n+z}^p\to 0 \ \mathrm{and} \ x_n^*(x_n)\to 1
  \end{equation*}
  imply that $x_n^*(z)\to 0$.
  \item For all sequences $(x_n)_{n\in \N}$ in $S_X$, $(x_n^*)_{n\in \N}$ in $S_{X^*}$ and every $z\in X$ the two conditions
  \begin{equation*}
  \norm{2x_n+z}^p-2^{p-1}\norm{x_n+z}^p\to 2^{p-1} \ \mathrm{and} \ x_n^*(x_n)=1 \ \forall n\in \N
  \end{equation*}
  imply that $x_n^*(z)\to 0$.
  \item For all sequences $(x_n)_{n\in \N}$ in $B_X$, $(y_n)_{n\in \N}$ in $X$, $(x_n^*)_{n\in \N}$ in $S_{X^*}$ and every $z\in X$ 
  the conditions
  \begin{equation*}
  \norm{x_n+y_n}\to 2, \norm{y_n}\to 1, x_n^*(x_n)\to 1 \ \mathrm{and} \ x_n-y_n=z \ \forall n\in \N
  \end{equation*}
  imply that $x_n^*(z)\to 0$.
  \item For all sequences $(x_n)_{n\in \N}$ in $S_X$, $(y_n)_{n\in \N}$ in $X$, $(x_n^*)_{n\in \N}$ in $S_{X^*}$ and every $z\in X$ 
  the conditions
  \begin{equation*}
  \norm{x_n+y_n}\to 2, \norm{y_n}\to 1, x_n^*(x_n)=1 \ \forall n\in \N \ \mathrm{and} \ x_n-y_n=z \ \forall n\in \N
  \end{equation*}
  imply that $x_n^*(z)\to 0$.
  \end{enumerate}
\end{proposition}

\begin{Proof}
We first prove $\mathrm{(i)} \Rightarrow \mathrm{(ii)}$. So let us fix two sequences $(x_n)_{n\in \N}$, $(y_n)_{n\in \N}$ in $B_X$,
a sequence $(x_n^*)_{n\in \N}$ in $S_{X^*}$ and a $z\in X\sm\set*{0}$ such that $\norm{x_n+y_n}\to 2$, $x_n^*(x_n)\to 1$ and $x_n-y_n\to z$.\par
If there is a subsequence $(x_{n_k}^*(z))_{k\in \N}$ such that $x_{n_k}^*(z)\leq0$ for every $k\in \N$ then because of $x_n-y_n\to z$ and 
$x_n^*(x_n)\to 1$ we have
\begin{equation*}
0\geq\limsup_{k\to \infty}x_{n_k}^*(z)=\limsup_{k\to \infty}x_{n_k}^*(x_{n_k}-y_{n_k})=1-\liminf_{k\to \infty}x_{n_k}^*(y_{n_k})\geq0,
\end{equation*}
so $\limsup_{k\to \infty}x_{n_k}^*(z)=0$ and analogously $\liminf_{k\to \infty}x_{n_k}^*(z)=0$, hence $\lim_{k\to \infty}x_{n_k}^*(z)=0$.\par
Otherwise the sequence $(x_n^*(z))_{n\in \N}$ is eventually positive. Since $\norm{x_n+y_n}\to 2$ and $\norm{x_n},\norm{y_n}\leq 1$ 
for each $n$ it follows that $\norm{x_n}, \norm{y_n}\to 1$. Because of $x_n-y_n\to z$ this implies $\norm{x_n-z}\to 1$.\par
It further follows from $\norm{x_n+y_n}\to 2$ and $x_n-y_n\to z$ that $\norm{2x_n-z}\to 2$.\par
Now as in the proof of \cite{dayjamesswaminathan}*{Theorem 1} we put
\begin{equation*}
\omega_n=\min\set*{1,\norm{x_n-z}^{-1}},\ a_n=\omega_n x_n,\ b_n=\omega_n(x_n-z) \ \ \forall n\in \N
\end{equation*}
and observe that $\norm{a_n}, \norm{b_n}\leq 1$ and $a_n-b_n=\omega_n z$ for each $n$ as well as $\omega_n\to 1$, $\norm{a_n+b_n}\to 2$
and $x_n^*(a_n)\to 1$.\par
Also as in the proof of \cite{dayjamesswaminathan}*{Theorem 1} we fix to sequences $(\alpha_n)_{n\in \N}$ and $(\beta_n)_{n\in \N}$ of 
non-negative real numbers such that 
\begin{equation*}
u_n:=a_n+\alpha_n z\in S_X \ \mathrm{and} \ v_n:=b_n-\beta_n z\in S_X \ \ \forall n\in \N.
\end{equation*}
Then $u_n-v_n=(\omega_n+\alpha_n+\beta_n)z$ for all $n\in \N$ and again as in the proof of \cite{dayjamesswaminathan}*{Theorem 1} one can
show that $\norm{u_n+v_n}\to 2$.\par
Because of $x_n^*(z)>0$ for sufficiently large $n$ it follows from $x_n^*(a_n)\to 1$ that $x_n^*(u_n)\to 1$. Since $X$ is a \ac{uacsed} space
it follows that $x_n^*(v_n)\to 1$. But $\omega_n+\alpha_n+\beta_n\geq \omega_n\to 1$, so we must have $x_n^*(z)\to 0$.\par
Thus we have shown that in any case there is a subsequence of $(x_n^*(z))_{n\in \N}$ that converges to one and the same argument works if
we start with an arbitrary subsequence of $(x_n^*(z))_{n\in \N}$. Hence the whole sequence must be convergent to one.\par
The implication $\mathrm{(ii)} \Rightarrow \mathrm{(iii)}$ is trivial. To prove $\mathrm{(iii)} \Rightarrow \mathrm{(i)}$ take sequences 
$(x_n)_{n\in \N}$, $(y_n)_{n\in \N}$ in $S_X$ and $(\alpha_n)_{n\in \N}$ in $\R$ as well as $z\in X\sm\set*{0}$ such that $x_n-y_n=\alpha_n z$ 
for all $n$ and $\norm{x_n+y_n}\to 2$. Also, take a sequence $(x_n^*)_{n\in \N}$ in $S_{X^*}$ with $x_n^*(x_n)\to 1$.\par
Since $(\alpha_n)_{n\in \N}$ is bounded by $2/\norm{z}$, by passing to subsequence we may assume that $\alpha_n\to \alpha$ for some 
$\alpha\in \R$.\par
Hence $x_n-y_n\to \alpha z$ and thus (iii) implies $x_n^*(y_n)\to 1$.\par
$\mathrm{(iv)} \Rightarrow \mathrm{(i)}$ is trivial and $\mathrm{(ii)} \Rightarrow \mathrm{(iv)}$ is proved exactly as we have just proved
$\mathrm{(iii)} \Rightarrow \mathrm{(i)}$. Thus the equivalence of (i)---(iv) is established.\par
$\mathrm{(iii)} \Rightarrow \mathrm{(v)}$ is trivial as well and $\mathrm{(v)} \Rightarrow \mathrm{(iii)}$ can be proved using the 
Bishop--Phelps--Bollob\'as theorem like in the proof of Proposition \ref{prop:char-uacs-sluacs}.\par
Next we prove $\mathrm{(ii)} \Rightarrow \mathrm{(vi)}$. Take a sequence $(x_n)_{n\in \N}$ in $B_X$ and an element $z\in X$ as well as
a sequence $(x_n^*)_{n\in \N}$ in $S_{X^*}$ such that
\begin{equation*}
2^{p-1}\paren*{\norm{x_n+z}^p+\norm{x_n}^p}-\norm{2x_n+z}^p\to 0 \ \mathrm{and} \ x_n^*(x_n)\to 1.
\end{equation*}
It follows that $\norm{x_n}\to 1$. As in the proof of \cite{dayjamesswaminathan}*{Theorem 1} we can make use of the inequality
\begin{equation*}
(a+b)^p+(a-b)^p\leq 2^{p-1}(a^p+b^p) \ \ \forall a\geq b\geq 0, \forall p\geq 2
\end{equation*}
to infer that $\norm{x_n+z}\to 1$ and $\norm{2x_n+z}\to 2$.\par
If we put $y_n=(x_n+z)/\norm{x_n+z}$ then $y_n\in S_X$, $\norm{x_n+y_n}\to 2$ and $x_n-y_n\to -z$, so (ii) implies $x_n^*(z)\to 0$
and we are done.\par
For the prove of $\mathrm{(vi)} \Rightarrow \mathrm{(ii)}$ fix two sequences $(x_n)_{n\in \N}$, $(y_n)_{n\in \N}$ in $B_X$ such that 
$\norm{x_n+y_n}\to 2$ and $x_n-y_n\to z\in X$ as well as a sequence $(x_n^*)_{n\in \N}$ in $S_{X^*}$ with $x_n^*(x_n)\to 1$.\par
It follows that $\norm{x_n}, \norm{y_n}\to 1$, $\norm{x_n-z}\to 1$ and $\norm{2x_n-z}\to 2$. Hence 
\begin{equation*}
2^{p-1}\paren*{\norm{x_n-z}^p+\norm{x_n}^p}-\norm{2x_n-z}^p\to 0
\end{equation*}
and (vi) implies $x_n^*(z)\to 0$.\par
The equivalence of (v) and (vii) can be proved analogously.\par
$\mathrm{(viii)} \Rightarrow \mathrm{(ix)}$ is trivial and $\mathrm{(vi)} \Rightarrow \mathrm{(viii)}$ is also obvious.
Let us finally prove $\mathrm{(ix)} \Rightarrow \mathrm{(vii)}$. If $(x_n)_{n\in \N}$ is a sequence in $S_X$, 
$(x_n^*)_{n\in \N}$ a sequence $S_{X^*}$ and $z\in X$ such that
\begin{equation*}
\norm{2x_n+z}^p-2^{p-1}\norm{x_n+z}^p\to 2^{p-1} \ \mathrm{and} \ x_n^*(x_n)=1 \ \forall n\in \N
\end{equation*}
then as before we can deduce that $\norm{x_n+z}\to 1$ and $\norm{2x_n+z}\to 2$. Thus $(x_n)_{n\in \N}$, $(y_n)_{n\in \N}:=(x_n+z)_{n\in \N}$,
$(x_n^*)_{n\in \N}$ and $-z$ meet the conditions of $(ix)$ and hence $x_n^*(z)\to 0$.
\end{Proof}

Let us also mention the following characterisation of the property \ac{uacsed} in terms of the space $X$ itself only. The proof is completely
analogous to the one for Proposition \ref{prop:char uacs without dual}.
\begin{proposition}\label{prop:char uacsed without dual}
For a Banach space $X$ the following assertions are equivalent.
\begin{enumerate}[\upshape(i)]
\item $X$ is \ac{uacsed}.
\item For every $z\in X\sm\set*{0}$ and every $\eps>0$ there exists some $\delta>0$ such that for every $t\in [0,\delta]$ and all $x,y\in S_X$
with $\norm{x+y}\geq 2(1-t)$ and $x-y\in \lin\set*{z}$ we have
\begin{equation*}
\norm{x+ty}+\norm{x-ty}\leq 2+\eps t.
\end{equation*}
\item For every $z\in X\sm\set*{0}$ and every $\eps>0$ there exists some $\delta>0$ such that for every $t\in [0,\delta]$ and all $x,y\in S_X$
with $\norm{x+y}\geq 2-\delta t$ and $x-y\in \lin\set*{z}$ we have
\begin{equation*}
\norm{x-ty}\leq 1+t(\eps-1).
\end{equation*}
\item For every $z\in X\sm\set*{0}$ there exists some $1\leq p<\infty$ such that for every $\eps>0$ there is some $\delta>0$ such that for all
$t\in [0,\delta]$ and all $x,y\in S_X$ with $\norm{x+y}\geq 2(1-t)$ and $x-y\in \lin\set*{z}$ we have
\begin{equation*}
\norm{x+ty}^p+\norm{x-ty}^p\leq 2+\eps t^p.
\end{equation*}
\item For every $z\in X\sm\set*{0}$ there exists some $1\leq p<\infty$ such that for every $\eps>0$ there is some $\delta>0$ such that for all
$t\in [0,\delta]$ and all $x,y\in S_X$ with $\norm{x+y}\geq 2-t\delta$ and $x-y\in \lin\set*{z}$ we have
\begin{equation*}
(1+t)^p+\norm{x-ty}^p\leq 2+\eps t^p.
\end{equation*}
\end{enumerate}
\end{proposition}
If we use the characterisation for \ac{uacsed} spaces with the condition $x_n-y_n\to z$ instead of $x_n-y_n\in \lin\set*{z}$ that was given
above, we also see that Proposition \ref{prop:char uacsed without dual} still holds true if we replace the condition $x-y\in \lin\set*{z}$
by $\norm{x-y-z}\leq\delta$ in the assertions (ii)--(v).\par
Next we consider quotient spaces. It is known that the quotient of a Banach space which is \ac{URED} by a finite-dimensional subspace 
is again \ac{URED} (cf. \cite{smith0}*{Remark before Problem 2}). By the same method of proof we can obtain the analogous result for 
the property \ac{uacsed}.
\begin{proposition}\label{prop:quot uacsed}
Let $X$ be a Banach space which is \ac{uacsed} and $U\ssq X$ a finite-dimensional subspace. Then $X/U$
is also \ac{uacsed}.
\end{proposition}

\begin{Proof}
As was implicitly mentioned in \cite{smith0}, if $W\ssq X/U$ is compact and $U$ is finite-dimensional then $\omega^{-1}(W)\cap(2B_X)$ 
is also compact, where $\omega$ is the canonical quotient map. For, if $(x_n)_{n\in \N}$ is a sequence in $\omega^{-1}(W)\cap(2B_X)$ 
then by compactness of $W$ we can pass to a subsequence such that $\omega(x_n)\to \omega(x)$ for some $x\in X$. Next fix a sequence 
$(y_n)_{n\in \N}$ in $U$ such that $\norm{x_n-x-y_n}\to 0$. It follows that $(y_n)_{n\in \N}$ is bounded and hence we can pass to a 
further subsequence such that $y_n\to y\in U$. Then $x_n\to x+y$ and since $\omega^{-1}(W)\cap(2B_X)$ is closed we must have 
$x+y\in \omega^{-1}(W)\cap(2B_X)$, so $\omega^{-1}(W)\cap(2B_X)$ is compact.\par
Now let us take two sequences $(z_n)_{n\in \N}$, $(w_n)_{n\in \N}$ in $S_{X/U}$ and an element $z\in X/U$ such that $\norm{z_n+w_n}\to 2$ 
and $z_n-w_n\to z\in X/U$. Also, fix a sequence $(\varphi_n)_{n\in \N}$ in $S_{(X/U)^*}$ with $\varphi_n(z_n)\to 1$.\par
Then $x_n^*:=\varphi_n\circ\omega\in S_{U^\perp}$ for all $n\in \N$. Since $U$ is in particular reflexive, we can find sequences 
$(x_n)_{n\in \N}$ and $(y_n)_{n\in \N}$ in $B_X$ such that $\omega(x_n)=z_n$ and $\omega(y_n)=w_n$ for each $n$ (cf. the proof of 
Proposition \ref{prop:quot acs}).\par
It follows that $\norm{x_n+y_n}\to 2$ and $x_n^*(x_n)\to 1$.\par
The set $W:=\set*{z_n-w_n:n\in \N}\cup\set*{z}$ is compact, hence by the introductory observation $\omega^{-1}(W)\cap(2B_X)$ is also compact
and so we can find a subsequence $(x_{n_k}-y_{n_k})_{k\in \N}$ that is convergent in $X$.\par
Since $X$ is \ac{uacsed} Proposition \ref{prop:char uacsed} implies $x_{n_k}^*(y_{n_k})\to 1$. Hence $\varphi_{n_k}(z)\to 0$.\par
The same argument shows that every subsequence of $(\varphi_n(z))_{n\in \N}$ possesses a subsubsequence which converges to zero, so we
have $\varphi_n(z)\to 0$ and hence by Proposition \ref{prop:char uacsed} the quotient space $X/U$ is also \ac{uacsed}.
\end{Proof}

Concerning absolute sums of \ac{uacsed} spaces, we have the following result (under the same general hypothesis as in 
section \ref{sec:abs sums}).
\begin{theorem}\label{thm:sums uacsed}
If $E$ is \ac{URED} and $(X_i)_{i\in I}$ is a family of \ac{uacsed} spaces then 
$\bigl[\bigoplus_{i\in I}X_i\bigr]_E$ is also \ac{uacsed}.
\end{theorem}

\begin{Proof}
The proof is analogous to the one of Smith's result on products of \ac{URED} spaces from \cite{smith1}.\par
Fix a non-zero element $z$ of $\bigl[\bigoplus_{i\in I}X_i\bigr]_E$ and two sequences $(x_n)_{n\in \N}$ and $(y_n)_{n\in \N}$ in 
the unit sphere of $\bigl[\bigoplus_{i\in I}X_i\bigr]_E$ such that $\norm{x_n+y_n}_E\to 2$ and $x_n-y_n\in \lin\set*{z}$, say 
$x_n-y_n=\alpha_n z$ for each $n\in \N$. Also, take a sequence $(x_n^*)_{n\in \N}$ in the dual unit sphere of 
$\bigl[\bigoplus_{i\in I}X_i\bigr]_E$ such that $x_n^*(x_n)\to 1$. As usual we write $x_n=(x_{n,i})_{i\in I}$,
$y_n=(y_{n,i})_{i\in I}$ and $x_n^*=(x_{n,i}^*)_{i\in I}$ for each $n\in \N$, as well as $z=(z_i)_{i\in I}$ and
as usual we conclude
\begin{equation}\label{eq:102}
\lim_{n\to \infty}\paren*{x_{n,i}^*(x_{n,i})-\norm{x_{n,i}^*}\norm{x_{n,i}}}=0 \ \ \forall i\in I.
\end{equation}
As in the proof from \cite{smith1} we put $f_n(i)=\norm{x_{n,i}}$ and $g_n^{(\beta)}(i)=\norm{x_{n,i}-\beta\alpha_n z_i}$ for all
$n\in \N$, all $i\in I$ and every $\beta\in \set*{1/2,1}$.\par
Then $\norm{f_n}_E=\norm{g_n^{(1)}}_E=1$ for every $n$ and 
\begin{equation}\label{eq:103}
\norm{g_n^{(1/2)}}_E=\frac{1}{2}\norm{x_n+y_n}_E\to 1.
\end{equation}
We also have 
\begin{equation}\label{eq:104}
\norm{f_n+g_n^{(1)}}_E=\norm*{\paren*{\norm{x_{n,i}}+\norm{y_{n,i}}}_{i\in I}}_E\to 2,
\end{equation}
as we have shown many times before, and furthermore
\begin{align*}
&1+\norm{g_n^{(1/2)}}_E\geq\norm{f_n+g_n^{(1/2)}}_E=\frac{1}{2}\norm*{\paren*{2\norm{x_{n,i}}+\norm{x_{n,i}+y_{n,i}}}_{i\in I}}_E \\
&\geq \frac{1}{2}\paren*{\norm*{\paren*{2\norm{x_{n,i}}+2\norm{y_{n,i}}+\norm{x_{n,i}+y_{n,i}}}_{i\in I}}_E-2} \\
&\geq \frac{1}{2}\paren*{3\norm{x_n+y_n}_E-2},
\end{align*}
hence
\begin{equation}\label{eq:105}
\norm{f_n+g_n^{(1/2)}}_E\to 2.
\end{equation}
Note that $\abs{f_n(i)-g_n^{(\beta)}(i)}\leq \beta\abs{\alpha_n}\norm{z_i}\leq 2\norm{z}_E^{-1}\norm{z_i}$ for all $n\in \N$, all
$i\in I$ and every $\beta\in \set*{1/2,1}$. Remember that we always assume that $\lin\set*{e_i:i\in I}$ is dense in $E$. This easily
implies that for each $f\in E$ the set $\set*{g\in \R^I: \abs{g(i)}\leq \abs{f(i)} \ \forall i\in I}$ is a compact subset of $E$.
Hence we can find a strictly increasing sequence $(n_k)_{k\in \N}$ in $\N$, two elements $h^{(1)}$ and $h^{(1/2)}$ of $E$ and
an $\alpha\in \R$ such that 
\begin{equation}\label{eq:106}
\norm{f_{n_k}-g_{n_k}^{(\beta)}-h^{(\beta)}}_E\to 0 \ \ \forall \beta\in\set*{1/2,1} \ \mathrm{and} \ \alpha_{n_k}\to \alpha.
\end{equation}
Since $E$ is \ac{URED} it follows from \eqref{eq:103}, \eqref{eq:104}, \eqref{eq:105} and \eqref{eq:106} together with 
\cite{dayjamesswaminathan}*{Theorem 1} that $h^{(\beta)}=0$ for $\beta\in\set*{1/2,1}$, thus
\begin{equation}\label{eq:107}
\lim_{k\to \infty}\paren*{\norm{x_{n_k,i}}-\norm{x_{n_k,i}-\beta\alpha_{n_k}z_i}}=0 \ \ \forall i\in I, \forall \beta\in\set*{1/2,1}.
\end{equation}
Now let us fix an arbitrary $i_0\in I$. We can find a subsequence $(\norm{x_{n_{k_j},i_0}})_{j\in \N}$ that is convergent to some 
$a\in \R$. From \eqref{eq:107} we get that 
\begin{equation*}
\norm{x_{n_{k_j},i_0}-\beta\alpha_{n_{k_j}}z_{i_0}}\to a \ \ \forall \beta\in\set*{1/2,1},
\end{equation*}
hence
\begin{equation*}
\norm{x_{n_{k_j},i_0}}, \norm{y_{n_{k_j},i_0}}\to a \ \mathrm{and} \ \norm{x_{n_{k_j},i_0}+y_{n_{k_j},i_0}}\to 2a.
\end{equation*}
Together with \eqref{eq:102}, $x_{n_{k_j},i_0}-y_{n_{k_j},i_0}\to \alpha z_{i_0}$ and Proposition \ref{prop:char uacsed} this easily implies 
$x_{n_{k_j},i_0}^*(\alpha z_{i_0})\to 0$.\par
The same argument works if we start with an arbitrary subsequence of $(\norm{x_{n_k,i_0}})_{k\in \N}$ thus we have
\begin{equation}\label{eq:108}
\lim_{k\to \infty} x_{n_k,i}^*(\alpha z_i)=0 \ \ \forall i\in I.
\end{equation}
Now fix an arbitrary $\eps>0$ and a finite subset $J\ssq I$ such that
\begin{equation*}
\norm*{\sum_{i\in J}\alpha\norm{z_i}-\alpha(\norm{z_i})_{i\in I}}_E\leq\eps.
\end{equation*}
By \eqref{eq:108} we have for all sufficiently large $k$
\begin{equation*}
\abs*{\sum_{i\in J}x_{n_k,i}^*(\alpha z_i)}\leq\eps.
\end{equation*}
These two inequalities together easily imply $x_{n_k}^*(\alpha z)\leq 2\eps$ (for sufficiently large $k$). Thus we have
$x_{n_k}^*(\alpha z)\to 0$ and hence $x_{n_k}^*(y_{n_k})\to 1$.\par
Again, the same argument works if we start with an arbitrary subsequence of $(x_n^*(y_n))_{n\in \N}$, thus we must 
have $x_n^*(y_n)\to 1$ and the proof is finished.
\end{Proof}

\begin{corollary}\label{cor:p-sum uacsed}
If $(X_i)_{i\in I}$ is a family of \ac{uacsed} spaces then $\bigl[\bigoplus_{i\in I}X_i\bigr]_p$ is also
\ac{uacsed} for every $1<p<\infty$.
\end{corollary}

\section{A few examples}\label{sec:examples}
In this section we collect some examples showing that no other arrows can be drawn in the Figures 
\ref{fig:3}, \ref{fig:4} and \ref{fig:6} (except, of course, for combinations of two or more existing arrows).\par
In what follows, we use the notation $x^\prime=(0,x(2),x(3),\dots)$ for every $x\in \R^{\N}$.

\begin{example}\label{ex:uacs not R not S}
{\em A \ac{uacs} space which is neither \ac{R} nor \ac{S}}.\par
\noindent We can simply take a norm on $\R^2$ which is neither \ac{R} nor \ac{S} but still \ac{acs} 
(since the space is finite-dimensional it will then be even \ac{uacs}). The unit ball  
of one such norm is sketched below.

\begin{figure}[H]
\begin{center}
\begin{tikzpicture}
  \draw[help lines] (-3.5,0)--(3.5,0);
  \draw[help lines] (0,-2)--(0,2);
  \draw (-3,0) to [out=65,in=180] (-1.5,1.5);
  \draw (-1.5,1.5)--(1.5,1.5);
  \draw (1.5,1.5) to [out=0,in=115] (3,0);
  \draw (-3,0) to [out=295,in=180] (-1.5,-1.5);
  \draw (-1.5,-1.5)--(1.5,-1.5);
  \draw (1.5,-1.5) to [out=0,in=245] (3,0);
\end{tikzpicture}
\CAP\label{fig:9}
\end{center}
\end{figure}

\noindent The norm is not rotund, because the unit sphere contains line segments, but the endpoints of these
line segments are smooth points of the unit ball. On the other hand, the norm has also non-smooth points,
but they are not the endpoints of any line segments on the unit sphere, so on the whole the space is
\ac{acs} but neither \ac{R} nor \ac{S}.
\end{example}

\begin{example}\label{ex:LUR URED not wuacs}
{\em A space which is \ac{LUR} and \ac{URED} but not \ac{wuacs}}.\par
\noindent In \cite{smith}*{Example 6} Smith defines an equivalent norm on $\ell^1$ as follows:
\begin{equation*}
\Norm{x}^2=\norm{x}_1^2+\norm{x}_2^2 \ \ \forall x\in \ell^1.
\end{equation*}
He shows that $(\ell^1,\Norm{\,.\,})$ is \ac{LUR} and \ac{URED} but not \ac{WUR}. In fact, $(\ell^1,\Norm{\,.\,})$ 
is not even \ac{wuacs}. To see this, put
\begin{align*}
&\beta_n=\frac{2}{\sqrt{4n^2+2n}}, \\
&x_n=(\underbrace{\beta_n,0,\beta_n,0,\dots,\beta_n,0}_{2n},0,0,\dots), \\
&y_n=(\underbrace{0,\beta_n,0,\beta_n,\dots,0,\beta_n}_{2n},0,0,\dots)
\end{align*}
for every $n\in \N$. Then it is easily checked that $\Norm{x_n+y_n}=2$ for every $n\in \N$ and
$\Norm{x_n}=\Norm{y_n}\to 1$.\par
Now let $x^*$ be the functional on $\ell^1$ represented by $(1,0,1,0,\dots)\in \ell^{\infty}$.
Then $\Norm{x^*}^*\leq1$ (where $\Norm{\,.\,}^*$ denotes the dual norm of $\Norm{\,.\,}$) and 
it is easy to see that $x^*(x_n)\to 1$. On the other hand, $x^*(y_n)=0$ for every $n\in \N$ thus
$(\ell^1,\Norm{\,.\,})$ is not \ac{wuacs}.
\end{example}

\begin{example}\label{ex:WUR not msluacs}
{\em A space which is \ac{WUR} and \ac{URED} but not \ac{msluacs}}.\par
\noindent The following equivalent norm $\Norm{\,.\,}$ on $\ell^2$ was also defined in \cite{smith}*{Example 2}.
Fix a sequence $(\alpha_n)_{n\in \N}$ in $]0,1]$ which decreases to 0 and define $T:\ell^2 \rightarrow \ell^2$
by $Tx=(x(1),\alpha_2x(2),\alpha_3x(3),\dots)$ for every $x\in \ell^2$. Then put
\begin{equation*}
\Norm{x}^2=\max\set{\abs{x(1)},\norm{x^\prime}_2}^2+\norm{Tx}_2^2 \ \ \forall x\in \ell^2.
\end{equation*}
It is shown in \cite{smith}*{Example 2} that $(\ell^2,\Norm{\,.\,})$ is \ac{WUR} and \ac{URED} but not \ac{MLUR}.\par
To prove the latter, Smith defines $\alpha=1/\sqrt{2}$, $x=\alpha e_1$, $x_n=\alpha(e_1+e_n)$ and 
$y_n=\alpha(e_1-e_n)$ for every $n$ and observes that $\Norm{x}=1$, $\Norm{x_n}=\Norm{y_n}\to 1$ and 
$x_n+y_n=2x$ for every $n$, but $\Norm{x_n-y_n}\to \sqrt{2}$.\par
If $x_n^*$ denotes the functional on $\ell^2$ represented by $\alpha(e_1+e_n)\in \ell^2$ for 
every in $n\in \N$ then we have $\Norm{x_n^*}^*\leq 1$ and $x_n^*(x_n)=1$ for every $n$ but $x_n^*(y_n)=0$
for every $n$, hence $(\ell^2,\Norm{\,.\,})$ is not even \ac{msluacs}.
\end{example}

\begin{example}\label{ex:R not mluacs}
{\em A Banach space which is \ac{R} but not \ac{mluacs}}.\par
\noindent This example is a slight modification of \cite{fabian}*{Exercise 8.52}. We define an equivalent 
norm on $\ell^1$ by
\begin{equation*}
\Norm{x}=\max\set{\abs{x(1)},\norm{x^\prime}_1}+\norm{x}_2 \ \ \forall x\in \ell^1.
\end{equation*}
Using the fact that $(\ell^2,\norm{\,.\,}_2)$ is \ac{R} it is easy to see that $(\ell^1,\Norm{\,.\,})$
is also \ac{R}. To see that $(\ell^1,\Norm{\,.\,})$ is not \ac{mluacs}, put $x=e_1$ and
\begin{align*}
&x_n=\Big(1,\underbrace{\frac{1}{n},\dots,\frac{1}{n}}_{n},0,0,\dots\Big), \\
&y_n=\Big(1,\underbrace{-\frac{1}{n},\dots,-\frac{1}{n}}_{n},0,0,\dots\Big)
\end{align*}
for every $n\in \N$. Then it is easy to see that $\Norm{x}=2$, $\Norm{x_n}=\Norm{y_n}\to 2$ and $x_n+y_n=2x$
for every $n$.\par
If $x^*$ is the functional on $\ell^1$ represented by $(1,1,\dots)\in \ell^{\infty}$ then $\Norm{x^*}^*\leq 1$
and $x^*(x_n)=2$ for every $n$, but $x^*(y_n)=0$ for every $n$, hence $(\ell^1,\Norm{\,.\,})$ is not \ac{mluacs}.
\end{example}

\begin{example}\label{ex:MLUR not luacs}
{\em A Banach space which is \ac{MLUR} but not \ac{luacs}}.\par
\noindent We define an equivalent norm $\norm{\,.\,}_M$ on $\ell^1$ as follows:
\begin{align*}
&\Norm{x}=\norm{x^\prime}_1+\norm{x}_2 \ \mathrm{and} \\
&\norm{x}_M^2=\norm{x}_1^2+\norm{x^\prime}_2^2+\Norm{x}^2 \ \ \forall x\in \ell^1.
\end{align*}
Then we have $\sqrt{2}\norm{x}_1\leq\norm{x}_M\leq\sqrt{6}\norm{x}_1$ for every $x\in \ell^1$.\par
We first show that $(\ell^1,\norm{\,.\,}_M)$ is not \ac{luacs}. To do so, we put $x=e_1$ and 
\begin{equation*}
x_n=\Big(0,\underbrace{\frac{1}{n},\dots,\frac{1}{n}}_{n},0,0,\dots\Big) \ \ \forall n\in \N.
\end{equation*}
Then it is easy to calculate $\norm{x}_M=\sqrt{2}$, $\norm{x_n}_M\to \sqrt{2}$ and $\norm{x_n+x}_M\to 2\sqrt{2}$.
Let $x^*$ be the functional on $\ell^1$ represented by $(0,\sqrt{2},\sqrt{2},\dots)\in \ell^{\infty}$. Then
$\norm{x^*}_M^*\leq1$ and $x^*(x_n)=\sqrt{2}$ for every $n\in \N$, but $x^*(x)=0$, thus $(\ell^1,\norm{\,.\,}_M)$ is not 
\ac{luacs}.\par
Now we prove that $(\ell^1,\norm{\,.\,}_M)$ is \ac{MLUR}. So let us fix two sequences $(x_n)_{n\in \N}, (y_n)_{n\in \N}$
in $\ell^1$ and $x\in \ell^1$ such that $\norm{x_n}_M=\norm{x}_M=\norm{y_n}_M=1$ for every $n\in \N$ and 
$\norm{x_n+y_n-2x}_M\to 0$.\par
Then we also have
\begin{equation}\label{eq:109}
\norm{x_n+y_n-2x}_1, \norm{x_n^\prime+y_n^\prime-2x^\prime}_2, \Norm{x_n+y_n-2x}\to 0
\end{equation}
and hence
\begin{equation}\label{eq:110}
\norm{x_n+y_n}_1\to 2\norm{x}_1, \norm{x_n^\prime+y_n^\prime}_2\to 2\norm{x^\prime}_2, \Norm{x_n+y_n}\to 2\Norm{x}.
\end{equation}
We further have
\begin{align*}
&\norm{x_n+y_n}_M\leq \paren*{(\norm{x_n}_1+\norm{y_n}_1)^2+(\norm{x_n^\prime}_2+\norm{y_n^\prime}_2)^2+(\Norm{x_n}+\Norm{y_n})^2}^{\frac{1}{2}} \\
&\leq \norm{x_n}_M+\norm{y_n}_M=2
\end{align*}
and because of $\norm{x_n+y_n}\to 2$ and the uniform rotundity of the euclidean norm on $\R^3$ this implies
\begin{align}
\label{eq:111}
&\paren*{\norm{x_n}_1+\norm{y_n}_1}^2-\norm{x_n+y_n}_1^2\to 0, \\
\label{eq:112}
&\paren*{\norm{x_n^\prime}_2+\norm{y_n^\prime}_2}^2-\norm{x_n^\prime+y_n^\prime}_2^2\to 0, \\
\label{eq:113}
&\paren*{\Norm{x_n}+\Norm{y_n}}^2-\Norm{x_n+y_n}^2\to 0
\end{align}
and 
\begin{equation}\label{eq:114}
\norm{x_n}_1-\norm{y_n}_1, \norm{x_n^\prime}_2-\norm{y_n^\prime}_2, \Norm{x_n}-\Norm{y_n}\to 0.
\end{equation}
Combining \eqref{eq:111}, \eqref{eq:112}, \eqref{eq:113}, \eqref{eq:110}  and \eqref{eq:114} we get that
\begin{equation}\label{eq:115}
\norm{x_n}_1,\norm{y_n}_1\to \norm{x}_1, \norm{x_n^\prime}_2, \norm{y_n^\prime}_2\to \norm{x^\prime}_2, \Norm{x_n},\Norm{y_n}\to \Norm{x}.
\end{equation}
Since $(\ell^2,\norm{\,.\,}_2)$ is \ac{UR} we can deduce from \eqref{eq:109} and \eqref{eq:114} that
\begin{equation}\label{eq:116}
\norm{x_n^\prime-x^\prime}_2, \norm{y_n^\prime-x^\prime}_2\to 0.
\end{equation}
By \eqref{eq:115} and the definition of $\Norm{\,.\,}$ we have $\norm{x_n^\prime}_1+\norm{x_n}_2\to \norm{x^\prime}_1+\norm{x}_2$,
which together with $\norm{x_n}_1\to \norm{x}_1$ implies
\begin{equation}\label{eq:117}
\abs{x_n(1)}-\norm{x_n}_2\to \abs{x(1)}-\norm{x}_2.
\end{equation}
Let us put $a_n=\abs{x_n(1)}$ and $b_n=\norm{x_n^\prime}_2$ for every $n$ as well as $a=\abs{x(1)}$ and $b=\norm{x^\prime}_2$.\par
Then \eqref{eq:117} reads
\begin{equation}\label{eq:118}
a_n-\sqrt{a_n^2+b_n^2}\to a-\sqrt{a^2+b^2}
\end{equation}
and by \eqref{eq:115} we have $b_n\to b$.\par
If $b\neq 0$ this easily implies $a_n\to a$. If $b=0$ then $x=x(1)e_1$ and because of $\norm{x}_M=1$ it follows
$\abs{x(1)}=1/\sqrt{2}$. But by \eqref{eq:109} we have $\abs{x_n(1)+y_n(1)}\to 2\abs{x(1)}=\sqrt{2}$ and since
$\abs{x_n(1)}\leq\norm{x_n}_1\leq\norm{x_n}_M/\sqrt{2}=1/\sqrt{2}$ (and likewise $\abs{y_n(1)}\leq1/\sqrt{2}$) for every $n$
it follows that $\abs{x_n(1)}, \abs{y_n(1)}\to 1/\sqrt{2}$.\par
Thus we have $\abs{x_n(1)}\to \abs{x(1)}$ in any case and analogously we can show that we always have $\abs{y_n(1)}\to \abs{x(1)}$.\par
Because of $x_n(1)+y_n(1)\to 2x(1)$ this implies $x_n(1)\to x(1)$ and $y_n(1)\to y(1)$. Taking into account \eqref{eq:116} it follows
\begin{equation}\label{eq:119}
x_n(i)\to x(i) \ \mathrm{and} \ y_n(i)\to x(i) \ \ \forall i\in \N.
\end{equation}
By \eqref{eq:115} we also have $\norm{x_n}_1, \norm{y_n}_1\to \norm{x}_1$ and it is well known that these two conditions together
imply $\norm{x_n-x}_1, \norm{y_n-x}_1\to 0$. Hence we have $\norm{x_n-y_n}_M\to 0$, as desired.
\end{example}

In the next section, we will use the notions of \ac{acs} and \ac{luacs} spaces to obtain some prohibitive results on
spaces with the so called alternative Daugavet-property.

\section[Prohibitive results for the alternative Daugavet property]{Prohibitive results for real Banach spaces with the alternative Daugavet property}\label{sec:alt daugavet}
Recall that a (real or complex) Banach space $X$ is said to have the \ac{aDP} if the so called \ac{aDE}
\begin{equation}\label{aDE}
\max_{\omega\in \T}\norm{\id+\omega T}=1+\norm{T} \tag{aDE}
\end{equation}
holds for every rank-1-operator $T\in L(X)$, where $\T$ denotes the set of all scalars of modulus one. This notion
was originally introduced in \cite{martin}. For more information and background on the \ac{aDP} and its relation to
the concept of numerical index of Banach spaces we refer the reader to \cite{martin}, \cite{kadetsmartin} and
references therein.\par
In \cite{kadetsmartin} it is shown that a Banach space with the \ac{aDP} (and its dual) cannot have certain rotundity
or smoothness properties, more precisely it is shown that if $X$ has the \ac{aDP} and dimension greater than one, then 
$X^*$ is neither \ac{R} nor \ac{S} (cf. \cite{kadetsmartin}*{Theorem 2.1}) and the unit ball of $X$ has no \ac{WLUR} points 
(cf. \cite{kadetsmartin}*{Proposition 2.4}). A point $x\in S_X$ is called a \ac{WLUR} point of $B_X$ if for every sequence 
$(x_n)_{n\in \N}$ in $S_X$ the condition $\norm{x_n+x}\to 2$ implies that $(x_n)_{n\in \N}$ converges weakly to $x$.
One can easily generalise these results to the \ac{acs} resp. \ac{luacs} case (for real Banach spaces), where the proofs 
stay almost exactly the same.
\begin{theorem}\label{thm:dual of aDP not acs}
Let $X$ be a real Banach space with the \ac{aDP} of dimension strictly greater than one. Then $X^*$ is not \ac{acs}.
\end{theorem}

\begin{Proof}
Exactly as in the proof of \cite{kadetsmartin}*{Theorem 2.1} we can find a sequence $(T_n)_{n\in \N}$ of compact norm-one operators 
from $X$ into $X$ such that $T_n^2\to 0$, a sequence $(\lambda_n)_{n\in \N}$ in $\T=\set{\pm 1}$ and sequences $(x_n^*)_{n\in \N}$
and $(x_n^{**})_{n\in \N}$ in $S_{X^*}$ and $S_{X^{**}}$ respectively, such that
\begin{equation*}
\lambda_nx_n^{**}(x_n^*)=1 \ \mathrm{and} \ x_n^{**}(T_n^*x_n^*)=1 \ \ \forall n\in \N.
\end{equation*}
But this implies $\norm{\lambda_nx_n^*+T_n^*x_n^*}=2$ and $\abs{[T_n^{**}x_n^{**}](\lambda_nx_n^*)}=1$ for every $n$.\par
If $X^*$ was \ac{acs} this would imply $\abs{x_n^{**}([T_n^*]^2x_n^*)}=\abs{[T_n^{**}x_n^{**}](T_n^*x_n^*)}=1$ and hence
$\norm{[T_n^*]^2}\geq 1$ for every $n\in \N$. But on the other hand we have $[T_n^*]^2\to 0$, a contradiction.
\end{Proof}

We call a point $x\in S_X$ a \ac{luacs} point of $B_X$ if for every sequence $(x_n)_{n\in \N}$ in $S_X$ and every $x^*\in S_{X^*}$
the two conditions $\norm{x_n+x}\to 2$ and $x^*(x_n)\to 1$ imply $x^*(x)=1$. Then we have the following result.
\begin{proposition}\label{prop:aDP no luacs points}
Let $X$ be a real Banach space with the \ac{aDP} of dimension strictly greater than one. Then $B_X$ has no \ac{luacs} points.
\end{proposition}

\begin{Proof}
Suppose $x\in S_X$ is a \ac{luacs} point of $B_X$. Exactly as in the proof of \cite{kadetsmartin}*{Proposition 2.4} we can
find $x^*\in S_{X^*}$ with $x^*(x)=0$ and sequences $(x_n)_{n\in \N}$ and $(x_n^*)_{n\in \N}$ in $S_X$ and $S_{X^*}$ respectively, 
such that 
\begin{equation*}
x_n^*(x_n)=1 \ \ \forall n\in \N \ \mathrm{and} \ \abs{x_n^*(x)},\abs{x^*(x_n)}\to 1.
\end{equation*}
Again as in the proof of \cite{kadetsmartin}*{Proposition 2.4} we take $\lambda_n\in \set{\pm1}$ with $x_n^*(x)=\lambda_n\abs{x_n^*(x)}$ 
for every $n$ and conclude that $\norm{x+\lambda_nx_n}\to 2$.\par
Because $\abs{x^*(\lambda_nx_n)}\to 1$ and since $x$ is a \ac{luacs} point of $B_X$ it follows that $\abs{x^*(x)}=1$, a contradiction.
\end{Proof}

\address
\email


\begin{bibdiv}
\begin{biblist}

\bib{abramovich}{article}{
  title={The Daugavet Equation in Uniformly Convex Banach Spaces},
  author={Abramovich, Y.A.},
  author={Aliprantis, C.D.},
  author={Burkinshaw, O.},
  journal={J. Funct. Anal.},
  volume={97},
  date={1991},
  pages={215--230},
  review={\mr{1105660}}
  }

\bib{albiac}{book}{
  title={Topics in Banach Space Theory},
  author={Albiac, F.},
  author={Kalton, N.J.},
  publisher={Springer},
  series={Graduate Texts in Mathematics},
  volume={233},
  date={2006},
  review={\mr{2192298}}
  }

\bib{anderson}{thesis}{
  title={Midpoint local uniform convexity, and other geometric properties of Banach spaces},
  author={Anderson, K.W.},
  organization={University of Illinois},
  type={Dissertation},
  date={1960}
  }

\bib{beauzamy}{book}{
  title={Introduction to Banach Spaces and their Geometry},
  author={Beauzamy, B.},
  publisher={North-Holland},
  edition={2},
  address={Amsterdam-New York-Oxford},
  date={1983}
  }

\bib{bollobas}{book}{
  title={Linear Analysis},
  subtitle={An Introductory Course},
  author={Bollob\'as, B.},
  publisher={Cambridge University Press},
  address={Cambridge--New York--Port Chester--Melbourne--Sydney},
  date={1990}
  }

\bib{clarkson}{article}{
  title={Uniformly convex spaces},
  author={Clarkson, J.A.},
  journal={Trans. Amer. Math. Soc.},
  volume={40},
  date={1936},
  pages={396--414}
  }

\bib{day1}{article}{
  title={Reflexive Banach spaces not isomorphic to uniformly convex spaces},
  author={Day, M.M.},
  journal={Bull. Amer. Math. Soc.},
  volume={47},
  number={4},
  date={1941},
  pages={313--317}
  }

\bib{day2}{article}{
  title={Some more uniformly convex spaces},
  author={Day, M.M.},
  journal={Bull. Amer. Math. Soc.},
  volume={47},
  number={6},
  date={1941},
  pages={504--507},
  review={\mr{0004068}}
  }

\bib{day3}{article}{
  title={Uniform Convexity III},
  author={Day, M.M.},
  journal={Bull. Amer. Math. Soc.},
  volume={49},
  number={10},
  date={1943},
  pages={745--750},
  review={\mr{0009422}}
  }

\bib{day-factorspaces}{article}{
  title={Uniform Convexity in Factor and Conjugate Sapces},
  author={Day, M.M.},
  journal={Ann. of Math.},
  volume={45},
  number={2},
  date={1944},
  pages={375--385},
  review={\mr{0010779}}
  }

\bib{day4}{article}{
  title={Strict Convexity and Smoothness of Normed Spaces},
  author={Day, M.M.},
  journal={Trans. Amer. Math. Soc.},
  volume={78},
  number={2},
  date={1955},
  pages={516--528},
  review={\mr{00067351}}
  }

\bib{dayjamesswaminathan}{article}{
  title={Normed linear spaces that are uniformly convex in every direction},
  author={Day, M.M.},
  author={James, R.C.},
  author={Swaminathan, S.},
  journal={Can. J. Math.},
  volume={23},
  number={6},
  date={1971},
  pages={1051--1059},
  review={\mr{0287285}}
  }

\bib{day5}{book}{
  title={Normed linear spaces},
  author={Day, M.M.},
  series={Ergebnisse der Mathematik und ihrer Grenz-gebiete},
  volume={21},
  edition={3},
  publisher={Springer},
  address={Berlin--Heidelberg--New York},
  date={1973},
  review={\mr{0344849}}
  }

\bib{deville}{book}{
  title={Smoothness and renormings in Banach spaces},
  author={Deville, R.},
  author={Godefroy, G.},
  author={Zizler, V.},
  series={Pitman Monongraphs and Surveys in Pure and Applied Mathematics},
  volume={64},
  publisher={Longman Scientific \& Technical},
  date={1993}
  }

\bib{dhompongsa0}{article}{
  title={On a generalized James constant},
  author={Dhompongsa, S.},
  author={Kaewkhao, A.},
  author={Tasena, S.},
  journal={J. Math. Anal. Appl.},
  volume={285},
  date={2003},
  pages={419--435}
  }

\bib{dhompongsa}{article}{
  title={Uniform smoothness and $U$-convexity of $\psi$-direct sums},
  author={Dhompongsa, S.},
  author={Kaewkhao, A.},
  author={Saejung, S.},
  journal={J. Nonlinear Convex Anal.},
  volume={6},
  number={2},
  date={2005},
  pages={327--338},
  review={\mr{2159843}}
  }

\bib{dowling1}{article}{
  title={The optimality of James's distortion theorems},
  author={Dowling, P.N.},
  author={Johnson, W.B.},
  author={Lennard, C.J.},
  author={Turett, B.},
  journal={Proc. Amer. Math. Soc.},
  volume={125},
  number={1},
  date={1997},
  pages={167--174},
  review={\mr{1346969}}
  }

\bib{dowling2}{article}{
  title={Extremal structure of the unit ball of direct sums of Banach spaces},
  author={Dowling, P.N.},
  author={Saejung, S.},
  journal={Nonlinear Analysis},
  volume={8},
  date={2008},
  pages={951--955},
  review={\mr{2382311}}
  }

\bib{dutta}{article}{
  title={Local $U$-Convexity},
  author={Dutta, S.},
  author={Lin, B.L.},
  journal={Journal of Convex Analysis},
  volume={18},
  number={3},
  date={2011},
  pages={811--821}
  }


\bib{fabian}{book}{
  title={Functional Analysis and Infinite-Dimensional Geometry},
  author={Fabian, M.},
  author={Habala, P.},
  author={H\'ajak, P.},
  author={Montesinos Santaluc\'{\i}a, V.},
  author={Pelant, J.},
  author={Zizler, V.},
  series={CMS Books in Mathematics},
  publisher={Springer},
  address={New York--Berlin--Heidelberg},
  date={2001}
  }

\bib{gao1}{article}{
  title={Normal structure and modulus of $u$-convexity in Banach spaces},
  author={Gao,J.},
  conference={
    title={Function Spaces, Differential Operators and Nonlinear Analysis},
    address={Paseky nad Jizerou},
    date={1995}
  },
  book={
    publisher={Prometheus},
    address={Prague},
    date={1996}
  },
  pages={195--199},
  review={\mr{1480939}}
  }

\bib{gao2}{article}{
  title={On two classes of Banach spaces with uniform normal structure},
  author={Gao, J.},
  author={Lau, K.S.},
  journal={Studia Math.},
  volume={99},
  number={1},
  date={1991},
  pages={41--56},
  review={\mr{1120738}}
  }

\bib{garkavi}{article}{
  title={The best possible net and best possible cross section of a set in a normed space},
  author={Garkavi, A.L.},
  journal={Izv. Akad. Nauk SSSR Ser. Mat.},
  volume={26},
  date={1962},
  pages={87--106},
  language={Russian},
  review={\mr{0136969}}
  }

\bib{goebel}{article}{
  title={Classical theory of nonexpansive mappings},
  author={Goebel, K.},
  author={Kirk, W.A.},
  book={
    title={Handbook of Metric Fixed Point Theory},
    editor={Kirk, W.A.},
    editor={Sims, B.},
    publisher={Kluwer Academic Publishers},
    address={Dordrecht--Boston--London},
    date={2001}
    },
  pages={49--91},
  review={\mr{1904274}}
  }

\bib{kadets82}{article}{
  title={Relation between some properties of convexity of the unit ball of a Banach space},
  author={Kadets, M.I.},
  journal={Funct. Anal. Appl.},
  volume={16},
  number={3},
  date={1982},
  pages={204--206}
  }

\bib{kadets}{article}{
  title={Banach spaces with the Daugavet property},
  author={Kadets, V.},
  author={Shvydkoy, R.},
  author={Sirotkin, G.},
  author={Werner, D.},
  journal={Trans. Amer. Math. Soc.},
  volume={352},
  number={2},
  date={2000},
  pages={855--873},
  review={\mr{1621757}}
  }

\bib{kadetsmartin}{article}{
  title={Convexity and smoothness of Banach spaces with numerical index one},
  author={Kadets, V.},
  author={Mart\'{i}n, M.},
  author={Mer\'{i}, J.},
  author={Pay\'a, R.},
  journal={Illinois J. Math.},
  volume={53},
  number={1},
  date={2009},
  pages={163--182},
  review={\mr{2584940}}
  }

\bib{klee}{article}{
  title={Some New Results on Smoothness and Rotundity in Normed Linear Spaces},
  author={Klee, V.},
  journal={Math. Ann.},
  volume={139},
  date={1959},
  pages={51--63},
  review={\mr{0115076}}
  }

\bib{lau}{article}{
  title={Best approximation by closed sets in Banach spaces},
  author={Lau, K.S.},
  journal={J. Approx. Theory},
  volume={23},
  date={1978},
  pages={29--36},
  review={\mr{0493114}}
  }

\bib{lee}{article}{
  title={Polynomial numerical indices of Banach spaces with absolute norms},
  author={Lee, H.J.},
  author={Mart\'{\i}n, M.},
  author={Mer\'{\i}, J.},
  journal={Linear Algebra and its Applications},
  volume={435},
  date={2011},
  pages={400--408},
  review={\mr{2782789}}
  }

\bib{lovaglia}{article}{
  title={Locally Uniformly Convex Banach Spaces},
  author={Lovaglia, A.R.},
  journal={Trans. Amer. Math. Soc.},
  volume={78},
  number={1},
  date={1955},
  pages={225--238}
  }

\bib{martin}{article}{
  title={An alternative Daugavet property},
  author={Mart\'{i}n, M.},
  author={Oikhberg, T.},
  journal={J. Math. Anal. Appl.},
  volume={294},
  number={1},
  date={2004},
  pages={158--180},
  review={\mr{2059797}}
  }

\bib{mazcunan-navvarro}{article}{
  title={On the modulus of $u$-convexity of Ji Gao},
  author={Mazcu\~{n}\'an-Navarro, E.M.},
  journal={Abstract and Applied Analysis},
  volume={2003},
  number={1},
  date={2003},
  pages={49-54},
  review={\mr{1954245}}
  }

\bib{mitani}{article}{
  title={$U$-convexity of $\psi$-direct sums of Banach spaces},
  author={Mitani, K.},
  journal={J. Nonlinear Convex Anal.},
  volume={11},
  number={2},
  date={2010},
  pages={199-213},
  review={\mr{2682863}}
  }

\bib{saejung}{article}{
  title={On the modulus of $U$-convexity},
  author={Saejung, S.},
  journal={Abstract and Applied Analysis},
  volume={2005},
  number={1},
  date={2005},
  pages={59--66},
  review={\mr{2142156}}
  }

\bib{sirotkin}{article}{
  title={New properties of Lebesgue-Bochner $L_p(\Omega,\Sigma,\mu;X)$ spaces},
  author={Sirotkin, G.G.},
  journal={Houston J. Math.},
  volume={27},
  number={4},
  date={2001},
  pages={897--906}
  }

\bib{smith0}{article}{
  title={Banach spaces that are uniformly rotund in weakly compact sets of directions},
  author={Smith, M.A.},
  journal={Can. J. Math.},
  volume={29},
  number={5},
  date={1977},
  pages={963--970},
  review={\mr{0450942}}
  }

\bib{smithsullivan}{article}{
  title={Extremely smooth Banach spaces},
  author={Smith, M.A.},
  author={Sullivan, F.},
  book={
    title={Banach Spaces of Analytic Functions},
    editor={Baker, J.},
    editor={Cleaver, C.},
    editor={Diestel, J.},
    series={Lecture Notes in Mathematics},
    volume={604},
    publisher={Springer},
    address={New York-Berlin},
    date={1977}
    },
  pages={125--137}
  }

\bib{smith1}{article}{
  title={Products of Banach spaces that are uniformly rotund in every direction},
  author={Smith, M.A.},
  journal={Pacific J. Math.},
  volume={73},
  number={1},
  date={1977},
  pages={215--219},
  review={\mr{0463892}}
  }

\bib{smith}{article}{
  title={Some Examples Concerning Rotundity in Banach Spaces},
  author={Smith, M.A.},
  journal={Math. Ann.},
  volume={233},
  date={1978},
  pages={155--161}
  }

\bib{smith2}{article}{
  title={A curious generalization of local uniform rotundity},
  author={Smith, M.A.},
  journal={Comment. Math. Univ. Carolinae},
  volume={25},
  number={4},
  date={1984},
  pages={659--665},
  review={\mr{0782015}}
 }

\bib{sullivan}{article}{
  title={Geometrical properties determined by the higher duals of a Banach},
  author={Sullivan, F.},
  journal={Illinois J. Math.},
  volume={21},
  number={2},
  date={1977},
  pages={315--331}
  }

\bib{yorke}{article}{
  title={Weak rotundity in Banach spaces},
  author={Yorke, A.C.},
  journal={J. Austral. Math. Soc.},
  volume={24},
  date={1977},
  pages={224--233},
  review={\mr{0482084}}
  }
  

\bib{zizler}{article}{
  title={Some notes on various rotundity and smoothness properties of separable Banach spaces},
  author={Zizler, V.},
  journal={Comment. Math. Univ. Carolin.},
  volume={10},
  number={2},
  date={1969},
  pages={195--206},
  review={\mr{0246095}}
  }
  
\end{biblist}
\end{bibdiv}
\end{document}